\documentclass[11pt]{amsart}

\usepackage{amssymb,amsthm,amsmath, graphicx, esvect, xcolor, tikz}
\usepackage{hyperref}
\usepackage[alphabetic]{amsrefs}
\usepackage[shortlabels]{enumitem}
\usepackage[margin=3cm]{geometry}
\usepackage{standalone}
\usepackage{microtype}

\usetikzlibrary{intersections, calc}

\usepackage[T1]{fontenc}

\newtheorem{theorem}{Theorem}[section]
\newtheorem{corollary}[theorem]{Corollary}
\newtheorem{lemma}[theorem]{Lemma}
\newtheorem{proposition}[theorem]{Proposition}
\newtheorem{definition}[theorem]{Definition}

\theoremstyle{remark}
\newtheorem{remark}[theorem]{Remark}

\newtheorem*{AttnProblem}{Attainability Problem}
\newtheorem*{BangBangProblem}{Bang-Bang Problem}

\numberwithin{equation}{section}

\newcommand{\N}{\mathbb{N}}

\newcommand{\R}{\mathbb{R}}
\newcommand{\C}{\mathcal{C}}

\newcommand{\area}{\mathrm{area}}
\newcommand{\co}{\mathrm{co}}
\newcommand{\PP}{\mathcal{P}}

\newcommand{\D}{\mathcal{D}}
\newcommand{\TT}{\mathcal{T}}
\newcommand{\V}{\mathcal{V}}

\newcommand{\arc}{\mathrm{arc}}
\newcommand{\cw}{\mathrm{cw}}
\renewcommand{\precsim}{\preccurlyeq}
\renewcommand{\epsilon}{\varepsilon}
\renewcommand{\leq}{\leqslant}
\renewcommand{\geq}{\geqslant}
\renewcommand{\le}{\leqslant}
\renewcommand{\ge}{\geqslant}

\title{Decreasing paths of polygons}

\author[I. Kulp]{Isaac Kulp}
\address[Isaac Kulp]{Department of Mathematics, Baylor University, 1410 S. 4th Street, Waco, TX 76706}
\email{isaac\_kulp1@baylor.edu}

\author[C. Ochanine]{Charlotte Ochanine}
\address[Charlotte Ochanine]{
Department of Mathematics,
University of Louisiana at Lafayette,
217 Maxim Doucet Hall,
Lafayette, LA 70503 USA}
\email{charlotte.ochanine@louisiana.edu}

\author[L. Richard]{Logan Richard}
\address[Logan Richard]{
Department of Mathematics,
Oregon State University,
Kidder Hall 368,
Corvallis, Oregon 97331-4605}
\email{richalog@oregonstate.edu}

\author[L. Robert]{Leonel Robert}
\address[Leonel Robert]{
Department of Mathematics,
University of Louisiana at Lafayette,
217 Maxim Doucet Hall,
Lafayette, LA 70503 USA}
\email{lrobert@louisiana.edu}

\author[S. Whitman]{Scott Whitman}
\address[Scott Whitman]{
Department of Mathematics,
University of Louisiana at Lafayette,
217 Maxim Doucet Hall,
Lafayette, LA 70503 USA}
\email{scott.whitman1@louisiana.edu}

\begin{document}

\begin{abstract}
We call a continuous path  of polygons decreasing if the convex hulls of the polygons form a decreasing family of sets. For an arbitrary  polygon of more than three vertices, we characterize the polygons contained in it that can be reached by a decreasing path (attainability problem), and we show that this can be done by a finite application of ``pull-in" moves (bang-bang problem).  
In the case of triangles, these problems were investigated by Goodman, Johansen, Ramsey, and Frydman among others, in connection with the embeddability 
problem for non-homogeneous Markov processes.
\end{abstract}


\maketitle
\setcounter{tocdepth}{1}
\tableofcontents

\section{Introduction}

Given $n\geq 3$, by a polygon of $n$ vertices, or $n$-gon, we understand  an $n$-tuple $P=(p_i)_{i=1}^n$ of points in the plane. We consider the  pre-order relation on polygons defined as follows: we say that $P'$ is contained in $P$, denoted by $P'\precsim P$,  if the convex hull of the vertices of $P'$  is contained in the convex hull of the vertices of $P$. We call a continuous path $\{P(t): 0\leq t\leq c\}$ of $n$-gons a decreasing path if it is decreasing with respect to this pre-order. We investigate in this paper the attainability problem for decreasing paths of polygons:

\begin{AttnProblem}
\label{attainabilityproblem}
Given $n\geq 3$ and an $n$-gon $P$, describe the set of $n$-gons $P'$ contained in $P$  and attainable from $P$ by a decreasing path, i.e., such that there exists a decreasing path  of $n$-gons $\{P(t): 0\leq t\leq c\}$ satisfying $P=P(0)$ and $P'=P(c)$.
\end{AttnProblem}

\begin{figure}[h]
\centering
\begin{tikzpicture}[line width=1pt, x=1.0cm, y=1.0cm, scale=0.9]

\coordinate (p1) at (0.48,-0.26);
\coordinate (p2) at (6.35,-0.27);
\coordinate (p3) at (4.5,4.8);
\coordinate (p4) at (1.15,3.7);

\coordinate (p1') at (1.47, 0.47);
\coordinate (p2') at (5.20, 0.16);
\coordinate (p3') at (4.18, 3.47);
\coordinate (p4') at (2.15, 3.40);

\draw (p1) -- (p2) -- (p3) -- (p4) -- cycle;

\draw (p1') -- (p2') -- (p3') -- (p4') -- cycle;

\draw[line width=0.5pt]  (p1) to [bend left = 10] coordinate[pos=0.33] (s1) coordinate[pos=0.66] (t1) (p1');
\draw[line width=0.5pt]  (p2) to [bend left=8] coordinate[pos=0.33] (s2) coordinate[pos=0.66] (t2) (p2');
\draw[line width=0.5pt]  (p3) to [bend left = 10] coordinate[pos=0.33] (s3) coordinate[pos=0.66] (t3) (p3');
\draw[line width=0.5pt]  (p4) to [bend left =10] coordinate[pos=0.33] (s4) coordinate[pos=0.66] (t4) (p4');

\draw[dashed] (s1) -- (s2) -- (s3) -- (s4) -- cycle;
\draw[dashed] (t1) -- (t2) -- (t3) -- (t4) -- cycle;

\foreach \point in {p1,p2,p3,p4,p1',p2',p3',p4'} {
    \draw[fill=black] (\point) circle (1pt);

\draw (p1) node[below] {$p_1$};
\draw (p2) node[right] {$p_2$};
\draw (p3) node[above] {$p_3$};
\draw (p4) node[left] {$p_4$};
\draw (p1') node[anchor = south west, inner sep = 2pt] {$p_1'$};
\draw (p2') node[above left] {$p_2'$};
\draw (p3') node[anchor=north east, inner sep = 1pt] {$p_3'$};
\draw (p4') node[anchor=north west, inner sep = 1pt] {$p_4'$};
}

\end{tikzpicture}

\caption{Decreasing path of quadrilaterals.}
\label{fig_path}
\end{figure}
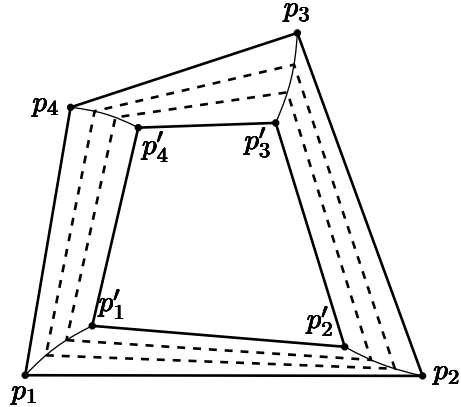

The special case of this problem concerning triangles ($n=3$) has previously been investigated in the context of the embeddability problem for stochastic matrices. The connection between these two problems goes roughly as follows (we shall discuss this in more detail in Section \ref{markovsec}): Let $[0,t_0]\ni t\mapsto (a_t,b_t,c_t)$  be  a decreasing path of non-degenerate triangles. For  $0\leq t\leq t'\leq t_0$ we obtain a
 3x3 row-stochastic matrix $E(t,t')$ by expressing the vertices of $(a_{t'},b_{t'},c_{t'})$ as convex combinations of the vertices of $(a_{t},b_{t},c_{t})$:
\[
\begin{pmatrix}
a_{t'}\\b_{t'}\\c_{t'}
\end{pmatrix}=E(t,t')
\begin{pmatrix}
a_{t}\\b_{t}\\c_{t}
\end{pmatrix}\hbox{ for }0\leq t\leq t'\leq t_0.
\]
These stochastic matrices form the transition matrices of 
a continuous-time non-homogeneous Markov process $\{D(t,t'):0\leq t\leq t'\leq t_0\}$ defined as 
\[
D(t,t')=E(t_0-t',t_0-t).
\] 
Conversely,  from a non-homogeneous Markov process of 3 states $\{D(t,t'):0\leq t\leq t'\leq t_0\}$ and a non-degenerate triangle $(a_0,b_0,c_0)$,
we obtain a decreasing path of triangles starting at $(a_0,b_0,c_0)$.
The embeddability problem seeks to determine whether a given row stochastic matrix $D$ can be obtained as a transition matrix in a (non-homogeneous) Markov process $\{D(t,t'):0\leq t\leq t'\leq t_0\}$ (\cite{goodmanCIME}).  
The geometric reformulation of a Markov process of 3 states as a decreasing path of triangles precisely leads to the attainability problem for triangles. Johansen and Ramsey made use of this reformulation to investigate the embeddability problem for Markov processes of 3 states (\cite{johansen-ramsey}).

We obtain a complete solution of the Attainability Problem for $n\geq 4$. The problem remains open for $n\geq 3$.

One way to create decreasing paths of polygons is to use pull-in moves. Given a polygon $P=(p_i)_{i=1}^n$ and two distinct indices  $i$ and $j$, a pull-in of vertex $p_i$ toward $p_j$ results in a polygon $P'$ whose vertices agree with the vertices of $P$  except for vertex $p_i'$,  which is chosen anywhere on the segment $[p_i,p_j]$.  Clearly, if $P'$ is obtained from $P$ by a pull-in move, there is a decreasing path going from $P$ to $P'$. Thus, if $P'$ can be obtained from $P$ by a finite application of pull-in moves, then $P'$ is attainable from $P$.

\begin{figure}[h]
\centering
\begin{tikzpicture}[line width = 1pt, x=1cm, y=1cm, scale=0.8]

 \coordinate (p1) at (0.,0.);
    \coordinate (p2) at (4.,0.);
    \coordinate (p3) at (6.3,1.95);
    \coordinate (p4) at (3.54,4.50);
    \coordinate (p5) at (0.,2.97);
\coordinate (p3') at (4.97,2.17);

    \draw (p1) -- (p2) -- (p3) -- (p4) -- (p5) -- cycle;

    \draw [dashed] (p3) -- (p5);
    \draw (p2) -- (p3');
    \draw (p3') -- (p4);

    \node[anchor=north west] at (p3) {$p_3$};
    \node[anchor=north east] at (p5) {$p_5$};
    \node[anchor=north west] at (p3') {$p_3'$};

   \foreach \point in {p1, p2, p3, p4, p5, p3'} {
        \draw [fill=black] (\point) circle (1pt);
    }
\end{tikzpicture}
\caption{Pull-in of $p_3$ toward $p_5$.}
\label{fig_pullin}
\end{figure}
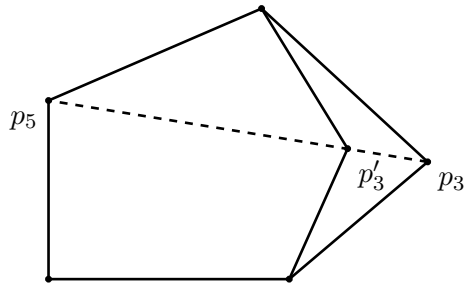

\begin{BangBangProblem}\label{movesproblem}
Is every attainable polygon $P'$ contained in $P$ attainable in finitely many pull-in moves? If so, is there an upper bound on the number of moves needed? 
\end{BangBangProblem}

Johansen and Ramsey showed in \cite{johansen-ramsey} that an attainable triangle is also attainable by finitely many pull-in moves, but did not obtain an upper bound on the number of moves needed. An example by Frydman shows that this number is at least $7$ (\cite{frydman83}). 
We solve here the Bang-Bang Problem affirmatively for $n\geq 4$.

Before stating our solutions of the Attainability and Bang-Bang Problems for $n\geq 4$, we need to introduce the concept of degenerate containment: Given $n$-gons $P$ and $P'$, we say that $P'$ is degenerately contained in $P$, or simply degenerate in $P$, if there exists an $m$-gon $Q$ with $m<n$ such that $P'\precsim Q\precsim P$.  For $n=3$, degenerate containment  entails that the vertices of $P'$ are collinear. This is not the case for $n\geq 4$, and degeneracy becomes a milder restriction as $n$ increases.

\begin{figure}[h]
\centering

\begin{tikzpicture}[line width= 1pt, x=1.0cm, y=1.0cm]

  \coordinate (p1) at (0,0);
  \coordinate (p2) at (3.88,0);
  \coordinate (p3) at (5.05,2.16);
  \coordinate (p4) at (2.17,4.31);
  \coordinate (p5) at (-0.64,3.40);

  \coordinate (p1') at (2.8,0.68);
  \coordinate (p2') at (3.90,1.3);
  \coordinate (p3') at (2.65,3.31);
  \coordinate (p4') at (0.93,2.86);
  \coordinate (p5') at (0.39,0.83);

  \coordinate (s1) at (-0.09,0.47);
  \coordinate (s2) at (4.10,0.41);
  \coordinate (s3) at (3.85,3.06);
  \coordinate (s4) at (1.01,3.93);

  \fill[color = black!10!white] (s1) -- (s2) -- (s3) -- (s4) -- cycle;
   \draw[dashed] (s1) -- (s2) -- (s3) -- (s4) -- cycle;

  \draw(p1') -- (p2') -- (p3') -- (p4') -- (p5') -- cycle;

  \draw (p1) -- (p2) -- (p3) -- (p4) -- (p5) -- cycle;

  \foreach \p in {p1, p2, p3, p4, p5, p1', p2', p3', p4', p5', s1, s2, s3, s4}
    \draw[fill=black] (\p) circle (1pt);    
    
    
    \node[below left] at (p1) {$p_1$};    
    \node[below right] at (p2) {$p_2$};
    \node[right] at (p3) {$p_3$};
    \node[above] at (p4) {$p_4$};
    \node[below left] at (p5) {$p_5$};    
    \node[above] at (p1') {$p_1'$};
    \node[left, yshift=4pt] at (p2') {$p_2'$};
    \node[below, xshift=-2pt] at (p3') {$p_3'$};
    \node[below right] at (p4') {$p_4'$};    
    \node[above right] at (p5') {$p_5'$};
    
  \end{tikzpicture}
\caption{Pentagon $P'$ degenerately contained in pentagon $P$.}
\label{fig_degenerate}
\end{figure}
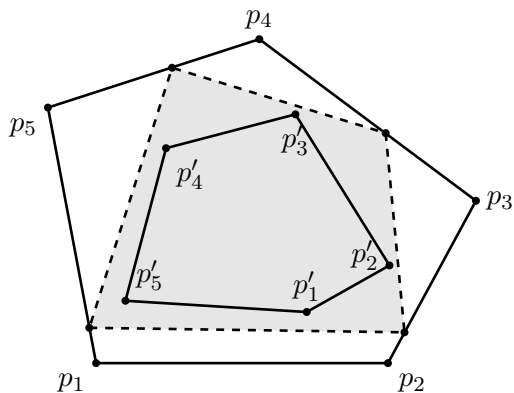

\begin{theorem}[Bang-Bang Theorem]\label{mainmoves}
Let $P$ be an $n$-gon. The following are true:
\begin{enumerate}[(i)]
\item
Any $n$-gon degenerately contained in $P$ is attainable from $P$ in fewer than $5n$ pull-in moves. 
\item
For $n\geq 4$, any $n$-gon non-degenerately contained in $P$ and attainable from $P$ is attainable in at most $2n$ pull-in moves.
\end{enumerate}
\end{theorem}
The number $5n$ in part (i) is unlikely to be optimal. For example, for $n=3$ degenerately contained triangles can always be attained in $5$ pull-in moves.  

We  now proceed to describe our characterization of attainable polygons. Let us first dispose of an easy case: if the $n$-gon $P$ is non-convex, then any $n$-gon $P'$  contained in $P$ is degenerately contained in $P$, and thus attainable from $P$ by Theorem \ref{mainmoves}. Assume thus that $P$ is convex and, after a re-indexing of its vertices, oriented counterclockwise. Fix $P'$ contained in $P$.  We consider a function $\pi\colon \partial P\to \partial P$, on the boundary of $P$, which we call the Poncelet map. Roughly described, for $x$ in the boundary of $P$, the point $\pi(x)$ is obtained by issuing a right tangent ray to $P'$ from $x$ and letting $\pi(x)$ be the new intersection of this ray with the boundary of $P$ (see Definition \ref{ponceletdef}). The \emph{broken line construction} with starting point $x$ consists in finitely many iterations of the Poncelet map 
 $x,\pi(x),\pi^2(x),\ldots$ allowed to go just once around $\partial P$ (Definition \ref{BLC}).

\begin{figure}[h]
\centering

\begin{tikzpicture}[line width=1pt, x=1.0cm, y=1.0cm]
    \coordinate (P1) at (0.,0.);
    \coordinate (P2) at (2.65,0.);
    \coordinate (P3) at (4.18,1.88);
    \coordinate (P4) at (2.5,4.34);
    \coordinate (P5) at (-0.18,3.97);
    \coordinate (P6) at (-1.21,2.46);

    \coordinate (Q1) at (0.06,0.76);
    \coordinate (Q2) at (1.66,0.44);
    \coordinate (Q3) at (3.09,1.65);
    \coordinate (Q4) at (2.42,3.4);
    \coordinate (Q5) at (0.30,3.09);
    \coordinate (Q6) at (-0.24,1.95);

    \draw (P1) -- (P2) -- (P3) -- (P4) -- (P5) -- (P6) -- cycle;
    \draw (Q1) -- (Q2) -- (Q3) -- (Q4) -- (Q5) -- (Q6) -- cycle;

    \coordinate (D1) at (1.94,0.);
    \coordinate (D2) at (3.72,2.55);
    \coordinate (D3) at (1.25,4.17);
    \coordinate (D4) at (-0.86,1.76);
    \coordinate (D5) at (0.76,0.);

    \draw[dashed] (D1) -- (D2);
    \draw[dashed] (D2) -- (D3);
    \draw[dashed] (D3) -- (D4);
    \draw[dashed] (D4) -- (D5);

    \foreach \point in {P1, P2, P3, P4, P5, P6, Q1, Q2, Q3, Q4, Q5, Q6, D1, D2, D3, D4, D5} {
        \draw [fill=black] (\point) circle (1pt);
    }

\draw (D1) node[below] {$x$};
\draw (D2) node[above right] {$\pi(x)$};
\draw (D3) node[above] {$\pi^2(x)$};
\draw (D4) node[below left] {$\pi^3(x)$};
\draw (D5) node[below] {$\pi^4(x)$};
\end{tikzpicture}

\caption{Broken line construction along $\partial P$ starting at $x$.}
\label{fig_poncelet}
\end{figure}
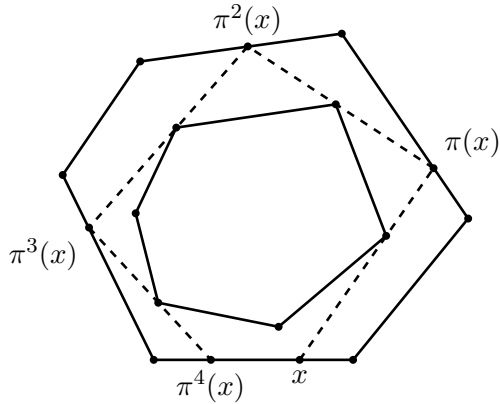

Equipped with the Poncelet map and the broken line construction (BLC), we obtain an effective way of testing for degenerate containment of $P'$ in $P$ and we formulate a description of the set of polygons attainable from $P$. We describe these results here in general terms.  

To test for the degenerate containment  of $P'$ in $P$ we ask that the BLC stops before $n$ steps when its initial  point $x$ is chosen from a finite set of test points (Theorem \ref{degeneracytestthm}). 

In order to describe the set of attainable polygons we first need to introduce two of its subsets:  the threshold and the vestibule. A polygon $P'$ contained in $P$ is said to be in the {\it threshold region}  if it is attainable from $P$ and at least one of its vertices belongs to $\partial P$. In Theorem \ref{thresholdthm} we characterize the polygons in the threshold and non-degenerately contained in $P$ in terms of the BLC. We further show that such polygons  can be attained from $P$ in at most $2n-1$ pull-in moves. (For triangles, this result was essentially obtained by Frydman in \cite[Theorem 3.1]{frydman80}, albeit not stated in this way). 
A polygon $P'$ is said to be in the {\it vestibule region} if it is an  attainable polygon  obtained from a polygon in the threshold by one pull-in move,  i.e., it is one pull-in move away from a polygon attainable from $P$ and with at least one vertex in $\partial P$. 

\begin{theorem}[Attainability Theorem]\label{mainattainable}
Let $n\geq 4$ and let $P$ be a convex $n$-gon. The set of polygons attainable from $P$ coincides with the set $\D_P\cup\V_P$ of $n$-gons that are either degenerately contained in $P$ or in the vestibule region.
\end{theorem}

This result,  combined with effective ways of testing whether a polygon is degenerately contained in $P$ and in the vestibule region (Theorems \ref{degeneracytestthm} and \ref{thresholdthm}), solves the Attainability Problem for $n\geq 4$. 
In contrast, Frydman's example of a triangle attainable in 7 moves and no less  does  not belong to $\D_P\cup\V_P$. Thus, the conclusion of Theorem \ref{mainattainable} fails for $n=3$.

\textbf{Acknowledgments}: We are grateful to George Turcu for useful discussions around the attainability problem for triangles
and the work of Johansen and Ramsey.

\section{Decreasing paths of polygons}\label{sec:decreasingpaths}
Let  $n\ge 3$.  By an $n$-gon, or polygon, we mean an ordered $n$-tuple $(p_1,\ldots,p_n)$ of points in $\R^2$. The points $p_1,\ldots,p_n$ are called vertices of the polygon. 
We denote the set of $n$-gons by  $\PP_n$ (thus, $\PP_n=(\R^2)^n$). Note that in the sense that we use the term polygon here a vertex may appear multiple times in a polygon, and that the indexing of the vertices in a polygon matters, i.e., a permutation of the vertices of a polygon may, in general, result in a different polygon. 

We will find it convenient to use arithmetic modulo $n$ in the index set of an $n$-gon, so that $p_{n+1}=p_1$ and  the list $i,i+1,\ldots,i+n-1$ runs through all indices for any $i$.

Given points $a,b\in \R^2$, we denote by $[a,b]$ the segment connecting them and including them as endpoints. We follow standard convention for half-open and open segments, e.g., $[a,b)$ denotes $[a,b]\backslash \{b\}$.

Let  $P=(p_i)_{i=1}^n$ be an $n$-gon.
\begin{itemize}
\item
We denote by $\{P\}$ the set of vertices of $P$, i.e., the set $\{p_1,\ldots,p_n\}$.

\item
 We call the segments $[p_1,p_2],[p_2,p_3],\ldots, [p_n,p_1]$ the edges of $P$. We denote by $\partial P\subseteq \R^2$ the union of the edges of $P$.

\item 
We denote by $\co(P)\subseteq \R^2$ the convex hull of the vertices of $P$. 

\item 
 We call $P$ \emph{set-convex}  if none of its vertices is a convex combination of the other ones. Equivalently, if the vertices of $P$ are pairwise distinct 
 and are all extreme points of $\co(P)$. 
 
 \item
 We call $P$ a \emph{simple polygon} if  $\partial P$ forms a simple closed curve. In this case, we say that $P$ is \emph{oriented counterclockwise}, if its vertices
 $p_1,p_2,\ldots,p_n$ follow the counterclockwise order of $\partial P$, and it is \emph{oriented clockwise} if its vertices follow the clockwise order
 of $\partial P$.
 
 \item 
 We call $P$ a \emph{convex polygon} if it is both set-convex and simple. Note that this in particular means that the vertices of $P$ are pairwise distinct.
 \end{itemize}

We define on the set of polygons the following preorder relation:

\begin{definition}
Given polygons $P$ and $Q$ (possibly with different numbers of vertices),
we say  that $P$ is contained in $Q$, denoted by $P\precsim Q$, if $\co(P)\subseteq \co(Q)$.
\end{definition}

Let $\{P(t):0\leq t\leq c\}$ be a continuous path in $\PP_n$. (The topology in $\PP_n=(\R^2)^n$ is always assumed to be the standard one.)
We call $P(\cdot)$ a \emph{decreasing path} if it is decreasing in the 
preorder $\precsim$, i.e., $P(t')\precsim P(t)$ for all $t'\geq t$. 

Given $n$-gons $P$ and $P'$, with $P'\precsim P$, we say that $P'$
is \emph{attainable} from $P$ if there exists a decreasing path in $\PP_n$ starting at $P$ and ending at $P'$.

Let $P=(p_i)_{i=1}^n$ be an $n$-gon. Let $1\leq i,j\leq n$ be distinct indices. A \emph{pull-in move} on $P$ of vertex $p_i$ toward vertex $p_j$ consists in replacing $p_i$ by a point $p_i'$ in the segment $[p_i,p_j]$, while leaving all other vertices of $P$ unchanged.  The result of a pull-in move is an $n$-gon $P'=(p_1,\ldots,p_i',\ldots,p_n)$, where  $p_i'=(1-c)p_i+cp_j$ for some  $0\leq c\leq 1$. We call $c$ the parameter of the move. Observe that $P'$ is attainable from $P$ by a decreasing path:
\[
[0,c]\ni t\mapsto (p_1,\ldots,p_{i-1},(1-t)p_i+tp_j,p_{i+1},\ldots,p_n).
\]
Thus, if a polygon is attainable from $P$ by the application of finitely many pull-in moves, then it is also attainable from $P$ by a decreasing path. 

\begin{remark}(On re-indexing.)\label{onreindexing} Let $P$ and $P'$ be $n$-gons, with $P'$ contained in $P$.
Let $\sigma$ be a permutation of the set of  indices $\{1,\ldots,n\}$, and set $P_{\sigma}=(p_{\sigma(i)})_{i=1}^n$, $P'_\sigma=(p'_{\sigma(i)})_{i=1}^n$.
We readily deduce from the definitions of attainability, and attainability by pull-in moves, that $P'$ is attainable from $P$ 
if and only if $P'_\sigma$ is attainable from $P_{\sigma}$. We can use this observation
to arrange for specific properties of either $P'$ or $P$, when examining the attainability question. For example, if $P$ is set-convex, then  $P_\sigma$
is convex and oriented counterclockwise for  a suitable permutation $\sigma$. In this case, the question of attainability of a given $P'$ from $P$ can be translated 
into the attainability of $P_\sigma'$ from the convex oriented counterclockwise polygon $P_\sigma$. 
\end{remark}

The inverse of a pull-in move is called a push-out move. More concretely, a push-out of vertex $p_i$ by vertex $p_j$ 
results in a polygon $P'=(p_1,\ldots,p_i',\ldots,p_n)$, where $p_i\in [p_j,p_i']$. Notice that in the special case when $p_i=p_j$ the new vertex
$p_i'$ can be chosen arbitrarily in $\R^2$.

Let us regard an $n$-gon as an $n\times 2$ matrix, the rows of the matrix simply being the vertices of the polygon. 
Let $P'$ and $P$ be $n$-gons such that $P'\precsim P$. Then we can write each vertex of $P'$ as a convex (not unique for $n\geq 4)$) combination of vertices of $P$. This can be expressed algebraically as
\[
P'=DP,
\]
where $D$ is an $n\times n$ (row) stochastic matrix and $P'$ and $P$ are regarded as $n\times 2$ matrices. The converse is also clearly true: if $P'=DP$, with $D$ stochastic, then $P'\precsim P$. 

Suppose that $P'$ is obtained from $P$ by a pull-in move, say, by the move of $p_i$ towards $p_j$ with parameter $c$. This can be expressed algebraically  as $P'=K^{ij}(c)P$, where $K^{ij}(c)$ is the matrix differing from the $n\times n$ identity matrix only in two entries of the $i$-th row: $k_{ii}=1-c$ and $k_{ij}=c$. We call any such  matrix an \emph{elementary stochastic matrix}. 
It is clear that $P'$ is attainable from $P$ in finitely many pull-in moves if and only if $P'=DP$ for $D$ expressible as a finite product of elementary stochastic matrices.

\begin{lemma}\label{commutingmoves}
A pull-in move of $p_i$ toward $p_j$ followed by a pull-in move of $p_k$ toward $p_l$ can always be expressed as 
a pull-in move of $p_k$ toward $p_l$ followed by a pull-in move of $p_i$ toward $p_j$, excepting the cases when
 $j=k$ and $i\neq l$, or  $i=l$ and  $j\neq k$. 
\end{lemma} 
\begin{proof}
Expressing pull-in moves in terms of multiplication by elementary stochastic matrices, this boils down to showing that given $c$ and $d$, there exist $c'$ and $d'$ such that
\[
K^{kl}(d)K^{ij}(c)=K^{ij}(d')K^{kl}(c').
\]

 If $i,j,k,l$ are pairwise distinct, then $K^{kl}(d)K^{ij}(c)$ is the matrix that agrees with the identity matrix except for the $ii$, $ij$, $kk$, $kl$-entries which are  $1-c,c,1-d, d,$ respectively. So $K^{ij}(c)$ and $K^{kl}(d)$ commute.

 If $i=k$ and $j=l$, then $K^{kl}(d)K^{ij}(c)=K^{ij}(c+d-cd).$ So we again have commutativity of $K^{ij}(c)$ with $K^{kl}(d)$.

  Sometimes $c',d'\neq d,c$, for example when $i=l$ and $j=k$. The full case when $i,j,k,l$ are not pairwise distinct (i,e., 
$p_i,p_j,p_k,p_l$ are vertices of a triangle) is handled in \cite[Lemma 3.2]{frydman2}. 
\end{proof}

\section{Degenerate polygons}\label{secdegenerate}

\begin{definition}
Given $n$-gons $P$ and $P'$,
we say that $P'$ is degenerately contained in $P$, or simply degenerate in $P$, if there exists an $m$-gon $Q$, with $m<n$, 
such that  $P'\precsim Q\precsim P$. In other words, there is a polygon $Q$ with fewer than $n$ vertices such that 
the convex hull of $P'$ is contained in that of $Q$, which in turn is contained in the convex hull of $P$.
\end{definition}

 Clearly, if $P''\precsim P'$, and $P'$ is degenerate in $P$, then $P''$ is also degenerate in $P$. Observe also that if $P'$ is not set-convex and contained in $P$, then it is degenerate in $P$, as in this case we can choose $Q$ to be any polygon formed by the vertices of the convex hull of $P'$. In particular, if $P$ itself is not set-convex, then $P$ is degenerate in $P$, and consequently, any $n$-gon $P'$ contained in $P$ is degenerate in $P$.

For a fixed $n$-gon $P$ we denote by $\mathcal D_P$ the set of all $n$-gons degenerately contained in $P$.
\begin{lemma}\label{compactD}
The set $\D_P$ is compact.
\end{lemma}
\begin{proof}
We can reformulate the existence of an $m$-gon $Q$ such that $P'\precsim Q\precsim P$ as the existence of row stochastic matrices
$D$ and $E$, of sizes $n\times m$ and $m\times n$ respectively, such that $P'=DEP$. 

Let $(P_k')_{k=1}^\infty$ be a sequence in  $\D_P$. Let us show the existence of a subsequence converging to an element of $\D_P$ (thus proving compactness). Clearly we can pass to a subsequence, which we relabel as $(P_k')_{k=1}^\infty$,  such that for a fixed $m<n$ and all $k=1,2\ldots$ we have $P_k'=D_kE_kP$ for row stochastic matrices $D_k$ and $E_k$ of sizes $n\times m$ and $m\times n$, respectively. Using the compactness of the sets of row stochastic matrices of sizes  $n\times m$ and $m\times n$, we can assume, after  passing again to subsequences and relabeling, that $D_k\to D$ and $E_k\to E$ entrywise. Then $P_k'\to P'$, with $P'=DEP$, which in turn implies that $P'\in \D_P$. 
\end{proof}

Our goal for the remainder of this section is to  show that if $P'$ is degenerate in $P$, then it is attainable from $P$ in fewer than $5n$ pull-in moves (Theorem \ref{in5nmoves}). Before doing this, we prove several lemmas.

Let us introduce some terms that will be used below. Let $P$ and $P'$ be polygons. 
\begin{itemize}
\item 
We say that a vertex  $p_i$ occupies a point $x\in \R^2$ if it coincides with it, i.e., $p_i=x$.

\item 
We call a point $x\in \R^2$ a double point of $P$ if there exist $i\neq j$ such that $p_i$ and $p_j$ both occupy $x$, i.e., $p_i=x=p_j$.

\item 
We say that a vertex $p_i'$ of $P'$  is stray on $\partial P$ if $p_i' \in \partial P$ and $p_i'$ does not occupy a vertex of $P$. 

\item 
We say that  $p_i'$ is stranded on $\partial P$ if it is stray on   $\partial P$ and no other vertices of $P'$
share the same edge with $p_i'$.

\item
We say that $P'$ is inscribed in $P$ if all the vertices of $P'$ belong to $\partial P$, i.e., $\{P'\}\subseteq \partial P$. 
\end{itemize}

\begin{figure}[h]
\centering

\begin{tikzpicture}[line width=1pt, x=1.0cm, y=1.0cm]

    \coordinate (A) at (0.,0.);
    \coordinate (B) at (2.62,0.);
    \coordinate (C) at (2.67,2.97);
    \coordinate (D) at (0.44,3.85);
    \coordinate (E) at (-1.06,1.86);
  \coordinate (A') at (2.63,0.86);
    \coordinate (B') at (2.66,2.14);
    \coordinate (C') at (-0.25,2.93);

    \draw   (A) -- (B) -- (C) -- (D) -- (E) -- cycle;

    \foreach \point in {B, C, D, E, A', B', C'} {
        \draw [fill=black] (\point) circle (1pt);
    }
    \draw [fill=black] (A) circle (2pt);

    \node[below] at (A) {$p_1'=p_2'$};
    \node[right] at (A') {$p_3'$};
    \node[right] at (B') {$p_4'$};
    \node[above left] at (C') {$p_5'$};
    
\end{tikzpicture}

\caption{The vertex $p_1'$ ($=p_2'$) is a double point of $P'$; the vertices $p_3'$ and $p_4'$ are stray (but not stranded); the vertex $p_5'$ is stranded}
\label{fig_stratstranded}
\end{figure}
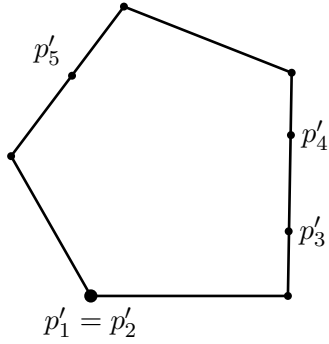

\begin{lemma}\label{fromP'toQ}
Let $P$ be a convex $m$-gon. Let $P'$ be an $n$-gon inscribed in $P$ with $n>m$. Then the vertices of $P'$ can be moved to the vertices of $P$ 
(not necessarily surjectively) in fewer than $\frac{3n}{2}$ push-out moves. That is, applying less than  $\frac{3n}{2}$ push-out moves to $P'$,
we can arrange that $\{P'\}\subseteq \{P\}$.
\end{lemma}


\begin{proof}
Assume that not all vertices of $P'$ occupy vertices of $P$, i.e., $\{P'\}\not\subseteq \{P\}$. To move vertices of $P'$ to vertices of $P$, we can proceed as follows. If a  stray vertex $p_i'$  shares some edge $[p_k,p_{k+1}]$ with another vertex $p_j'$, push $p_i'$ out with $p_j'$ onto either $p_k$ or $p_{k+1}$.
Continue moving any vertices of $P'$ that are stray but not stranded to vertices of $P$  until perhaps all remaining stray vertices of $P'$ are  stranded. In this situation, we are guaranteed the existence of a  vertex of $P$ that is a double point of $P'$.
Indeed, partitioning  $\partial P$ in half-open edges $\partial P= \bigsqcup_{l=1}^{m} [p_l,p_{l+1})$, we see that  at least one pair of vertices $p_i',p_j'$ belong to the same  half-open edge of $P$.  Since neither of them is stray, they occupy the same vertex of $P$.
We can use this double point to unstrand any stray vertex of $P'$. Say such a stray vertex exists and belongs to  some edge $[p_k,p_{k+1}]$. We first push $p_i'$ out onto $p_k$ with $p_j'$, then push the stray vertex onto $p_{k+1}$ with $p_i'$.  We can continue to move any stranded vertices using double points until perhaps we again have unstranded stray vertices. Alternating between single push-out moves and use of double points, we can move all vertices of  $P'$ to vertices of $P$.

Let's estimate the number of moves in this algorithm.
Let $M$ be the number of times single push-out moves are used to move unstranded stray vertices of $P'$, and let  $N$ be the number of times that a double point is used to move stranded  vertices of $P'$. The total number of push-out moves in the algorithm is $M+2N$.
 Let $s$ denote the initial number of stray vertices of $P'$, and let $u$ denote the initial number of unoccupied vertices of $P$.  Since at every step the number of stray vertices is lowered by 1, we have that $M+N\leq s$.  Since using a double point lowers the number of unoccupied vertices by 2, we have $2N\leq u$. We can now  estimate the total number number of push-out moves used:
\[
M+2N = (M+N)+N\leq s+ \lfloor \frac{u}{2}\rfloor\leq n+\lfloor \frac{m}{2}\rfloor<\frac{3n}{2}.\qedhere
\]
\end{proof}


\begin{lemma}\label{whenmeqn}
Let $P$ be a convex $n$-gon. Let $P'$ be an $n$-gon inscribed in $P$, non-degenerate in $P$, and such that at least one vertex of $P'$ occupies a vertex of $P$. Then the  vertices of $P'$ can be moved to occupy all the vertices of $P$ in at most $n-1$ push-out moves. That is, applying at most $n-1$ push-out moves to $P'$, we can arrange that 
$\{P'\}=\{P\}$.
\end{lemma}

\begin{proof}
As in the previous lemma, we search for stray vertices of $P'$ that are not stranded, and we move them to  the vertices of $P$ with one push-out move. Let us show that
this process can only stop when every vertex of $P'$ occupies a vertex of $P$. Suppose that we reach a polygon $P'$ such that each vertex of $P'$ is either a vertex of $P$ or is stranded on an edge of $P$. Then each half-open edge $[p_i,p_{i+1})$ of $P$ must contain at most one vertex of $P'$. For if two vertices of $P'$ occupy the same half-open edge, then 
either they both occupy the same vertex of $P$, contradicting that $P'$ is non-degenerate in $P$ (it interpolates between the original polygon and $P$), or one of the vertices is stray but not stranded. The half-open edges of $P$ are pairwise disjoint, by the convexity of $P$. Since the number of vertices of $P'$ agrees with the number of half-open edges of $P$, it follows that each half-open edge  contains exactly one vertex of $P'$. By assumption, at least one vertex $p_i$ of $P$ is occupied by a vertex of $P'$. Then $[p_{i-1},p_i)$ cannot contain a stranded vertex of $P'$. Thus, $p_{i-1}$ is occupied by a vertex of $P'$. Continuing in this way we see that all vertices of $P$ are occupied by vertices of $P'$, equivalently, every vertex of $P'$ occupies a vertex of $P$. 

At each step of the algorithm the number of stray vertices is decreased by 1. Thus, the total number of push-out moves used is at most $n-1$.
\end{proof}

\begin{lemma}\label{permutation}
Let $P$ and $P'$ both be $n$-gons such that every vertex of $P'$ occupies a vertex of $P$, i.e., $\{P'\}\subseteq \{P\}$, and $P$ has at least one double point. Then we can go from $P'$ to $P$ in no more than $\displaystyle\frac{3n}{2}$ push-out moves. That is, applying at most $\displaystyle\frac{3n}{2}$ push-out moves to $P'$, we can arrange that $P'=P$.
\end{lemma}

\begin{proof}
Let $P'=(p_i')_{i=1}^n$ and $P=(p_i)_{i=1}^n$.
We want an algorithm that, for each $i$, moves $p_i'$ to $p_i$ using push-out moves. We call a vertex $p_i'$ of $P'$ correctly placed if  $p_i'=p_i$, and we call it incorrectly placed otherwise.

Let us denote by $S$ be the set of vertices of $P$, i.e., $S=\{P\}$. By the assumption that $P$ has a double point, we have $|S|<n$. Since $P'$ is an $n$-gon, 
and $\{P'\}\subseteq S$, there exists at least one point in  $S$ occupied by multiple vertices of $P'$. Let $q\in S$ be occupied by multiple vertices of $P'$. 
If one of these vertices, say $p_k'$, is incorrectly placed,  we can move it to its correct destination $p_k$ in a single push-out move. Let us continue performing this kind of move until there are no double points among the incorrectly placed vertices of $P'$ (see Figure \ref{fig_permutation}). Let $d\geq 1$ be the total number of moves this requires.

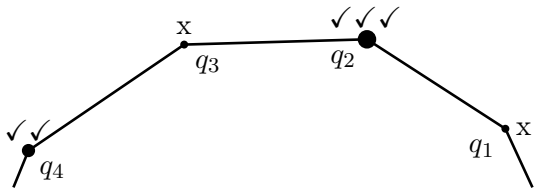
\begin{figure}[h]\label{fig_permutation}
\centering

\begin{tikzpicture}[line width=1pt,x=1.0cm,y=1.0cm]

\coordinate (p1) at (0.32,0.83);
\coordinate (p2) at (0.52,1.32);
\coordinate (p3) at (2.59,2.73);
\coordinate (p4) at (5.02,2.80);
\coordinate (p5) at (6.86,1.61);
\coordinate (p6) at (7.22,0.83);

\draw (p1) -- (p2) -- (p3) -- (p4) -- (p5) -- (p6);

    \draw (p2) node[above] {$\checkmark\checkmark$};
    \draw (5,3.10) node {$\checkmark\checkmark\checkmark$};
    \draw (p3) node[above] {x};
    \draw (p5) node[right] {x};

\draw [fill=black] (p2) circle (2pt) node[below right] at (p2) {$q_4$};
\draw [fill=black] (p3) circle (1pt) node[below right] at (p3) {$q_3$};
\draw [fill=black] (p4) circle (3pt) node[below left] at (p4){$q_2$};
\draw [fill=black] (p5) circle (1pt) node[below left] at (p5){$q_1$};

\end{tikzpicture}

\caption{$q_2,q_4$ are multiply occupied by correctly placed vertices of $P'$; $q_1,q_3$ are singly occupied by 
incorrectly placed vertices of $P'$.}
\label{fig_permutation}
\end{figure}

Let $T\subseteq S$ be the set of all $q\in S$ such that a single incorrectly placed vertex of $P'$ occupies $q$. 
We note, for later use, that  $T$ is a proper subset of $S$, since at least one point of $S$ is multiply occupied by vertices of $P'$, which are correctly placed at the current stage of the algorithm.
Define $\sigma\colon T \to T$ as follows: Given $q\in T$, find $p_k$ such that $p_k=q$.  Let $q'$ be the point of $S$ occupied by  $p_k'$. We cannot have $q'=q$, for if this were the case then $p_k'$ would be correctly placed at $q$, contradicting that $q\in T$.
Thus, $q'\neq q$. Since $p_k'$ is incorrectly placed (at $q'$), we have that $q'\in T$. Define $\sigma(q)=q'$. Note that this relation means that the vertex of $P'$ occupying $q'$ is correctly placed when moved to $q$.  Since there is only one correct placement for each vertex of $P'$, $\sigma$ is an injective map from $T$ to $T$ and therefore a permutation of $T$. 


Let $\sigma^{-1}=\tau_r\cdot\cdot\cdot\tau_2\tau_1$ be the cycle decomposition of $\sigma^{-1}$. Each cycle permutes some subset of the vertices of $P'$ to their correct places. Note that, if we had a double point to work with, this cyclic permutation could be thought of as sequence of push-out moves. So, let us  borrow a vertex  $p$ of $P'$ from some multiply occupied vertex of $S$ (the existence of one such vertex follows from $|S|<n$). We send $p$ to one of the vertices permuted by $\tau_1$. Then we correctly place the subset of vertices of $P'$ related to $\tau_1$ using a finite number of push-out moves that cycle the vertices forward. The number of push-out moves is exactly the length $l(\tau_1)$ of the cycle. Then we send $p$ to one of the vertices permuted by $\tau_2$ and cycle those vertices forward, and continue like this until we have cycled the vertices of $\tau_r$ forward. We then return $p$ to where it was originally as our final push-out move.  The number of push-out moves used in this operation is 
\begin{align*} 
1+l(\tau_1)+1+l(\tau_2)+...+1+l(\tau_r)+1 &= r+ (l(\tau_1)+...l(\tau_r))+1 \\
& = r+|T|+1.
\end{align*}

Recall that we denote by $d$ the number of moves used to correctly place all double points of $P'$. 
The total number of push-out moves used to go from $P'$ to $P$ is $N=d+r+|T|+1$. Let us estimate this number.
Since each cycle of $\sigma$ has length at least 2 ($\sigma$ has no fixed points), $r\leq \frac{|T|}{2}$. 
Thus, $N\leq d + \frac{3}{2}|T| + 1$.
Since, after the first $d$ moves, at least $d$ vertices of $P'$ are correctly placed,
and the number of of correctly placed vertices is $n-|T|$,  we have  that
$d\leq n-|T|$.
Thus,
\[
N\leq n-|T|+ \frac{3}{2}|T| + 1=n+ \frac{|T|}{2}+1\leq n+\frac{n-2}{2}+1=\frac{3n}{2}.
\]
Here we have used that  $|T|\leq n-2$, since $|T|< |S|<n$.
\end{proof}

Let $P$ be a convex $n$-gon. 
We call a polygon $Q$ contained in  $P$ a \emph{maximal degenerate} polygon if 
\begin{itemize}
\item
$Q$ is an $(n-1)$-gon inscribed in $P$,
\item
each vertex of $Q$ either singly occupies a vertex of $P$ (no other vertex of $Q$ occupies that same vertex of $P$) or lies stranded on an edge of $P$.
\end{itemize}

\begin{figure}[h]
\centering

\begin{tikzpicture}[line width=1pt, x=1.0cm, y=1.0cm]
 
    \coordinate (A) at (0.,0.);
    \coordinate (B) at (2.65,0.);
    \coordinate (C) at (2.84,3.37);
    \coordinate (D) at (0.49,3.92);
    \coordinate (E) at (-1.12,2.38);

    \draw (A) -- (B) -- (C) -- (D) -- (E) -- cycle;

    \coordinate (F) at (2.72,1.22);
    \coordinate (G) at (1.79,3.62);
    \coordinate (H) at (-0.52,2.96);

    \draw (A) -- (F) -- (G) -- (H) -- cycle;

    \foreach \point in {A, B, C, D, E, F, G, H} {
        \draw [fill=black] (\point) circle (1pt);
    }
\end{tikzpicture}

\caption{Maximal degenerate quadrilateral in a convex pentagon}
\label{fig_maxdeg}
\end{figure}
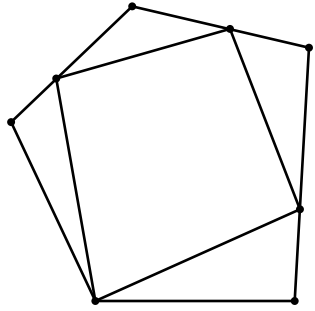

\begin{lemma}\label{maxdeg}  
Let $P$ be a convex $n$-gon.
\begin{enumerate}[(i)]
\item
For every $m$-gon $Q$ contained in $P$, with $m<n$, there exists a maximal degenerate  $(n-1)$-gon $Q'$ such that $Q\precsim Q'\precsim P$. 
\item
If $Q\precsim P$ is a maximal degenerate polygon then $Q$ has the following maximality property: if  a polygon $Q'$ is degenerately contained in $P$ and such that $Q\precsim Q'\precsim P$, then $\co(Q)=\co(Q')$.
\end{enumerate} 
\end{lemma}

\begin{proof}
(i) Let $Q$ be an  $m$-gon contained in $P$, with $m<n$. If $m<n-1$, we can add vertices 
to  $Q$  chosen from $\co(P)$ so  that we get an $(n-1)$-gon contained in $P$. Assume thus that $Q$ is an $(n-1)$-gon. Let us now apply a sequence of push-out
moves on $Q$ that result in a polygon maximal degenerate in $P$. First, using push-out moves, move all the vertices of $Q$ to $\partial P$ 
(choose a vertex $q$ of $Q$ to push the other vertices onto $\partial P$, then push $q$ out as well onto $\partial P$ using any of the other vertices). 
Next, whenever two vertices of $Q$ share an edge $[p_i, p_{i+1}]$ of $P$, use one to push the other out onto either $p_i$ or $p_{i+1}$.
If a double point of $Q$ occurs at a vertex of $P$, split the double point by sending one point to an unoccupied vertex of $P$, which must exist by the pigeonhole principle.
Continue in this way until all vertices of the resulting polygon either are single occupants of a vertex of $P$ or are stranded on edges of $P$. The resulting polygon  is maximal degenerate and interpolates between the original polygon and $P$.

(ii) Let $Q$ be maximal degenerate in $P$, and suppose that $Q\precsim Q'\precsim P$ for some $Q'$ degenerate in $P$. 
Then there exists an $m$-gon $Q''$, with $m<n$, interpolating between $Q'$ and $P$. By (i), there  exists  $Q'''$ interpolating between $Q''$ and $P$ that is maximal degenerate. It will suffice to  show that $\co(Q)=\co(Q''')$.
Renaming $Q'''$ as $Q'$, let us show that if $Q\precsim Q'\precsim P$, with $Q$ and $Q'$ both maximal degenerate, then $\co(Q)=\co(Q')$. Since $Q$ and $Q'$ are both $(n-1)$-gons, and the vertices of $Q$ are pairwise distinct, it will suffice to show that every vertex of $Q$ is also a vertex of $Q'$.

Suppose, for the sake of contradiction, that $\{Q\}\backslash \{Q'\}\neq \varnothing$. Let $q\in \{Q\}\backslash \{Q'\}$.  Then $q$ cannot be a vertex of $P$,  since  any vertex of $Q$ that is a vertex of $P$ is an extreme point of $\co(P)$,  hence also an extreme point of $\co(Q')$, and hence a vertex of $Q'$. Thus, $q$ is stranded on some edge $[p_i,p_{i+1}]$ of $P$.  We have that 
\[
\co(Q')\cap [p_i,p_{i+1}]=\co(Q'\cap [p_i,p_{i+1}]),
\] 
since $[p_i,p_{i+1}]$ is a face of $\co(P)$. Thus, $q\in  \co(Q'\cap [p_i,p_{i+1}])$. But  $\{Q'\}\cap [p_i,p_{i+1}]$ is either a single stranded vertex of $Q'$, one of $p_i$ or $p_{i+1}$, or both of them. Since  $q$ is a stranded vertex of $Q$ that is not a vertex of $Q'$, the only possible case is the third: both $p_i$ and $p_{i+1}$ are vertices of $Q'$.  Note that neither $p_i$ nor $p_{i+1}$ can be vertices of $Q$. In summary,  if $q\in \{Q\}\backslash \{Q'\}$, then $q$ is stranded on an edge of $P$ and the vertices of $P$ adjacent to  $q$ belong to $\{Q'\}\backslash \{Q\}$.

We have that $0<|\{Q\}\backslash \{Q'\}|$, by assumption, and $|\{Q\}\backslash \{Q'\}|\leq n-1$, since $Q$ is an $(n-1)$-gon. Thus, the edges of $P$ on which points in 
$\{Q\}\backslash \{Q'\}$ are stranded form a non-empty proper collection of the edges of $P$.  The number of vertices of $P$ incident with these edges  is at least $|\{Q\}\backslash \{Q'\}|+1$. Since these vertices belong to  $\{Q'\}\backslash \{Q\}$, we obtain that $|\{Q'\}\backslash \{Q\}|>|\{Q\}\backslash \{Q'\}|$,   in contradiction 
with the fact that $|\{Q\}|=|\{Q'\}|=n-1$. It follows that  $\{Q\}\backslash \{Q'\}$ is empty, i.e., all vertices of $Q$ are vertices of $Q'$.
\end{proof}

\begin{theorem}\label{in5nmoves}
Let $P$ and $P'$ be $n$-gons, with $P'$ degenerately contained in $P$. Then we can go from $P'$ to $P$ in fewer than $5n$ push-out moves.
\end{theorem}

\begin{proof}
We deal separately with the cases that $P$ is set-convex and that it is not.

\textbf{Case that $P$ is not set-convex}. Define a convex $m$-gon $Q$ oriented counterclockwise created by suitably indexing the  vertices of $P$ which are extreme points of  $\co(P)$. 
Since $P$ is not set-convex, $m<n$.

We will move $P'$ to $P$ via push-outs in stages: first, using push-outs we arrange that $\{P\}\subseteq \partial Q$. With another sequence of push-outs we arrange that $\{P\}\subseteq \{Q\}$. Next,  we redistribute the vertices of $P'$ within the vertices of $Q$ in a suitable fashion. Finally we apply push-outs to obtain $P'=P$.

Choose some vertex $p_i'$ of $P'$ and use it to push out all other vertices of $P'$ onto edges of $Q$, then use another vertex of $P'$ to push $p_i'$ out onto an edge of $Q$. 
This requires a total of $n$ moves. Continue to denote by $P'$ the resulting polygon. At this point, we have arranged that $\{P'\}\subseteq \partial Q$.
Using Lemma \ref{fromP'toQ},   move the vertices of  $P'$ to the vertices of $Q$ in fewer than $ \frac{3n}{2}$ push-out moves. We have thus arranged that $\{P'\}\subseteq \{Q\}$.

 Define an $n$-gon $P''$ as follows: 
Choose a vertex $p_i$ of $P$ that is also a vertex of $Q$ (thus, an extreme point of $\co(P)$),  and for each 
$p_j$ of $P$ that is not a vertex of $Q$, move $p_j$ to $p_i$ in one pull-in move. Denote the resulting $n$-gon by $P''$.  
By Lemma \ref{permutation}, we can go from $P'$ to $P''$ in $3n/2$ push-out moves. 
By construction, we can go from $P''$ to $P$ in $n-m$ push-out moves (undoing the pull-in moves used to construct $P''$). In total,
we can go from the original $n$-gon $P'$ to $P$
in fewer than
\[
n+ \frac{3n}{2}+\frac{3n}{2}+(n-m)< 5n
\] 
push-out moves.

\textbf{Case that $P$ is set-convex}. Assume that $P$ is set-convex. After re-indexing its vertices, along with the vertices of $P'$, we may assume that $P$ is convex and oriented counterclockwise (see Remark \ref{onreindexing}; note that degenerate containment is preserved after re-indexing both polygons). 

Since $P'$ is degenerate in $P$, there exists an $m$-gon $Q$ interpolating between $P'$ and $P$, with $m<n$. By Lemma \ref{maxdeg} (i), we can choose $Q$ to be a  maximal degenerate  $(n-1)$-gon. 

We will move $P'$ to $P$ via push-outs going through the same stages as before: first, we arrange that $\{P\}\subseteq \partial Q$, then we  arrange that $\{P\}\subseteq \{Q\}$. Next,  we redistribute the vertices of $P'$ within the vertices of $Q$. Finally, we apply push-outs to obtain $P'=P$.

We can move the vertices of $P'$ to the edges of $Q$ in at most $n$ push-out moves: one vertex of $P'$ pushes all other vertices out, and then this vertex is pushed by another one.  Let us continue to denote by $P'$ the resulting polygon. At this point, we have arranged that $\{P'\}\subseteq \partial Q$.
By Lemma \ref{fromP'toQ}, we can then move the vertices of  $P'$ to the vertices of $Q$ in fewer than $\frac{3n}{2}$ push-out moves. Let us continue to denote by $P'$ the resulting polygon. We have thus arranged that $\{P'\}\subseteq \{Q\}$.

Define the $n$-gon $\tilde Q=(q_1,q_1,q_2,\ldots,q_{n-1})$, with a double point at $q_1$.

\emph{Claim}: The vertices of $\tilde Q$ can be moved to occupy every vertex of $P$ (i.e., $\{\tilde Q\}=\{P\}$) in at most $n$ push-out moves. Let us prove this.  Choose a vertex $p_i$ of $P$ that is not  a vertex of $Q$ (at least one such vertex exists since $n>n-1$).  In the first move, push $q_1$ out with $q_1$ onto $p_i$. Denote the resulting polygon by $\tilde Q'$. Since $\co(\tilde Q')$ is  strictly larger than $\co(Q)$, the $n$-gon  $\tilde Q'$ is non-degenerately contained in $P$, by the maximal degeneracy of $Q$ in $P$. By Lemma \ref{whenmeqn}, we can then move the vertices of $\tilde Q'$ to the vertices of $P$ in at most $n-1$ push-out moves. By the non-degeneracy of  $\tilde Q'$ in $P$, every vertex of $P$ is occupied by (exactly) one vertex of   $\tilde Q'$ in this process.

By the claim just established, there is an $n$-gon $P''$ attainable from $P$ by $n$ pull-in moves and agreeing with $\tilde Q$ up to a re-indexing of the vertices. Indeed, $P''$ is obtained starting from $P$ and undoing the push-out moves that take the vertices of $\tilde Q$ to occupy all vertices of $P$. 
Since the vertices of $P'$ at the current stage of the process occupy the vertices of $P''$, we can go $P'$ to $P''$ in no more than $\frac{3n}{2}$ push-out moves,
by Lemma \ref{permutation}. We then go from $P''$ to $P$ in $n$ push-out moves.  The total number of push-out moves taking us from the original $P'$ to $P$ is thus less than 
\[
n+\frac{3n}{2}+\frac{3n}{2}+n=5n.\qedhere
\]
\end{proof}

\section{Decreasing paths and Markov processes}\label{markovsec}
In this section we discuss the connection between decreasing paths of polygons and non-homogeneous Markov processes, and we prove a ``chattering principle" type of result (Corollary \ref{Pchattering}) that will be needed later on.

We start by reviewing the formulation of a continuous-times non-homogeneous Markov processes in $n$ states as a parametrized collection of row stochastic matrices. We refer the reader
to \cite{goodman} for further details.

By a continuous-time non-homogeneous Markov process of $n$ states we understand a continuous map $(s,t)\mapsto D(s,t)$ defined for $0\le s\le t<\infty$, taking values in the set of $n\times n$ non-singular row stochastic matrices, and such that
\begin{align}
\label{transitiveD}
D(r,t)&=D(r,s)D(s,t) \hbox{ for all $0\leq r\le s\le t$,}\\
D(t,t) &=I\hbox{ (the identity matrix) for all $t$}. 
\end{align}

If the function $D(t)=D(0,t)$ is Lipschitz, then it is in particular a.e. differentiable, and it satisfies that
\[
\frac{d}{dt} D(t)=D(t)Q(t)
\]
where $Q(t)$ is a measurable, a.e. bounded \emph{intensity matrix}, i.e., a matrix whose  rows add to 0 and whose off-diagonal entries are nonnegative. 
Conversely, if $t\mapsto Q(t)$ is  measurable, a.e. bounded on finite intervals, and  such that $Q(t)$ is an intensity matrix for all  $t\geq 0$, then the ODE above with $D(0)=I$ has a unique solution $D(t)$ ranging in the nonsingular stochastic matrices (explicitly given by a Peano-Baker series). Upon setting  $D(s,t)=D(s)^{-1}D(t)$ for $0\leq s\leq t$ we obtain a Markov process.

We call a non-singular row stochastic matrix $D$  embeddable if $D=D(s,t)$ for some $0\le s\le t$ and Markov process $\{D(s,t):0\leq s\leq t<\infty\}$. 
The main result about embeddable stochastic matrices that we will use below is the following:
\begin{theorem}[\cite{johansen73}*{Theorem 1.9}]\label{johansen}
A nonsigular stochastic matrix $D$ is embeddable if and only if it is a limit of products of elementary stochastic matrices. 
\end{theorem}
We remark on two consequences of this theorem that will be used below:
\begin{enumerate}   
\item a product of embeddable stochastic matrices is embeddable, 

\item the set of embeddable stochastic matrices is closed in the set of nonsingular stochastic matrices.
\end{enumerate} 

In order to relate  Markov processes to decreasing paths of polygons, we will find it convenient to first introduce a time-reversed form of the Markov process $D(\cdot,\cdot)$. Let $t_0>0$. Define 
\begin{align}\label{reversedMarkov}
E(t,s) &=D(t_0-t,t_0-s)\hbox{ for $0\le s\le t\le t_0$},\\
E(t) &=E(t,0)\hbox{ for $0\le t\leq t_0$}.
\end{align}
Note that \eqref{transitiveD} in terms of $E(\cdot, \cdot)$ becomes 
\[E(t,r)=E(t,s)E(s,r),\] for $t\geq s\geq r$.
If the function $t\mapsto E(t)$ is Lipschitz, then 
\[
\frac{d}{dt} E(t)=Q(t) E(t)
\]
where $Q(t)$ is a measurable, a.e. bounded intensity matrix.  Conversely, if $t\mapsto Q(t)$ is measurable, a.e. bounded on an interval $[0,t_0]$, and  such that $Q(t)$ is an intensity matrix for all  $t$, then the ODE above with $E(0)=I$ has a unique solution $t\mapsto E(t)$ ranging in the nonsingular stochastic matrices. Upon setting  $E(t,s)=E(t)E(s)^{-1}$ and $D(s,t)=E(t_0-t,t_0-s)$ we obtain a Markov process defined for $0\leq s\leq t\leq t_0$.

Now let $P$ be an $n$-gon regarded as an $n\times 2$ matrix, and define
\[
P(t)=E(t)P\hbox{ for $0\leq t\leq t_0$}.
\]
This is a decreasing path of polygons starting at $P$. Indeed,   $P(0)=E(0)P=P$ and 
\[
P(t)=E(t)P=E(t,s)E(s)P=E(t,s)P(s)
\] 
for $0\le s\le t\le t_0$. This shows that $P(t)\precsim P(s)$ for $t\ge s$. 

\begin{remark}
In \cite{goodmanCIME}, Goodman discusses the geometric interpretation of a Markov process as an increasing path of non-degenerate simplices. In the case of $n=3$, these are increasing paths of triangles.
\end{remark}

Our goal in the remainder of this section is to prove the following theorem:
\begin{theorem}\label{Plift}
Let $P$ be a convex $n$-gon oriented counterclockwise, and let $P'$ be an $n$-gon non-degenerately contained in $P$ and attainable from $P$. Then
$P'=DP$ for some  embeddable nonsingular stochastic matrix $D$.
\end{theorem}

\begin{corollary}\label{Pchattering}
Let $P$ and $P'$ be as in the previous theorem. Then $P'=\lim P_k'$ with $P_k'$  attainable from $P$ in finitely many pull-in moves for all $k$.
\end{corollary}

\begin{proof} 
By the previous theorem, $P'=DP$ for some embeddable stochastic matrix $D$. By Theorem \ref{johansen}, $D=\lim D_k$, with each $D_k$ a finite product of elementary stochastic matrices. Then $P'=\lim P_k'$, where $P_k'=D_kP$ and $P_k'$ is attainable from $P$ in finitely many pull-in moves for all $k$.
\end{proof}

The proof of Theorem \ref{Plift} will follow after some lemmas. 
For the remainder of this section we assume that $P$ is a convex $n$-gon oriented counterclockwise and that $P'$ is an $n$-gon non-degenerately contained in $P$ and attainable from $P$.

We call a path of polygons $ t\mapsto P(t)$ a Lipschitz path if  there exists a constant $L>0$ such that $\|p_i(t)-p_i(t')\|\leq L|t-t'|$ for all $t,t'$ and all $i$. (Here $\|\cdot\|$ denotes the Euclidean norm in $\R^2$.)

\begin{lemma}\label{lipschitzpath}
Suppose that there exists a Lipschitz decreasing  path of polygons  $\{P(t):0\leq t\leq t_0\}$ such that $P(0)=P$ and $P(t_0)=P'$. Then there exists a reversed Markov process $\{E(t):0\leq t\leq t_0\}$ such that $P(t)=E(t)P$ for all $0\leq t\leq t_0$. 
\end{lemma}

\begin{proof} 
Let $1\leq i\leq n$.
Since Lipschitz functions are a.e. differentiable,  $\frac{d p_i}{dt} (t)$ exists for almost all $t\in [0,t_0]$ and $\frac{d p_i}{dt}(t)$ is a bounded measurable function bounded by the Lipschitz constant of $p_i(t)$. 
For each  $0\leq t< t_0$ such that $\frac{d p_i}{dt}(t)$ exists, let   $\alpha_i(t)$ and $\beta_i(t)$ be the unique scalars such that
\begin{equation}\label{alphaibetai}
\frac{d p_i}{dt}(t)=\alpha_{i}(t)(p_{i+1}(t)-p_i(t))+\beta_i(t)(p_{i-1}(t)-p_i(t)).
\end{equation}
The functions $t\mapsto \alpha_i(t)$ and $t\mapsto \beta_i(t)$ are measurable and a.e. bounded on $[0,t_0]$, since they can be explicitly calculated as
\[
\begin{pmatrix}
\alpha_i\\
\beta_i\end{pmatrix}=\Gamma_i^{-1}
\cdot \frac{d p_i}{dt},
\]
where $\Gamma_i$ is the $2\times 2$ matrix  with column vectors  $p_{i+1}(t)-p_i(t)$ and $p_{i-1}(t)-p_i(t)$. 
(Note that $\det(\Gamma_i(t))\neq 0$ for all $0\leq t\leq t_0$, since $P(t)$ is non-degenerately contained in $P$ for all $t$.)

Let us now argue that $\alpha_i(t)$ and $\beta_i(t)$ are nonnegative scalars. Since the path $P(\cdot)$ is decreasing, $p_i(t')\in \co(P(t))$ for $t'\geq t$. Thus,  $p_{i}(t')-p_i(t)$ is a linear combination of   $p_{i+1}(t)-p_i(t)$ and $p_{i-1}(t)-p_i(t)$ with nonnegative coefficients. So for $t'>t$, we can write
\[
\frac{1}{t'-t}(p_{i}(t')-p_i(t))=\alpha_i(t',t)(p_{i+1}(t)-p_i(t)+\beta_i(t',t)(p_{i-1}(t)-p_i(t)),
\]
for nonnegative  scalars $\alpha_i(t',t)$ and $\beta_i(t',t)$. Since
\[
\begin{pmatrix}
\alpha_i(t',t)\\
\beta_i(t',t)\end{pmatrix}=\frac{1}{t'-t}\Gamma_i^{-1}(p_i(t')-p_i(t)),
\]
by letting $t'\to t^+$  we see that $\alpha_i(t',t)\to \alpha_i(t)$ and $\beta_i(t',t)\to \beta_i(t)$. Thus,
$\alpha_i(t)$ and $\beta_i(t)$ are nonnegative.
 
Let us rewrite \eqref{alphaibetai} as
\[
\frac{dp_i}{dt}(t)=\alpha_{i}(t)p_{i+1}(t) -(\alpha_i(t)+\beta_i(t))p_i(t)+\beta_i(t)p_{i-1}(t).
\]
Let $t\mapsto Q(t)$ be the $n\times n$ intensity  matrix 
\[
\begin{pmatrix}
-(\alpha_1+\beta_1) & \alpha_1 & 0 & \cdots & 0 & \beta_1 \\[6pt]
\beta_2 & -(\alpha_2+\beta_2) & \alpha_2 & \cdots & 0 & 0 \\[6pt]
0 & \beta_3 & \ddots & \ddots & \vdots & \vdots \\[6pt]
\vdots & \vdots & \ddots & \ddots & \alpha_{n-1} & 0 \\[6pt]
0 & 0 & \cdots & \beta_{n-1} & -(\alpha_{n-1}+\beta_{n-1}) & \alpha_{n-1} \\[6pt]
\alpha_n & 0 & 0 & \cdots & \beta_n & -(\alpha_n+\beta_n)
\end{pmatrix},
\]
defined so that  $\frac{d}{dt} P(t)=Q(t)P(t)$. 
Consider now the ODE on the $n\times n$ matrices
\[
\frac{d}{dt} E(t)=Q(t)E(t),
\]
with initial condition $E(0)=I$. This ODE has a unique solution $E(t)$ for $0\le t\le t_0$ (given by a Peano-Baker series) ranging in the nonsingular stochastic matrices. Upon defining 
$E(t,s)=E(t)E(s)^{-1}$ for $t_0\geq t\geq s\geq 0$ we obtain a reversed Markov process. Moreover, we have $P(t)=E(t)P$ for all $t$,  as both sides agree at $t=0$ and have the same derivative a.e.
\end{proof}

Lemma \ref{lipschitzpath} proves Theorem \ref{Plift} in the case that $P'$ is attainable from $P$ by a Lipschitz decreasing path. We will show next that if $P'$ is attainable from $P$, and not too far  from $P$ in a suitable sense, then it can be attained by a Lipszhitz decreasing path.
More explicitly, 
we work with the following notion of proximity of $P'$ to $P$: Let us say that $P'$ is \emph{in the midpoint pockets} of $P$ if for each $i$ the vertex $p_i'$ is contained in the triangle $(m_i,p_i,m_{i+1})$, where $m_i$ and $m_{i+1}$ are the midpoints of $[p_{i-1},p_i]$ and $[p_i,p_{i+1}]$ respectively, and $p_i\notin [m_i,m_{i+1}]$.

\begin{lemma}\label{midpointpockets}
Suppose that $P'$ is in the midpoint pockets of $P$. Then 
\[
\|p_i'-p_i\|\leq \frac4{\delta}(\area(P)-\area(P')),\hbox{ for all }i=1,\ldots,n,
\]
where $\delta$ is the minimum distance from any vertex of $P$ to a line determined by any two other vertices of $P$.
\end{lemma}
\emph{Note}: The midpoint pockets constraint implies that $P'$ is also a convex polygon. In particular, it has a well defined area. 

\begin{proof}
We assume without loss of generality that $i=1$.  

Extend the segment connecting $p_1$ to $p_1'$, until it intersects the segment  $[m_1,m_2]$ at the point $o$. See Figure \ref{fig_areas}.

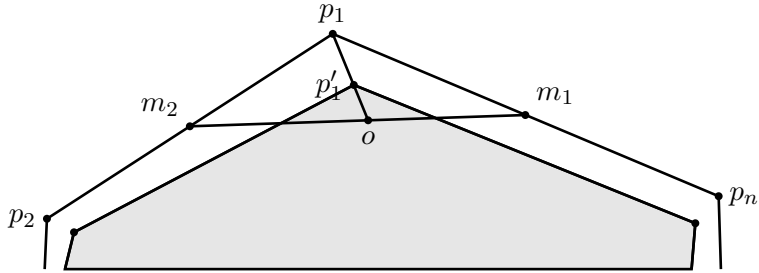
\begin{figure}[h]
\centering

\begin{tikzpicture}[line width=1pt,x=1.0cm,y=1.0cm]
 
    \coordinate (p3) at (-0.21,1.52);
    \coordinate (p2) at (-0.18,2.19);
    \coordinate (p1) at (3.62,4.65);
    \coordinate (pn) at (8.75,2.49);
    \coordinate (p1') at (3.9, 3.97);

    \coordinate (m2) at (1.72,3.42);
    \coordinate (m1) at (6.18,3.57);
    \coordinate (o) at (4.09,3.50);

    \coordinate (p2') at (0.18,2.01);
    \coordinate (pn') at (8.44,2.13);

    \coordinate (p3') at (0.06,1.52);
    \coordinate (p') at (8.39,1.52);
    \coordinate (p) at (8.78,1.52);

 \draw[fill=black!10!white] (p3') -- (p2') -- (p1') -- (pn') -- (p') -- cycle;
 
    \draw (p3) -- (p2) -- (p1) -- (pn) -- (p);
    \draw (p3') -- (p2') -- (p1') -- (pn') -- (p');
   
    \draw (m2) -- (m1);
    \draw (p1) -- (o);

    \foreach \point in {p2, p1, pn, p2', p1', pn', m1, m2, o} {
        \draw [fill=black] (\point) circle (1pt);
    }
    
    \node[left] at (p2) {$p_2$};
    \node[above] at (p1) {$p_1$};
    \node[right] at (pn) {$p_n$};
    \node[above left] at (m2) {$m_2$};
    \node[above right] at (m1) {$m_1$};
    \node[below] at (o) {$o$};
    \node[left] at (p1') {$p_1'$};

\end{tikzpicture}

\caption{}
\label{fig_areas}
\end{figure}

Then
\begin{align*}
\|p_1-p_1'\| &=\|o-p_1\| - \|o-p_1'\|\\
& = \|o-p_1\|\Big(1-\frac{\|o-p_1'\|}{\|o-p_1\|}\Big)\\
& = \|o-p_1\|\Big(1-\frac{\area(m_1,p_1',m_2)}{\area(m_1,p_1,m_2)}\Big) \\
& = \frac{\|o-p_1\|}{\area(m_1,p_1,m_2)} (\area(m_1,p_1,m_2)- \area(m_1,p_1',m_2)).
\end{align*}
We now estimate the RHS from above.
The points inside the triangle $(m_1,p_1,m_2)$ and outside of $(m_1,p_1',m_2)$ belong to $\co(P)\backslash \co(P')$. Thus,
\[
\area(m_1,p_1,m_2)- \area(m_1,p_1',m_2)\leq \area(P)-\area(P').
\]
On the other hand, since 
\[
\|o-p_1\| \leq \frac12\max(\|p_n-p_1\|,\|p_1-p_2\|)
\]
and 
\[
4\cdot \area(m_1,p_1,m_2) =\area(p_n,p_1,p_2),
\]
we have that
\[
\frac{\|o-p_1\|}{\area(m_1,p_1,m_2)} \leq 2\cdot \frac{\max(\|p_n-p_1\|,\|p_1-p_2\|)}{\area(p_n,p_1,p_2)}\leq  \frac4h,
\]
where $h$ is the shortest  altitude of the triangle $(p_n,p_1,p_2)$. We thus get that
\[
\|p_1-p_1'\| \leq \frac{4}{h}(\area(P)-\area(P'))\leq \frac{4}{\delta}(\area(P)-\area(P')).\qedhere
\]
\end{proof}

\begin{lemma}\label{reparametrize} 
Suppose that $P'$ is attainable from $P$ by a decreasing path $\{P(t):0\leq t \leq t_0\}$ such that  $P(t')$ is in the midpoint pockets of $P(t)$ for all $t'\geq t$. Then $P'$ is attainable from $P$ by a Lipschitz path.
\end{lemma}

\begin{proof} 
We obtain the new path as a reparametrization of the path $P(\cdot)$.
Observe that $t\mapsto \area( P(t))$ is decreasing and nonzero for all $t\geq 0$, as $ P(t)$ is  non-degenerate in $P$ for all $t$. Moreover,
\[
\area(P(t))=\area(P(t')) \Leftrightarrow P(t)=P(t').
\]
Define 
\[
h(t) := \area(P)-\area(P(t))\hbox{ for }t\in [0,t_0].
\]
Then, $h$ is non-decreasing, $h(0) = 0$, and 
$h(t') = h(t)$ if and only if  $P(t') = P(t)$. Set $s_0=h(t_0)$, and define a new path of polygons $\{R(s):0\leq s\leq s_0\}$ by
\[
 R(s) :=  P(h^{-1}(s))\hbox{ for  $s\in [0, s_0]$.} 
 \]
(Note that $h^{-1}(s)$ may possibly be an interval, but in this case the path $P(\cdot)$ is constant on that interval). We have that $R(0)=P$, $R(s_0)=P'$, and 
\begin{equation*}
\area(R(s))= \area( P(h^{-1}(s))) = \area(P)-h(h^{-1}(s))=\area(P)-s,
\end{equation*}
for $0\leq s\leq s_0$.

Let us show that $R(\cdot)$ is a Lipschitz path. Let $\delta(s)$ denote the minimum distance
from any vertex of $R(s)$ to a line determined by any two other vertices of $R(s)$. Let $\delta=\min\{\delta(s):s\in [0,s_0]\}$.
Since $P'=R(s_0)$ is non-degenerate in $P$, we must  have that $\delta>0$. 
Choose  $0\leq s < s' \leq s_0$. Since $R(s')$ is in the midpoint pockets of $R(s)$, Lemma \ref{midpointpockets} applied to the polygons  $R=R(s)$ and $R'=R(s')$ implies that
\[
\|r_i(s')-r_i(s)\|\leq \frac{4}{\delta}(\area(R(s))-\area(R(s')))=\frac{4}{\delta}(s'-s),
\]
where we have used that $\area(R(s))=\area(P)-s$ for all $s$. Thus, the path $R(\cdot)$ is Lipschitz with Lipschitz constant $4/\delta$.
\end{proof}

\begin{remark}
In \cite{goodman} Goodman showed that given a Markov process $D(\cdot,\cdot)$, a suitable reparametrization of $D(\cdot,\cdot)$  results in a Lipschitz Markov process. Thus,  if a stochastic matrix $D$ is embeddable, then it is also embeddable in a Lipschitz Markov process. 
\end{remark}

We are now ready to prove Theorem \ref{Plift}. 

\begin{proof}  [Proof of Theorem \ref{Plift}]
Let $\{P(t):0\leq t\leq t_0\}$ be a decreasing path starting at $P$ and ending at $P'$. Since $P'$ is non-degenerate in $P$, the distance from any vertex of $P(t)$ to a line formed by any two other vertices of $P(t)$ is bounded below for all $t$ by a positive number $\delta$. Using the uniform continuity of the path $P(\cdot)$, let us choose $\epsilon>0$ such that $\|p_{i}(t)-p_i(t')\|<\delta/2$ for all $|t-t'|<\epsilon$ and all $i$. 
Observe then that $\epsilon$ has the property that    $P(t')$ is in the midpoint pockets of $P(t)$ whenever  $0<t'-t<\epsilon$.

Let us subdivide the interval $[0,t_0]$ into times $0\leq s_1\leq s_2\leq\cdots\leq s_m=t_0$ such that $s_{i+1}-s_i<\epsilon$ for all $i$.
For each $1\leq i\leq m$, the path  $\{P(t):s_i\leq t\leq s_{i+1}\}$ is such that $P(t')$ is in the midpoint pockets of $P(t)$ for $s_i\leq t\leq t'\leq s_{i+1}$. Thus, by Lemma \ref{reparametrize}, $P(s_{i+1})$ is attainable from $P(s_i)$ by a Lipschitz path. By Lemma \ref{lipschitzpath}, 
there exists an embeddable stochastic matrix $D_i$ such that $P(s_{i+1})=D_iP(s_i)$. It follows that  $P'=D P$ with $D=\prod_{i=1}^{m-1} D_i$.  Since a finite product of embeddable stochastic matrices is again embeddable, the theorem  follows.
\end{proof}

It is possible to show, refining the methods used in this section, that the set of polygons attainable from a fixed polygon $P$ is a closed set.
We will not need to prove this en route to establishing our main results, but rather it will follow from those results in combination with the following
proposition. 

\begin{proposition}\label{compactNmoves}
For each $N\in \N$, the set of polygons attainable from $P$ in at most $N$ pull-in moves is a compact set.
\end{proposition}

\begin{proof}
Fix a vector of indices $\vec i=(i_1,j_1,i_2,j_2,\ldots,i_N,j_N)$, with $1\leq i_k,j_k\leq n$ and $i_k\neq j_k$ for all $k$.
Let $f_{\vec i}\colon [0,1]^N\to \PP_n$ be the function such that $f_{\vec i}(c_1,\ldots,c_N)$ is the polygon obtained from $P$
by successively applying pull-in moves of  $i_k$ toward $j_k$ with parameter $c_k$ for $k=1,\ldots,N$. Then $f_{\vec i}$ is a continuous function whose image is the compact set $K_{\vec i}$ of all polygons obtained from $P$ by performing pull-in moves as prescribed by the indices $\vec i$. The union $\bigcup_{\vec i} K_{\vec i}$, with $\vec i$ ranging through all $2N$-tuples of indices as above, is again a compact set, now  consisting of all polygons attainable from $P$ in $N$ pull-in moves.
\end{proof}

\section{Geometric preliminaries}\label{geoprelims}
In this section we lay down  some  geometric preliminaries that will be needed later on.

\subsection{Directed lines} 
Let $\ell$ be a line in $\R^2$ and $(a,b)$ an ordered pair of distinct points in $\ell$.  We call the map $\gamma\colon \R\to \ell$ defined by 
$\gamma(t)=(1-t)a+tb$ for all $t\in \R$ the affine parametrization of $\ell$ determined by $(a,b)$. This parametrization defines an order 
on $\ell$ induced by the order on $\R$. Thus, the choice of an ordered pair of distinct point $(a,b)$ in  $\ell$ induces an order on $\ell$ such that $a<b$.
 We call $\ell$ a \emph{directed line} if it is endowed with either one of the two possible orders induced by its affine parametrizations. We denote the directed line 
 determined by $(a,b)$ by $\ell_{ab}$.

 Let $\ell$ be a directed line and choose points $a,b\in\ell$ such that $a<b$ in the order of $\ell$. We say that a point $c\in \R^2$ lies to the left of $\ell$ if the triangle $(a,b,c)$ is either degenerate or oriented counterclockwise, i.e., if $\Delta(a,b,c)\geq 0$ where 
\[
\Delta(a,b,c)=\mathrm{det}
\begin{pmatrix}
a_1 & a_2 & 1\\
b_1 & b_2 & 1\\
c_1 & c_2 & 1
\end{pmatrix},
\]
and $a=(a_1,a_2)$, $b=(b_1,b_2)$, and $c=(c_1,c_2)$.
If $c$ is additionally not on $\ell$, then we say that $c$ lies strictly to the left of  $\ell$.
Similarly, $c$ lies to the right of $\ell$ if $(a,b,c)$ is oriented clockwise, and strictly to the right of $\ell$ if it lies to the right of $\ell$ and not on $\ell$. The points on a directed line lie both to the left and to the right of that line.

\subsection{Tangent rays}\label{righttangent}
A ray determines a directed line, so we may speak of the left and right sides of a ray, meaning the left and right sides of the directed line determined by the ray. Given distinct points $a,b\in \R^2$, we denote by $r_{ab}$ the ray emanating from $a$ and passing through $b$.

Let $S$ be a finite subset of $\R^2$. Given a point $x$ not in the interior of the convex hull of $S$, there exists a unique ray emanating from $x$, henceforth denoted by $r_x$, such that
\begin{enumerate}
 \item the set $S$ is entirely to the left of $r_x$, and
\item $r_x$ contains at least one point $p\in S\backslash \{x\}$.
\end{enumerate}
We call $r_x$ the \emph{right tangent ray} to $S$ emanating from $x$.
 
\begin{figure}[hp]

\begin{tikzpicture}[line width = 1pt,x=1.0cm,y=1.0cm]
\coordinate (A) at (0.61,3.10);
\coordinate (B) at (2.09,2.95);
\coordinate (C) at (0.51,4.44);
\coordinate (D) at (1.45,3.59);
\coordinate (E) at (2.17,4.25);
\coordinate (F) at (0.82,2.23);
\coordinate (H) at (3.2,3.58);

\draw (F) -- (H);

\node[left] at (F) {$x$};
\node[below] at (H) {$r_x$};

\foreach \p in {(A), (B), (C), (D), (E), (F)}{
\draw[fill=black] \p circle (1pt);}
\end{tikzpicture}\hspace{2cm}
\begin{tikzpicture}[line width = 1pt,x=1.0cm,y=1.0cm]

\coordinate (A) at (0.61,3.10);
\coordinate (B) at (2.09,2.95);
\coordinate (C) at (0.51,4.44);
\coordinate (D) at (1.45,3.59);
\coordinate (E) at (2.17,4.25);
\coordinate (F) at (-0.33,3.2);
\coordinate (H) at (3.6,2.8);
\coordinate (I) at (1.42,3.02);

\draw (I) -- (H);
\draw[dashed] (I) -- (F);

\node[below] at (I) {$x$};
\node[below] at (H) {$r_x$};
\foreach \p in {(A), (B), (C), (D), (E), (I)}{
\draw[fill=black] \p circle (1pt);}
\end{tikzpicture}

\caption{Tangent rays.}
\label{fig_perspectivities}
\end{figure}
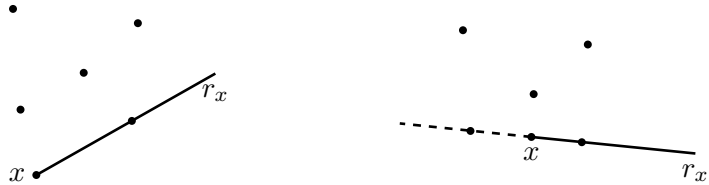 
 
\subsection{Perspectivities and projectivities}
Let $\ell_1$ and $\ell_2$ be distinct lines in $\R^2$ and $o$ a point not belonging to either one of them. The \emph{perspectivity from $\ell_1$ to $\ell_2$ with center $o$} is defined as the map $\alpha\colon \ell_1\to \ell_2$
such that 
\[
\alpha(x)=\ell_{ox}\cap \ell_2
\]
for all $x\in \ell_1$.  In the event that $\ell_{ox_0}$ is parallel to $\ell_2$ for some (unique) $x_0\in \ell_1$, we regard $\alpha(x_0)$ as undefined  and call $x_0$ the pole of $\alpha$. We use the notation $\ell_1\stackrel{o}{\wedge}\ell_2$ to denote  the perspectivity from $\ell_1$ to $\ell_2$ with center $o$.

The composition of two or more perspectivities between multiple lines in the plane is called a \emph{projectivity}, or projective transformation. 

\begin{remark} 
The natural domain for perspectivities and projectivities is the projective plane, as in this setting they are defined everywhere. 
More formally then, a projectivity is a composition of perspectivities in the projective plane. We will refrain however from explicitly working in the projective plane, since the concept of a decreasing path of polygons is naturally rooted in the affine plane.  Thus, we simply allow projectivities to be undefined at (at most) one point. 
\end{remark}

Let $\alpha\colon\ell_1\to\ell_2$ be a projectivity between directed lines. Let $\gamma_1\colon \R\to \ell_1$ and $\gamma_2\colon \R\to \ell_2$ be affine parametrizations of $\ell_1$ and $\ell_2$, respectively, compatible with their orders. Then 
\[
f=\gamma_2^{-1}\circ \alpha\circ\gamma_1
\]
is a fractional linear transformation from $\R$ to $\R$, i.e., it has the form 
\[
f(t)=\frac{at+b}{ct+d}
\]
for all $t\in \R$.  From this, the following facts are readily verified:
\begin{enumerate}
\item 
If $\alpha$ has no pole, then it is an affine transformation.

\item 
If $\alpha$ has a pole at $x_0\in\ell_1$, then we have one of the two following cases:
\begin{enumerate}
\item[(a)]
$\alpha$ is increasing and concave up on
$(-\infty,x_0)\subseteq \ell_1$, increasing and concave down on $(x_0,\infty)$, and
$\alpha(x)>\alpha(y)$ for all $x<x_0<y$. In this case we call $\alpha$ an \emph{orientation 
preserving} projectivity. 
\item[(b)]
$\alpha$ is decreasing and concave down on $(-\infty,x_0)\subseteq \ell_1$, decreasing and concave up on $(x_0,\infty)$, and $\alpha(x)<\alpha(y)$ for all $x<x_0<y$.  In this case we call $\alpha$ an \emph{orientation reversing} projectivity. 
\end{enumerate}
\item
If the segment $(a,b)\subset  \ell_1$ does not contain the pole of $\alpha$, then $\alpha$ maps $[a,b]$
bijectively to $[\alpha(a),\alpha(b)]$.

\item
If the segment $(a,b)\subset  \ell_1$ contains the pole of $\alpha$, then $\alpha$ maps $[a,b]$ onto $\ell_2\backslash (\alpha(a),\alpha(b))$ and $\ell_1\backslash (a,b)$  into 
$[\alpha(a),\alpha(b)]$.
\end{enumerate}
Given a function $\gamma\colon \ell_1\to \ell_2$ between directed lines $\ell_1$ and $\ell_2$, the properties of increasing and convex are defined relative to the orders on $\ell_1$ and $\ell_2$. To wit, $\gamma$ is increasing on a given interval of $\ell_1$ if $x\leq y$ implies $\gamma(x)\leq \gamma(y)$ for all $x,y$ in the interval, and convex if $\gamma(tx+(1-t)y)\leq t\gamma(x)+(1-t)\gamma(y)$ 
for all $t\in [0,1]$ and $x,y$ in the interval.

We will make use of the following two simple facts about perspectivities, which we state as  lemmas for ease of reference, and whose proofs we leave  to the reader:

\begin{lemma}\label{aboutperspectivities}
Let $\alpha=\ell_1\stackrel{o}{\wedge}\ell_2$, where $\ell_1$ and $\ell_2$ are (distinct) directed lines and $o$ lies strictly to the left of both of them. Then $\alpha$ is orientation preserving. If $o$ lies strictly to the left of one line and strictly to the right of the other line, then $\alpha$ is orientation reversing.
\end{lemma}

\begin{lemma}\label{pocketperspectivities}
Let $\alpha=\ell_1\stackrel{o}{\wedge}\ell_2$. Let $a,b\in \ell_1$ be points neither of which is the pole of $\alpha$. 

\begin{enumerate}[(i)]
\item 
If either $o\in [a,\alpha(a)]$ and $o\in [b,\alpha(b)]$, or $o\notin [a,\alpha(a)]$ and $o\notin [b,\alpha(b)]$, then the segment
$[a,b]$ does not contain the pole of $\alpha$, and consequently $[a,b]$ is mapped bijectively onto $[\alpha(a),\alpha(b)]$ by $\alpha$.
(See Figure \ref{fig_perspectivities}.)

\item
If either $o\in [a,\alpha(a)]$ and $o\notin [b,\alpha(b)]$, or $o\in [a,\alpha(a)]$ and $o\notin [b,\alpha(b)]$, then the segment
$[a,b]$  contains the pole of $\alpha$. Consequently, $\ell_1\backslash (a,b)$ is mapped into $[\alpha(a),\alpha(b)]$ by $\alpha$,
and $[a,b]$ is mapped onto $\ell_2\backslash (\alpha(a),\alpha(b))$.
(See Figure \ref{fig_perspectivities}.)
\end{enumerate}
\end{lemma}

\begin{figure}[h]
\centering

\begin{tikzpicture}[line width = 1pt,x=1.0cm,y=1.0cm]

\coordinate (a) at (0,2.5);
\draw[fill=black] (a) circle (1pt);

\coordinate (b) at (1.5,2.5);
\draw[fill=black] (b) circle (1pt);

\coordinate (aa) at (5.5,4.5);
\draw[fill=black] (aa) circle (1pt);

\coordinate (bb) at (6.5,6);
\draw[fill=black] (bb) circle (1pt);

%

\coordinate (b') at ($(b) + (b) - (a)$);

\draw [name path = ell1, shorten >= -4cm, shorten <= -1cm] (a) -- (b');
\draw [name path = ell2, shorten <= -1cm, shorten >= -3.8cm] (bb) -- (aa);
\draw [name path = ell3, shorten >= -1cm, shorten <= -1cm] (a) -- (aa);
\draw [name path = ell4, shorten <= -1cm, shorten >= -1cm] (b) -- (bb);

\path [name intersections={of=ell3 and ell4,by=o}];
\draw[fill=black] (o) circle (1pt);

\coordinate (o') at ($(o) + (aa)-(bb)$);

\draw[line width = 0.5pt, dashed, name path = dashell, shorten >= -0.5cm, shorten <= -4cm] (o) -- (o');
\path [name intersections={of=ell1 and dashell,by=p}];
\draw[fill=black] (p) circle (1pt);

\node [below] at (a) {$a$};
\node [below] at (b) {$b$};
\node [below] at (p) {pole};
\node [below] at (o) {$o$};
\node [right] at (aa) {$\alpha(a)$};
\node [right] at (bb) {$\alpha(b)$};
\end{tikzpicture}
\quad
\begin{tikzpicture}[line width = 1pt,x=1.0cm,y=1.0cm]

\coordinate (a) at (0,2.5);
\draw[fill=black] (a) circle (1pt);

\coordinate (b) at (1.5,2.5);

\coordinate (newb) at (3.25,2.5);
\draw[fill=black] (newb) circle (1pt);

\coordinate (aa) at (5.5,4.5);
\draw[fill=black] (aa) circle (1pt);

\coordinate (bb) at (6.5,6);
\coordinate (farbb) at ($(aa) - 2.75*(bb) + 2.75*(aa)$);
\draw [name path = ell1, shorten >= -3cm, shorten <= -1cm] (a) -- (newb);

\draw [name path = ell2, shorten <= -1cm] (aa) -- (farbb);

\draw [name path = ell3, shorten >= -1cm, shorten <= -1cm] (a) -- (aa);
\draw [draw = none, name path = ell4, shorten <= -1cm, shorten >= -1cm] (b) -- (bb);

\path [name intersections={of=ell3 and ell4,by=o}];
\draw[fill=black] (o) circle (1pt);

\coordinate (o') at ($(o) + (aa)-(bb)$);
\draw[line width = 0.5pt, dashed, name path = dashell, shorten >= -1cm, shorten <= -2cm] (o) -- (o');
\path [name intersections={of=ell1 and dashell,by=p}];
\draw[fill=black] (p) circle (1pt);

\coordinate (farbp) at ($(o) + 2.9*(newb) - 2.9*(o)$);

\draw[name path = ell5, shorten < = -1cm] (o) -- (farbp);
\path [name intersections={of=ell2 and ell5,by=newbb}];
\draw[fill=black] (newbb) circle (1pt);

\node [below] at (a) {$a$};
\node [below left] at (newb) {$b$};
\node [below] at (p) {pole};
\node [below] at (o) {$o$};
\node [right] at (aa) {$\alpha(a)$};
\node [right] at (newbb) {$\alpha(b)$};

\end{tikzpicture}

\caption{Left and right figures illustrate the situations in Lemma \ref{pocketperspectivities} (i) and (ii).}
\label{fig_perspectivities}
\end{figure}
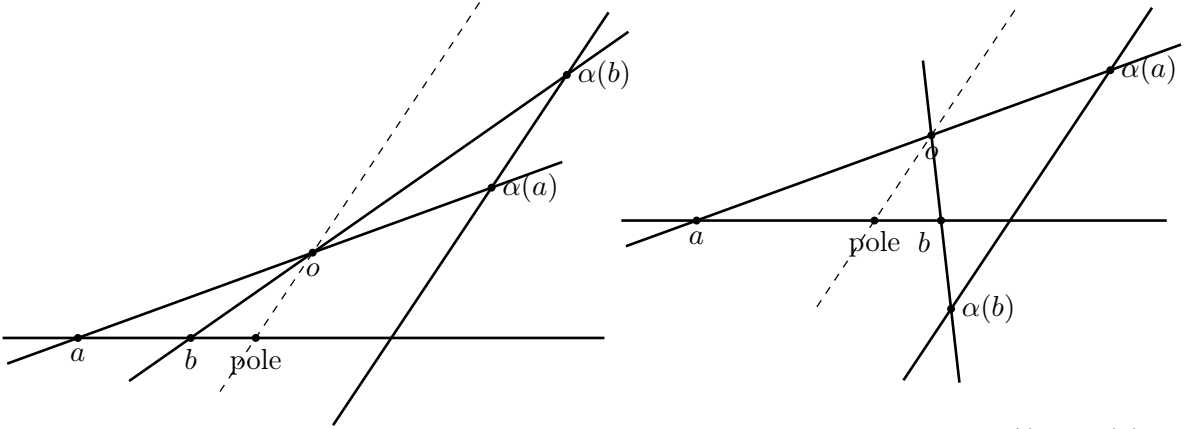

\section{The Poncelet map}
Let $n\geq 3$. If an $n$-gon $P$ is not set-convex, then every $n$-gon $P'\precsim P$ is degenerate in $P$, and so attainable in at most $5n$ moves, by Theorem \ref{in5nmoves}.
This solves both the Attainablity and Bang-Bang problems for a non-set-convex $P$.  In this section we initiate the investigation of these problems assuming that $P$ is set-convex. After a suitable re-indexing of the vertices, a set-convex polygon can be turned into a convex polygon oriented counterclockwise (see Remark \ref{onreindexing}). We thus assume throughout this section that  $P=(p_k)_{k=1}^n$ is a convex $n$-gon oriented counterclockwise.

We fix an $n$-gon $P'=(p_i')_{i=1}^n$ contained in $P$. We henceforth assume that not all vertices of $P'$ are collinear. 

Recall that $\partial P$ denotes the boundary of $P$, i.e, the union of its edges. Next we define a function  $\pi\colon\partial P \to \partial P$. This is the Poncelet map that titles this section.

\begin{definition}[The Poncelet map]\label{ponceletdef}
For any $x\in \partial P$, let $r_x$ be the unique right tangent ray to $P'$ emanating from $x$. We define $\pi(x)\in\partial P$ as the point in  $r_x\cap \partial P$ which is furthest from $x$. 
\end{definition}

\begin{remark}
We can define a a function $\pi_{\cw}\colon \partial P\to \partial P$ by simply replacing ``right tangent ray" with ``left tangent ray" in the definition of $\pi$. We call $\pi_{\cw}$ the clockwise Poncelet map. Our focus will be largely on $\pi$. With the aid of a reflection on $\R^2$, results on $\pi$ are easily transferred to $\pi_{\cw}$.  
\end{remark}

We will now further characterize $\pi$. Let $x\in \partial P$. Since $P'$ is contained in $P$, and $P$ is convex, $x$ does not belong to the interior of the convex hull of $P'$. Thus, the existence of the right tangent ray $r_x$ is guaranteed (see \ref{righttangent}). There are two cases to consider
in the evaluation of $\pi$:

\noindent\textbf{Interior Case}: The ray $r_x$ intersects the interior of $P$. In this case, $\pi(x)$ is the unique 
point distinct from $x$ that is both on $r_x$ and on $\partial P$.

\noindent\textbf{Boundary Case}:
The ray $r_x$ does not intersect the interior of $P$. In this case  $r_x$ is collinear with the half-open edge $[p_i, p_{i+1})$ containing $x$, and so either $r_x=r_{xp_i}$ (with $x\neq p_i$) or $r_x=r_{xp_{i+1}}$. However, if $r_x=r_{xp_i}$, since $P'$ lies to the left of both $r_x$ and $\ell_{p_ip_{i+1}}$, the entire polygon $P'$ would have to lie on the edge $[p_i, p_{i+1}]$, a scenario that we have ruled out by assumption. Therefore,  $r_x=r_{xp_{i+1}}$ and $\pi(x)=p_{i+1}$.

\begin{figure}[h]
\centering

\begin{tikzpicture}[line width=1pt, x=1.0cm,y=1.0cm]
\draw  (0,0) coordinate (p1) -- (2.65,0) coordinate (p2) -- (2.84,3.37) coordinate (p3) -- (0.49,3.92) coordinate (p4) -- (-1.12,2.38) coordinate (p5) -- cycle;

\draw  (2.71,1.08) coordinate (p1') -- (2.78,2.34) coordinate (p2') -- (0.66,3.03) coordinate (p3') -- (-0.01,0.91) coordinate (p4') -- cycle;

\foreach \p in {p1, p2, p3, p4, p5, p1', p2', p3', p4'} {
\draw[fill=black] (\p) circle (1pt);}
\draw[fill=black] (0.84, 1.5) circle (1pt);

\coordinate (x) at (1.4, 3.71);
\draw[fill=black] (x) circle (1pt);
\node[above] at (x) {$x$};

\coordinate (pix) at (-0.79, 1.69);
\draw[fill=black] (pix) circle (1pt);
\node[left] at (pix) {$\pi(x)$};

\draw[dashed] (x) -- (pix);

\node[left] at (p3') {$p(x)$};

\end{tikzpicture}
\hspace{1cm}
\begin{tikzpicture}[line width=1pt, x=1.0cm,y=1.0cm]
\draw  (0,0) coordinate (p1) -- (2.65,0) coordinate (p2) -- (2.84,3.37) coordinate (p3) -- (0.49,3.92) coordinate (p4) -- (-1.12,2.38) coordinate (p5) -- cycle;

\draw  (2.71,1.08) coordinate (p1') -- (2.78,2.34) coordinate (p2') -- (0.66,3.03) coordinate (p3') -- (-0.01,0.91) coordinate (p4') -- cycle;

\foreach \p in {p1, p2, p3, p4, p5, p1', p2', p3', p4'} {
\draw[fill=black] (\p) circle (1pt);}
\draw[fill=black] (0.84, 1.5) circle (1pt);

\coordinate (x) at (2.68, 0.56);
\draw[fill=black] (x) circle (1pt);
\node[right] at (x) {$x$};
\node[right] at (p3) {$\pi(x)$};
\node[right] at (p2') {$p(x)$};
\end{tikzpicture}

\caption{Interior and boundary cases in the evaluation of the Poncelet map.}
\label{fig_cases}
\end{figure}
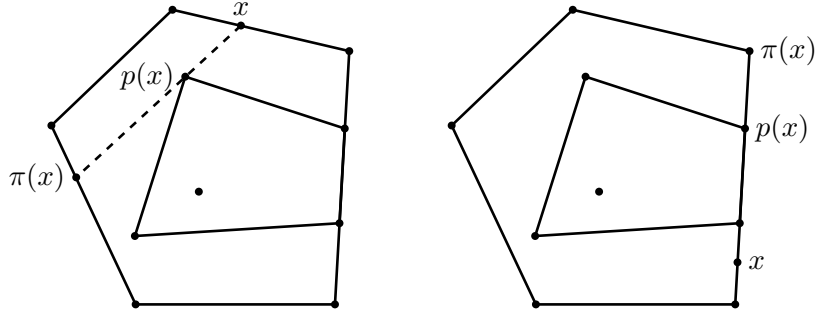

\begin{remark}\leavevmode
\begin{enumerate}[(i)]
\item
In either the interior or boundary cases of the evaluation of $\pi$, $x$ and $\pi(x)$ belong to different half-open edges $[p_i,p_{i+1})$ of $P$. In particular, $\pi$ has no fixed points.
\item (Pivots.)
In either the interior or boundary cases of the evaluation of $\pi$, the open ray $r_x\backslash\{x\}$ passes through one or two vertices of $\co(P')$. Any such vertex of $\co(P') $ will be called a  \emph{pivot} in the evaluation of $\pi$ at $x$. 
Pivots are necessarily vertices of $P'$. We shall denote by $p(x)$ the pivot furthest from $x$ along $r_x$.
\end{enumerate}
\end{remark}

Given distinct points $a,b\in\partial P$, let us denote by $\arc[a,b]$ the points on the counterclockwise arc along $\partial P$ from $a$ to $b$ including both $a$ and $b$. (We follow standard notational conventions for arcs excluding one or both endpoints.)  Note that the complement of $\arc[a,b]$ in $\partial P$ is $\arc(b,a)$. By the convexity of $P$, if the line $\ell_{ab}$ crosses the interior of $P$, then $\arc[a,b]$ coincides with the points in $\partial P$ that lie to the right of $\ell_{ab}$.

\begin{lemma} \label{Fablemma} 
Let $a,b\in \partial P$ be distinct points such that $P'$ lies to the left of the directed line $\ell_{ab}$. Let $c\in\partial P$.
\begin{enumerate}[(i)]
\item 
$\pi$ maps   $\arc[a,b]$ into   $\arc[b,a)$. 

\item 
$\pi$ maps $\arc[c,\pi(c)]$ into $\arc[\pi(c),c)$. 

\item 
If $z\in \arc(\pi(c),c]$, then $\pi(z)\in \arc(z,\pi(c)]$.
\end{enumerate}
\end{lemma}

\begin{proof}
(i) Let us show first that $\pi$ maps $\arc[a,b]$ into   $\arc[b,a]$. Assume that $a$ and $b$ belong to the same edge, say $[p_i,p_{i+1}]$, of $P$. We cannot have $b<a$ in $\ell_{p_ip_{i+1}}$, since this would mean that $P'$, which lies to the left of $\ell_{ab}$, is contained in the segment 
$[p_i,p_{i+1}]$, but the vertices of $P'$ are not all collinear. Thus, $\arc[a,b]\subseteq [p_i,p_{i+1}]$. Then, for each $x\in \arc[a,b]$, 
either $\pi(x)$ belongs to a different edge $[p_j,p_{j+1}]$ of $P$  or $\pi(x)=p_{i+1}$. In either case  $\pi(x)\in \arc[p_{i+1},p_i]\subseteq \arc[b,a]$, as desired.  Assume now that $a$ and $b$ belong to different edges of $P$. In this case   $\arc[a,b]$ agrees with the points of $\partial P$ that lie to the right of $\ell_{ab}$. Let $x\in \arc[a,b]$ and   suppose for contradiction that $\pi(x)$ lies strictly to the right of $\ell_{ab}$. Then every point on the half-open segment $(x,\pi(x)]$ lies also strictly to the right of $\ell_{ab}$. However, the segment $(x,\pi(x)]$ contains at least one vertex of $P'$, namely the pivot $p(x)$. This contradicts the assumption that $P'$ lies to the left of $\ell_{ab}$. Thus, $\pi(x)$ lies to the  left of $\ell_{ab}$, i.e., $\pi(x)\in \arc[b,a]$.

Finally, let us show that $\pi$ does not attain the value $a$  on $\arc[a,b]$. Suppose that $\pi(x)=a$ for some $x\in \arc[a,b]$. Then, on the one hand,  $\pi(a)\in \arc[b,a]$, as shown in the previous paragraph.
On the other hand, since $a=\pi(x)\in \arc[x,\pi(x)]$, we have $\pi(a)\in \arc[\pi(x),x]=\arc[a,x]$, where again we have applied the result established in the previous paragraph to  $\arc[x,\pi(x)]$. It follows that $\pi(a)\in \arc[b,a]\cap \arc[a,x]$. If $x\neq b$, this means that $\pi(a)=a$, contradicting that $\pi$ has no fixed points. If $x=b$, then  $P'$ lies to the left of both $\ell_{ab}$ and $\ell_{ba}$, hence it is contained in $\ell_{ab}$, contradicting that not all vertices of $P'$ are collinear. Thus, $\pi(x)\neq a$ for all $x\in \arc[a,b]$.

(ii) Since $c\neq \pi(c)$ (as $\pi$ has no fixed points), and 
$P'$ is to the left of the directed line $\ell_{c,\pi(c)}$, we can apply (i) with $a=c$ and $b=\pi(c)$.

(iii) Note that $\arc(z,\pi(c)]$ is well defined, since  $z\neq \pi(c)$ by assumption. 

As we know, $\pi(z)\neq z$. Suppose for a contradiction that $\pi(z)\in \arc(\pi(c),z))$. 
This implies that   $c\in [z,\pi(z)]$ and $\pi(c)\in \arc(z,\pi(z))$. However, from  $c\in [z,\pi(z)]$  and (ii) we deduce that
$\pi(c)\in \arc(\pi(z), z)$. This contradicts our assumption.
\end{proof}

A function $f\colon \partial P\to \partial P$ is called orientation preserving if for all  $x\neq y$ such that $f(x)\neq f(y)$, $f$ 
maps $\arc[x,y]$ into $\arc[f(x),f(y)]$. 

\begin{proposition}\label{piorientation}
Any function $f\colon\partial P\mapsto \partial P$ that maps $\arc[x,f(x)]$ into $\arc[f(x),x]$ for all nonfixed points $x$ of $f$ is orientation preserving. In particular, the Poncelet map is orientation preserving.    
\end{proposition}

\begin{proof}
Let $f\colon\partial P\mapsto \partial P$ be such that it maps $\arc[x,f(x)]$ into $\arc[f(x),x]$ for all $x$ not fixed by $f$.
Observe that the property assumed for $f$ is precisely (ii) of the previous lemma. In the previous lemma, (iii) follows from (ii). We thus have that $f(z)\in\arc[z,f(x)]$ for all $x$ not fixed by $f$ and $z\in \arc(f(x),x]$. 

We will make use of another property of $f$ regarding fixed points: Let $y\in \partial P$ be fixed by $f$.
If $x$ is not fixed by $f$, then $\arc(x,f(x))$ does not contain any fixed points of $f$ (as it is mapped into its complement).
It follows that $y\in \arc[f(x),x]$, i.e., $f(x)\in \arc[x,y]$. This conclusion holds  for all $x\neq y$ (fixed or not fixed by $f$).

Let $x\neq y$ be such that $f(x)\neq f(y)$, and let $z\in \arc(x,y)$. We wish to show that  $f(z)\in \arc[f(x),f(y)]$. 
We examine several cases.

\textbf{Case that $y$ is fixed by $f$}. As remarked above, in this case  $f(x)\in\arc[x,y]$ and $f(z)\in \arc[z,y]$. We either have that  $z\in \arc(x,f(x))$ or that $z\in [f(x),y)$. If $z\in\arc[x,f(x)]$, then 
\[
f(z)\in \arc[f(x),x]\cap\arc[z,y]=\arc[f(x),f(y)].
\] 
If $z\in \arc [f(x),y]$, then $f(z)\in \arc[z,y]\subseteq \arc[f(x),f(y)$. 

\textbf{Case that $x$ is fixed by $f$}. If $f(z)=f(x)$, then we certainly have that $f(z)\in\arc[f(x),f(y)]$, as desired.
Assume thus that $f(z)\neq f(x)$. Notice that  $y\in \arc(z,x)$ and that $x$ is fixed by $f$. We can thus apply the previous case 
to the triple of points $z,y,x$. We deduce that $f(y)\in \arc[f(z),f(x)]$, i.e., that $f(z)\in \arc[f(x),f(y)]$, as desired.

In the remaining cases we assume that neither $x$ nor $y$ is a fixed point of $f$.

\textbf{Case  $x\in \arc[y,f(y)]$}. In this case $f(x)\in \arc[f(y),y]$. If $z\in \arc[x,f(x)]$, then
\[
f(z)\in \arc[f(x),x]\subseteq \arc[f(x),f(y)].
\]
If on the other hand $z\in \arc(f(x),y]\subseteq \arc(f(y),y]$, then 
\[
f(z)\in \arc[z,f(y)]\subseteq \arc[f(x),f(y)].
\]

\textbf{Case $x\in \arc(f(y),y)$}. 
In this case, $f(x)\in [x,f(y)]$. We  either have that $f(x)\in (x,y]$ or that $f(x)\in (y,f(y))$. We analyze each of these cases next.

\textbf{Subcase $f(x)\in (x,y]$}. If $z\in \arc[x,f(x)]$, then
\[
f(z)\in \arc[f(x),x]\cap \arc[z,f(y)]=\arc[f(x),f(y)],
\] 
and if $z\in \arc (f(x),y]\subseteq (f(y),y]$, then 
\[
f(z)\in \arc[z,f(y)]\subseteq \arc[f(x),f(y)].
\]

\textbf{Subcase $f(x)\in (y,f(y))$}. Then
\[
f(z)\in \arc[f(x),x]\cap \arc[z,f(y)]=\arc[f(x),f(y)].\qedhere
\]
\end{proof}

\begin{lemma}\label{rst}
Let $r,s,t$ be consecutive vertices of $\co(P')$, listed in the counterclockwise order. Suppose that $s$ is located in the interior of $P$. The following statements are true:

\begin{enumerate}[(i)]
\item
The line $\ell_{rs}$ intersects $\partial P$ at points $a,\pi(a)$, and $\ell_{st}$ intersects $\partial P$ at points $b,\pi(b)$ respectively. Moreover, these points can be listed in counterclockwise order as $a,b,\pi(a),\pi(b)$.

\item
The set of points $x\in\partial P$ such that $p(x)=s$ is the half-open $\arc[a,b)$, and $s$ is the unique pivot of the points in the open arc $\arc(a,b)$. 

\item
$\pi$ maps the  $\arc[a,b]$ bijectively onto  $\arc[\pi(a),\pi(b)]$.

\end{enumerate}
\end{lemma}

\begin{figure}[h]
\centering

\begin{tikzpicture}[line width=1pt, x=1.0cm,y=1.0cm, scale=1.15]
\draw  (0,0) coordinate (p1) -- (2.65,0) coordinate (p2) -- (2.84,3.37) coordinate (p3) -- (0.49,3.92) coordinate (p4) -- (-1.12,2.38) coordinate (p5) -- cycle;

\coordinate (v) at (0.45, 0.76);
\coordinate (r) at (1.81, 1.06);
\coordinate (s) at (1.44, 2.53);
\coordinate (t) at (0.42, 2.85);
\coordinate (u) at (-0.03, 1.71);

\filldraw[fill=black!10!white]  (r) -- (s) -- (t) -- (u) -- (v) -- cycle;

\foreach \point in {t, r, s, u, v}
    \draw[fill=black] (\point) circle (1pt);
    
\foreach \point in {p1, p2, p3, p4, p5}
    \draw[fill=black] (\point) circle (1pt);    

\coordinate (a) at (2.08, 0);
\coordinate (pia) at (1.12, 3.77);
\coordinate (b) at (2.77, 2.12);
\coordinate (pib) at (-0.38, 3.09);

\foreach \point in {a, pia, b, pib} {
\draw[fill=black] (\point) circle (1pt);}

\draw[dashed, shorten < = -0.5cm, shorten > = -1cm] (a) -- (pia);
\draw[dashed, shorten < = -0.5cm, shorten > = -1cm] (b) -- (pib);

\node[below] at (a) {$a$};
\node[above] at (pia) {$\pi(a)$};
\node[above right] at (b) {$b$};
\node[above] at (pib) {$\pi(b)$};
\node[right] at (r) {$r$};
\node[above right] at (s) {$s$};
\node[above] at (t) {$t$};
\end{tikzpicture}

\caption{The points on $\arc[a,b]$, with common pivot $s$, are mapped bijectively onto $\arc[\pi(a),\pi(b)]$ by the Poncelet map.}
\end{figure}
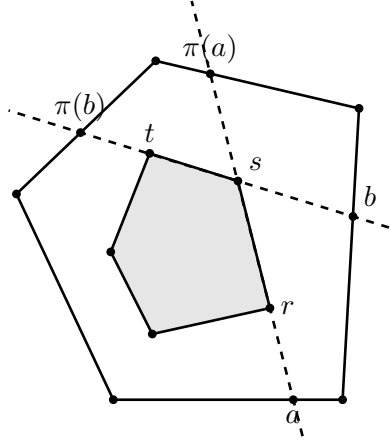

\begin{proof}
(i) Since $s$ is an interior point, $\ell_{rs}$ intersects $\partial P$ at exactly two points (\cite[Corollary 8.16]{lee}). Let $a$ and $a'$ be the two points on  $\partial P$ such that $\ell_{aa'}=\ell_{rs}$ (as directed lines). Now, $P'$ lies to the left of $\ell_{rs}$, by the convexity of $\co(P')$ and the fact that  $r$ and $s$ are consecutive vertices in the counterclockwise order. Thus, the ray $r_{as}$ is the right tangent ray to $P'$ emanating from $a$
(note that $s\neq a$).  Since  $s$ is in the interior of $P$, we are in the interior case of the evaluation of $\pi$, and  $a'=\pi(a)$. 

Similarly, the line $\ell_{st}$ intersects $\partial P$ at exactly two points $b$ and $b'$, with $\ell_{bb'}=\ell_{st}$ and $b'=\pi(b)$. Moreover, since $t$
lies to the left of $\ell_{rs}$, $\pi(b)$ lies to the left of $\ell_{a\pi(a)}$ and $b$ lies to the right of this line too. Thus, $b\in \arc[a,\pi(a)]$
and $\pi(b)\in \arc[\pi(a),a]$, which implies the stated counterclockwise ordering of the points.

(ii) Since $r\in [a,s]$, $r$ and $s$ are the two pivots for the point $a$, with $s$ being the pivot furthest from $a$. Hence, $p(a)=s$. Similarly, $s$ and $t$ are pivots for $b$, but in this case $t$ is furthest from $b$, so $p(b)=t$.

Let $x\in \arc(a,b)$. Let $x'$ be the unique point distinct from $x$ where $r_{xs}$ intersects $\partial P$. Since $x$ lies strictly to the right of $\ell_{rs}$ and strictly to the left of $\ell_{st}$,  $x'$ lies strictly to the left of $\ell_{rs}$ and strictly to the right of $\ell_{st}$. Therefore, from (i), $x'\in \arc(\pi(a),\pi(b))$. Note also that $[x,s)$
is strictly right of $\ell_{rs}$ and $(s,x']$ is strictly right of $\ell_{st}$. It follows that $r_{xs}$ is the right tangent ray to $P'$ at $x$, and $s$ is its unique point of contact with $P'$. Thus, $x'=\pi(x)$ and $s$ is the only pivot in the evaluation of $\pi$ at $x$. In particular, $p(x)=s$.

Let us show that $s$ is not a pivot for $x\notin \arc[a,b]$. If $x\in \arc(b,\pi(b))$, then $t$ lies strictly to the right of $\ell_{xs}$. So $s$ cannot be a pivot for such $x$. Similarly, if $x\in \arc[\pi(b),a)\subseteq \arc(\pi(a),a)$, then $r$ lies strictly to the right of $\ell_{xs}$, and again $s$ cannot be a pivot for such points.

(iii) Injectivity is clear, since for distinct points $x,y\in \arc[a,b]$, we have that $\pi(x)\in r_{xs}$, $\pi(y)\in r_{ys}$, and the rays $r_{xs}$ and $r_{ys}$ intersect only at the interior point $s$. To show surjectivity, consider any $x'\in \arc(\pi(a),\pi(b))$. The line $\ell_{sx'}$ intersects $\partial P$ at exactly one other point $x\in \arc(a,b)$. The ray $r_{xs}$ is the right tangent ray to $P'$ emanating from $x$. Since $r_{xs}$ intersects $\partial P$ at $x'$, we must have that $x'=\pi(x)$, thus showing surjectivity.
\end{proof}

\begin{definition}\label{Gammasets}
Define the following subsets of $\partial P$:
\begin{enumerate}[(i)]
\item
$\Gamma_1$ is the set of all $x\in \partial P$ such that the right tangent ray $r_x$ to $P'$ intersects the interior of $P$ (interior case of the evaluation of $\pi$) and $r_x$ contains at least two vertices of $P'$.
\item
$\Gamma_2$ is the set of all $x\in \partial P$ such that $r_x$ intersects the interior of $P$ (interior case of the evaluation of $\pi$) and $\pi(x)$ is a vertex of $P$.
\item 
$\Gamma=P\cup \Gamma_1\cup \Gamma_2$. (Here $P$ is regarded as  its set of vertices.)
\end{enumerate}
\end{definition}

Note that $\Gamma$ is a finite set. Note also that
if $\gamma$ and $\gamma'$ are consecutive points of $\Gamma$ in the counterclockwise order,  then the segment $[\gamma,\gamma']$ is contained in a single edge of $P$, since the vertices of $P$ belong to $\Gamma$.

The set $\Gamma$ consists of all critical junctures for $\pi$ as $x$ travels counterclockwise along $\partial P$, in the sense specified by Theorem \ref{gamma} below.
We remind the reader that by the perspectivity with center $o$ from line  $\ell_1$ to line $\ell_2$  we understand  the map 
$\ell_1\ni x\mapsto  \ell_{ox}\cap \ell_2\in \ell_2$.

\begin{theorem}\label{gamma}
Let $\gamma$ and $\gamma'$ be two consecutive points of $\Gamma$ in the counterclockwise orientation. One of the following cases occurs:
\begin{enumerate}[(i)]
\item
$\pi$ is constant on the segment $[\gamma,\gamma')$ and equal to a vertex of $P$. 
\item
$\pi$ is constant on the segment $[\gamma,\gamma')$ and equal to a vertex of $P'$ in $\partial P$. 
\item
$\pi$ is a perspectivity from $[\gamma,\gamma']$ to another edge of $P$ with center of perspectivity a vertex of $\co(P')$ in the interior of $P$. 
Moreover, $\pi$ is increasing  on $[\gamma,\gamma']$.
\end{enumerate}
\end{theorem}

\begin{proof}
 Let us assume that $\gamma$ belongs to the half-open edge $[p_j,p_{j+1})$ of $P$.  
Then  $\gamma'\in (\gamma,p_{j+1}]$. 

We will prove the theorem in three cases.

Suppose that the evaluation of $\pi$ on $\gamma$ is in the boundary case. Then  $\pi(\gamma)=p_{j+1}$ and  $p(\gamma)$ is the vertex of $P'$ on the segment $(\gamma,p_{j+1}]$ that is furthest from $\gamma$. If $p(\gamma)=p_{j+1}$, then $p(\gamma)\in \Gamma$. If, on the other hand, $p(\gamma)\in (\gamma,p_{j+1})$, then the evaluation of $\pi$ on $p(\gamma)$ is in the interior case. Furthermore, the right tangent ray $r_{p(\gamma)}$ contains a vertex of $P'$ (a pivot) other than $p(\gamma)$. Thus, $p(\gamma)\in \Gamma_1$, and again we have that $p(\gamma)\in \Gamma$. It follows that $\gamma'\in (\gamma,p(\gamma)]$. It is then clear that $p(x)=p(\gamma)$ for all $x\in [\gamma,\gamma')$ and that $\pi$ is constant equal to $p_{j+1}$ on $[\gamma,\gamma')$.  This confirms case (i) of the statement of the theorem, thus proving the theorem when the evaluation of $\pi$ on $\gamma$ is in the boundary case.

Suppose now that the evaluation of $\pi$ on $\gamma$ is in the interior case, i.e., the right tangent ray $r_{\gamma}$ intersects the interior of $P$, and that $p(\gamma)\in \partial P$. In this case $\pi(\gamma)=p(\gamma)$.  Let $[p_k,p_{k+1})$ be the half-open edge of $P$ containing $p(\gamma)$. Since $\gamma$ and $\pi(\gamma)$  belong to different half-open edges of $P$, $k\neq j$. 
Thus, $p_{j+1}$ is strictly to the right of $\ell_{\gamma,p(\gamma)}$. 
For each point $x\in \arc[\gamma,p_{j+1})$, it is not hard to show that  $r_{xp(\gamma)}$ is the right tangent ray to $P'$ emanating from $x$ and that this ray intersects the interior of $P$, and so  $p(\gamma)$ is the unique point of intersection 
of $r_{xp(\gamma)}$ and $\partial P$ distinct from $x$. It follows that $\pi(x)=p(\gamma)$.
In this case, $\pi$ is constant on $[\gamma,p_{j+1})$ equal to $p(\gamma)$. Since $\gamma'\in (\gamma,p_{j+1}]$, $\pi$ is constant on $[\gamma,\gamma')$ equal to $p(\gamma)$.  This is case (ii) of the statement of the theorem.

Finally, let us suppose that the evaluation of $\pi$ at $\gamma$ is in the interior case and that $p(\gamma)$ is in the interior of $P$. 
Then, by Lemma \ref{rst}, the set of all $x$ such that $p(x)=p(\gamma)$ is a half-open $\arc[\gamma,\gamma'')$, where $\gamma''\in \Gamma_1\subset \Gamma$. It follows that $\gamma'\in \arc[\gamma,\gamma'']$. Thus, $p(x)=p(\gamma)$ for all $x$ in $[\gamma,\gamma')$.

The image of the open segment $(\gamma,\gamma')$ under $\pi$ is the open arc $\arc(\pi(\gamma),\pi(\gamma'))$. This arc cannot contain vertices of $P$, for if the evaluation of $\pi$ at some $x\in (\gamma,\gamma')$ agrees with a vertex of $P$, then  $x\in\Gamma_2\subseteq \Gamma$, contradicting that $\gamma$ and $\gamma'$ are consecutive points of $\Gamma$. 
Therefore, $\arc(\pi(\gamma),\pi(\gamma'))$ is contained within a single edge $[p_k,p_{k+1}]$ of $P$. We have already established that $[\gamma,\gamma']$ is contained in the edge $[p_j,p_{j+1}]$. Thus, on the segment $[\gamma,\gamma']$, the function
$\pi$ coincides with the perspectivity 
\[
\alpha=\ell_{p_jp_{j+1}}\stackrel{p(\gamma)}{\wedge} \ell_{p_kp_{k+1}}
\] 
from $\ell_{p_jp_{j+1}}$ to $\ell_{p_kp_{k+1}}$ with center $p(\gamma)$.
 Additionally, since $p(\gamma)$ is an interior point of $P$, it lies to the left both of $\ell_{p_jp_{j+1}}$ and of $\ell_{p_kp_{k+1}}$.
 Thus, the perspectivity $\alpha$ is increasing  on $[\gamma,\gamma']$, by  Lemma \ref{aboutperspectivities}. 
\end{proof}

Consider  a function $f\colon \partial P\to \partial P$. Let us say that $f$ is piecewise projective if
there exists a finite set of \emph{juncture points} $S\subseteq \partial P$ containing the vertices of $P$ and such that for any two $s,s'\in S$, consecutive in the counterclockwise order, $f|_{[s,s']}$ agrees with a projectivity (without a pole on $[s,s']$) from the line $\ell_{ss'}$ to a line containing an edge of $P$. Observe that piecewise projective functions are always continuous.

\begin{corollary}\label{pwprojectivepi}
If all the vertices of $P'$ belong to the  interior of $P$, then  $\pi$ is piecewise projective with set of juncture points $\Gamma$,
and bijective with inverse the clockwise Poncelet map $\pi_{\cw}$.  
\end{corollary}

\begin{proof}
Suppose that the vertices of $P'$ belong to the interior of $P$.
Then we only encounter case (iii) of Theorem \ref{gamma}, as cases (i) and (ii) entail the existence of a vertex of $P'$ on $\partial P$. Thus, $\pi$ is piecewise projective with set of juncture points $\Gamma$. 

Let $x\in\partial P$ and set $x'=\pi(x)$. Since the ray  $r_{xx'}$
emanating from $x$ and passing through $x'$ is right tangent to $P'$, $P'$ lies to the left of this ray. Further, the segment $(x,x')$ contains at least one vertex of $P'$, e.g., $p(x)$. (Note that $p(x)\neq x'$ since $P'$ is in the interior of $P$.) Thus, the ray $r_{x'x}$, emanating from $x'$ and passing through $x$, is the left tangent ray to $P'$. Since this ray passes through the interior of $P$, it intersects $\partial P$ at exactly two points: $x$ and $x'$. It follows that $\pi_{\cw}(x')=x$. A symmetric argument shows that $\pi$ is a left inverse of $\pi_{\cw}$. 
\end{proof}

\begin{remark}
If some vertices of $P'$ belong to the boundary of $P$,
then $\pi$ is neither continuous, nor injective, nor surjective.
\end{remark}

 In the following lemma we make use  of the following observation whose simple proof we omit:  If $f_1\colon \partial P\to \partial P$ and $f_2\colon \partial P\to \partial P$  are bijective piecewise projective functions, then $f_2\circ f_1$ is also piecewise projective. Moreover, if $S_1\subseteq \partial P$ and $S_2\subseteq \partial P$ are  sets of juncture points for $f_1$ and $f_2$ respectively, then $S=f_2^{-1}(S_2)\cup S_1$ is a set of juncture points for $f_2\circ f_1$.

\begin{lemma}\label{junctures}
Suppose that every vertex of $P'$ belongs to the interior of $P$. Then for each $k\geq 1$, the function  $\pi^k$ is piecewise projective with set of juncture points
\[
\Gamma\cup \pi^{-1}(\Gamma)\cup \cdots \pi^{-(k-1)}(\Gamma).
\] 
\end{lemma}

\begin{proof}
By Corollary \ref{pwprojectivepi}, $\pi$ is piecewise projective with set of juncture points $\Gamma$. Repeatedly applying the observation made in the preceding paragraph, we arrive at the set of juncture points for  $\pi^k$ in the statement of the lemma.
\end{proof}

\section{The broken line construction}
In this section we  continue to assume that $P$ is a convex  $n$-gon oriented counterclockwise, and that $P'$ is an $n$-gon contained in $P$ and such that  not all vertices of $P'$ are collinear. We continue to denote by  $\pi\colon \partial P\to \partial P$ the Poncelet map  relative to $P'$ defined in the previous section. 

\begin{definition}[Broken Line Construction]\label{BLC}
We define the broken line construction (BLC) relative to $P'$ with starting point $x\in\partial P$ as follows: 
Set $x_1=x$ and $x_2=\pi(x_1)$. For $k\geq 2$, if $\pi(x_k)$ lies on $\arc(x_k,x_1)$, i.e., the open counterclockwise arc from $x_k$ to $x_1$, set $x_{k+1}=\pi(x_k)$; otherwise, stop the construction. 
\end{definition}

We note that, by Lemma \ref{Fablemma} (ii), $\pi^2(x)\in \arc(\pi(x),x)$ for all $x\in \partial P$. This guarantees that
 there are always at least three  distinct points in the broken line construction.

\begin{remark}
We can similarly construct points along $\partial P$  from an initial $x\in \partial P$ by applying the clockwise Poncelet map $\pi_{\cw}$ instead of $\pi$,  following a similar iterative process, and using clockwise arcs for the stopping condition. We refer to this construction as the clockwise  BLC.  
\end{remark}

Let $a_1,...,a_m$ be distinct points on the boundary $\partial P$, indexed in the counterclockwise order. These  points partition  $\partial P$ into disjoint arcs which we call partition arcs:
\[
\partial P=\arc[a_1,a_2)\sqcup \arc[a_2,a_3)...\sqcup \arc[a_m,a_1).
\]

\begin{lemma}[Partition Arcs Lemma]\label{pockets} Let $m\geq 3$.
Consider distinct points $a_1,a_2,...,a_m$  on  $\partial P$, listed in the counterclockwise order, such that $P'$ lies to the left of  $\ell_{a_ia_{i+1}}$ for all $i$. The broken line construction, started at a point $x_1\in \arc[a_1,a_2)$, has the following property: For each $k\geq 1$, if $x_k$ belongs to $\arc[a_i,a_{i+1})$, with $1\leq i \leq m$, and if $x_{k+1}$ is defined and belongs to $\arc[a_j,a_{j+1})$, then either
\begin{enumerate}[(1)]
\item
 $i<j\leq m$, or
\item 
$j=1$ and  $x_{k+1}$ is the final point in the construction.
\end{enumerate}
In the second case, $x_{k+1}\in \arc[a_1,x_1)$ and $\pi(x_{k+1})\in \arc[x_1,x_2]$.
\end{lemma}

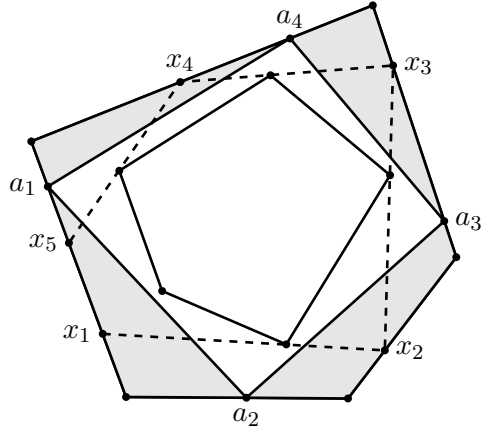
\begin{figure}[h]
\centering

\begin{tikzpicture}[line width = 1pt, x=1.0cm, y=1.0cm]

\draw (-0.10, -1.30) coordinate (p1) -- (2.85, -1.32) coordinate (p2) -- (4.29, 0.56) coordinate (p3) -- (3.18, 3.91) coordinate (p4) -- (-1.36, 2.10) coordinate (p5) -- cycle;

\coordinate (a1) at (-1.14, 1.50);
\coordinate (a2) at (1.50, -1.31);
\coordinate (a3) at (4.13, 1.04);
\coordinate (a4) at (2.08, 3.47);

\draw[fill=black!10!white] (a2) -- (p2) -- (p3) -- (a3) -- cycle;
\draw[fill=black!10!white] (a3) -- (p4) -- (a4) -- cycle;
\draw[fill=black!10!white] (a4) -- (p5) -- (a1) -- cycle;
\draw[fill=black!10!white] (a1) -- (p1) -- (a2) -- cycle;

\draw (1.82, 2.98) coordinate (p1') -- (-0.19, 1.71) coordinate (p2') -- (0.38, 0.11) coordinate (p3') -- (2.03, -0.60) coordinate (p4') -- (3.41, 1.65) coordinate (p5') -- cycle;

\draw[dashed] (-0.41, -0.46) coordinate (x1) -- (3.34, -0.68) coordinate (x2)  -- (3.45, 3.11) coordinate (x3) -- (0.62, 2.89) coordinate (x4) -- (-0.86, 0.75) coordinate (x5);

\foreach \point in {p1, p2, p3, p4, p5, a1, a2, a3, a4, p1', p2', p3', p4', p5', x1, x2, x3, x4, x5}{
\draw[fill=black] (\point) circle (1pt);
}

\node[left] at (x1) {$x_1$};
\node[right] at (x2) {$x_2$};
\node[right] at (x3) {$x_3$};
\node[above] at (x4) {$x_4$};
\node[left] at (x5) {$x_5$};

\node[left] at (a1) {$a_1$};
\node[below] at (a2) {$a_2$};
\node[right] at (a3) {$a_3$};
\node[above] at (a4) {$a_4$};
\end{tikzpicture}

\caption{BLC with four partition arcs. Only the arc containing the first BLC point may contain another BLC point, namely the last one.}
\end{figure}

\begin{proof} 
By Lemma \ref{Fablemma} (i), $x_{k+1}=\pi(x_k)$ does not belong to $\arc[a_i,a_{i+1})$, that is, $j\neq i$. Also, since $x_k$ is not the last point of the construction,  $x_{k+1}\in \arc(x_k,x_1)$. For $i'=2,\ldots,i-1$, the partition arc $\arc[a_{i'},a_{i'+1})$ is contained in the arc $\arc(x_1,a_i)$, which in turn is contained in  $\arc[x_1,x_k]$. Thus, $x_{k+1}$ does not belong to any of these partition arcs. We thus have that either $i<j\leq m$ or $j=1$. 

Suppose that $j=1$, i.e., $x_{k+1}\in \arc[a_1,a_2)$. Since $x_{k+1}$ is defined in the construction, we further have that $x_{k+1}\in \arc[a_1,x_1)$. By Lemma \ref{Fablemma} (i), $\pi(x_{k+1})\in \arc[a_2,a_1)$. This implies that  $\pi(x_{k+1})\notin \arc[x_{k+1},x_1)$; hence, $x_{k+1}$ is the final point in the construction. Since $x_{k+1}\in \arc(x_2,x_1)$, we have $\pi(x_{k+1})\in \arc(x_{k+1},x_2]$ by Lemma \ref{Fablemma} (iii). But $\pi(x_{k+1})\notin \arc[x_{k+1},x_1)$, so $\pi(x_{k+1})\in \arc[x_1,x_2]$.
\end{proof}

If we choose the points $a_1,\ldots,a_m$ in the above lemma to be the vertices of $P$, so that the partition arcs they create  are half-open edges of $P$, we see that the points  $x_1,x_2,\ldots$ obtained in the BLC  belong to different half-open edges of $P$, except possibly for the first and last points in the construction. In particular, the  BLC always yields a finite sequence of points $x_1,x_2,\ldots,x_l$.  We refer to $Q=(x_i)_{i=1}^l$ as the $l$-gon resulting from the BLC starting at $x$.

\begin{theorem}\label{allaboutQ}
Let $x\in \partial P$, and assume that $x$ belongs to the half-open edge $[p_i,p_{i+1})$. Let $Q=(x_k)_{k=1}^{l}$ be the $l$-gon obtained from the BLC relative to $P'$ starting at the point $x\in \partial P$. Then $Q$ has the following properties:
\begin{enumerate}[(i)]
\item
$Q$ is inscribed in $P$ and its vertices follow the counterclockwise order along $\partial P$.

\item
Every half-open edge $[p_j,p_{j+1})$ of $P$ contains at most one vertex of $Q$, except possibly  $[p_i,p_{i+1})$, which may contain $x_l$ and $x_1$. In particular,
$l\leq n+1$.

\item  
$Q$ is convex, unless $x_l\in [p_i,p_{i+1})$ and $x_2=p_{i+1}$, in which case  the polygon $Q'=(x_2,x_3,...,x_l)$ is convex.

\item
$P'\precsim Q$, i.e., $P'$ is contained in $Q$.

\item
For $k=1,\ldots, l-1$, the half-open edge $(x_{k},x_{k+1}]$ of $Q$ contains the vertex $p(x_k)$ of $P'$. Since the half-open edges of $Q$ are pairwise disjoint,  the BLC  uses at least $l-1$ distinct vertices of $P'$ as pivots.
\end{enumerate}
\end{theorem}
\begin{proof}
Let us assume without loss of generality that $i=1$, i.e., that $x\in [p_1,p_2)$. 

(i) It is clear from the construction of $Q$ that its vertices belong to $\partial P$, i.e., that it is inscribed in $P$. It is also clear from construction that since $x_{k+1}\in \arc(x_k,x_1)$ for all $k$,  the vertices of $Q$ follow the counterclockwise order along $\partial P$. 

(ii) This follows immediately from  the partition arcs lemma applied with $a_1,\dots,a_m$ the vertices of $P$.

(iii) Let $1\leq k \leq l$. Due to the counterclockwise ordering of the vertices of $Q$,  $x_j\in \arc(x_{k+1},x_k)$ for $j\neq k,k+1$, which by the convexity of $P$ means that $x_j$ lies to the left of $\ell_{x_kx_{k+1}}$. Let us argue that the only three consecutive vertices of $Q$ that can be collinear are $x_l,x_1,x_2$. 

Indeed, for $k\neq l,1$  the vertices $x_{k-1},x_k,x_{k+1}$ of $Q$ belong to 
different half-open edges  $[p_i,p_{i+1})$ of $P$, by (ii), which implies that not all three of them can belong to the same edge.

If $k=l$, and $x_{l-1}$ and $x_l$ belong to the same edge $[p_j,p_{j+1}]$ of $P$, then  $x_l=p_{j+1}$, as $x_{l-1}$ and $x_l$ must belong to different half-open edges. This prevents $x_1$ from being on this same edge. Thus, the only possibly collinear vertices of $Q$ are $x_{l},x_1,x_2$.

Suppose that $x_l,x_1,x_2$ are not collinear. Then  for each $k$ the vertices of $Q$  lie  to the left of $\ell_{x_kx_{k+1}}$ and no three consecutive vertices are collinear. This precisely means that $Q$ is convex  and oriented counterclockwise. 

If $x_l,x_1,x_2$ are collinear, then they belong to the same edge $[p_1,p_2]$. Then $P'$ necessarily  contains a vertex on $[p_1,p_2]$, namely the pivot $p(x_1)$ in the evaluation of $\pi$ at $x_1$. So $x_2=\pi(x_1)=p_2$. In this case we see that $x_l,x_2,x_3$ cannot be collinear, as $x_3\notin [p_1,p_2]$. Hence, as before, we conclude that $Q'$ is convex.

(iv)  By construction, $P'$ lies to the left of $\ell_{x_kx_{k+1}}$ for $1\leq k<l$, since  $x_{k+1}=\pi(x_k)$ for such $k$ and $P'$
lies to the left of $\ell_{x\pi(x)}$ for any $x\in \partial P$. For the same reason, $P'$ lies to the left of the line $\ell_{x_l,\pi(x_l)}$. But, by the stopping condition of the BLC, and since $l\geq 3$, $x_1 \in \arc(x_{l},\pi(x_{l})]$. Thus, $P'$ also lies to the left of $\ell_{x_lx_1}$. 

From part (iii), either $Q$ is convex or $Q'$ is convex. If $Q$ is a convex polygon, then $\co(Q)$ is precisely the set of points left of each of its edges. 
Hence  $P'\precsim Q$. If $Q$ is not convex, we saw that $Q'$ is convex.  Then   $P'$ is left of every edge of   $Q'$ (since in this case $\ell_{x_lx_2}=\ell_{p_1p_2}$ and $P$ is convex). So again we obtain that $P'\precsim Q'\precsim Q$, as desired.

(v) This is clear from the definition  of $Q$ and the properties of $\pi$: if $1\leq k < l$ then $x_{k+1}=\pi(x_{k+1})$, and in this case the half-open edge $(x_k,x_{k+1}]$ contains the pivot  $p(x_k)$.
\end{proof}

\section{Degeneracy test}
Throughout this section we let $P$ be a convex $n$-gon oriented counterclockwise and we let  $P'$ be an $n$-gon contained in $P$. As in the previous section, we assume that not all vertices of $P'$ are collinear. We denote by  $\pi\colon \partial P\to \partial P$ the Poncelet map relative to $P'$.

In Theorem \ref{degeneracytestthm} below we identify a finite set of ``test"  points on $\partial P$ from which to start the BLC in order to determine whether $P'$ is degenerately contained in $P$. Before proving this theorem we establish a couple of lemmas. 

Let us call $x\in \partial P$ a \emph{good point} if the BLC  (relative to $P'$) starting at $x$ results in an $l$-gon with $l<n$, and call it bad otherwise.

\begin{lemma}\label{Fxgood}
If $x$ is good, then $\pi^k(x)$ is good for all $k\in \N$. 
\end{lemma}

\begin{proof}
Note that $x$ is good if and only if $\pi^l(x)\in \arc[x,\pi(x)]$ for some $1<l<n$. Let $x$ be a good point and set $x'=\pi(x)$. Let us show that $x'$ is good too. We have:
\[
\pi^l(x')=\pi^l(\pi(x))=\pi(\pi^l(x)).
\]
Since $\pi$ is orientation preserving (Proposition \ref{piorientation}) and $\pi^l(x)\in \arc[x,\pi(x)]$, it follows that 
\[
\pi^l(x')\in \arc[\pi(x),\pi^2(x)]=\arc[x',\pi(x')].
\]
Thus, $x'$ is a good point.
\end{proof}

\begin{lemma}\label{m+1verticeslemma}
Let $Q$ be an $m$-gon inscribed in $P$, with pairwise distinct vertices,  and such that  $P'\precsim Q\precsim P$. The following statements are true:
\begin{enumerate}[(i)]
\item
The BLC starting at any point on $\partial P$ produces an interpolating polygon with at most $m+1$ vertices.
\item
The BLC starting at any of the vertices of $Q$ produces an interpolating polygon  with at most $m$ vertices. 
\end{enumerate} 
\end{lemma}

\begin{proof} 
Let $Q=(q_k)_{k=1}^m$, and assume, after re-indexing if necessary,  that the vertices of $Q$ follow the counterclockwise order along $\partial P$. Let $x\in \partial P$, and let $(x_k)_{k=1}^l$ be the polygon resulting from the BLC with initial point $x_1=x$.

(i)  Since $Q$ lies to the left of $\ell_{q_kq_{k+1}}$ for all $k$, $P'$  lies also to the left of $\ell_{q_kq_{k+1}}$ for all $k$.  The vertices of $Q$ create $m$ partition arcs  $\arc[q_k,q_{k+1})$, with $k=1,\ldots,m$, along $\partial P$. By Lemma \ref{pockets}, one of these partition arcs may contain $x_1$ and $x_l$, while each of the others contains at most one point from the BLC.  It follows that $l\leq m+1$, thus proving (i). 
 
 (ii) Suppose now that $x$ is a vertex of $Q$, say $x=q_1$. In this case $x_l$ cannot belong to  $\arc[q_1,q_2)$, since   $x_l\in \arc(x_{l-1},x_1)$.  Thus each of the $m$ partition arcs created by $Q$  contains at most one point from the BLC, implying that $l\leq m$.
\end{proof}

Define a set  $T\subseteq \partial P$ as 
\[
T= P\cup \Gamma_1,
\]
where  $P$ stands for its set of vertices and $\Gamma_1$ is as in Definition \ref{Gammasets}. 
Each point in  $\Gamma_1$ is the push-out onto $\partial P $ of some vertex $v$ of $\co(P')$ by its successor $v'$ (in the counterclockwise order), with $v$ and $v'$ not both belonging to the same  edge of $P$.

\begin{theorem}\label{degeneracytestthm} 
The $n$-gon $P'$ is degenerately contained in $P$ if and only if there exists $x\in T$ such that the BLC starting at $x$ results in an $l$-gon with $l<n$.
\end{theorem}

\begin{proof}
One direction of the theorem is obvious: Any $l$-gon  obtained by the BLC is an interpolating polygon between $P'$ and $P$. 
If $l<n$, then this verifies that $P'$ is degenerately contained in $P$. 

Suppose now that $P'$ is degenerately contained in $P$.
Let us further assume that $n\geq 4$, since in the case $n=3$ degenerate containment is equivalent  to collinearity, which we have assumed not to hold for $P'$. Our goal is to show that the set $T$ contains at least one good point.

We first dispose of several cases in which the theorem can be settled quickly.

Suppose that there exists an interpolating $m$-gon  $P'\precsim Q\precsim P$, with $m\leq n-2$. By pushing out the vertices of $Q$ onto $\partial P$ if necessary
(with push-out moves), we may assume that $Q$ is inscribed in $P$. Then, by Lemma \ref{m+1verticeslemma} (i), the  polygon resulting from the BLC starting at any $x\in \partial P$ has at most $m+1\leq n-1$ vertices. Thus, every point in $\partial P$, including the vertices of $P$, is a good point,  proving the theorem in this case.

Suppose that there exists an interpolating $m$-gon  $P'\precsim Q\precsim P$, where $m<n$ and   at least one vertex of $Q$ coincides with a vertex of $P$.
By pushing out the remaining vertices of $Q$ onto $\partial P$ if necessary, we may assume that $Q$ is inscribed in $P$. 
By Lemma \ref{m+1verticeslemma} (ii), every vertex of $Q$ is a good point. In particular, the common vertex of $Q$ and $P$, which belongs to the set $T$, is a good point, thus proving the theorem in this case.

Suppose that there exists an interpolating $m$-gon  $P'\precsim Q\precsim P$, where $m<n$ and  at least two vertices of $Q$ belong to the same edge of $P$.  
We can push out one of these two vertices with the other one onto a vertex of $P$. We are now in the previously considered case.

Suppose that  at least one vertex of $P'$ belongs to $\partial P$.  Say that $p_i'$ is on the edge $[p_j,p_{j+1}]$. 
Choose an interpolating $m$-gon  $P'\precsim Q\precsim P$  inscribed in $P$ and with $m<n$.  
Since $[p_j,p_{j+1}]$ is a face of $\co(P)$, $p_i'$ is a convex combination of the vertices of $Q$
in $[p_j,p_{j+1}]$. If at least two  vertices of $Q$ belong to $[p_j,p_{j+1}]$, we are  in a  previously considered case, hence done. 
Assume, then, that only one vertex of $Q$ belongs to $[p_j,p_{j+1}]$, which must be $p_i'$. Then $p_i'$ is a a good point, as  every vertex of $Q$  is good (by Lemma \ref{m+1verticeslemma} (ii)). If $p_i'=p_{j+1}$, then $T$ contains good points, as $P\subseteq T$.
If the evaluation  $\pi(p_i')$ is in the interior case, then $p_{i}'\in \Gamma_1$ is a good point in $T$. Finally, if $p_i'\neq p_{j+1}$ and the evaluation $\pi(p_i')$ 
is in the boundary case, the $\pi(p_i')=p_{j+1}$, and it follows that $p_{j+1}\in T$ is a good point (by Lemma \ref{Fxgood}).

From this point on we assume that $P'$ does not fall in any of the four cases covered above. We thus assume that the vertices of $P'$ are all in the interior of $P$, and that any interpolating polygon $P'\precsim Q\precsim P$ witnessing the degeneracy of $P'$ in $P$ is an $(n-1)$-gon not sharing  a vertex with $P$ and with no two vertices on the same edge of $P$.

Since the vertices of $P'$ lie in the interior of $P$, the Poncelet map $\pi$ is a piecewise projective function with  set of juncture points   $\Gamma$ (Corollary \ref{pwprojectivepi}). By Lemma \ref{junctures}, $\pi^{n-1}$ is also piecewise projective, and has set of juncture points 
\[
T' = \Gamma \cup \pi^{-1}(\Gamma) \cup \cdots \cup \pi^{-(n-2)}(\Gamma).
\]
If $T'$ contains a good point,  then $\Gamma$ or a preimage of $\Gamma$ contains a good point. Since the set of good points is closed under evaluation by $\pi$ (Lemma \ref{Fxgood}),  it follows that $\Gamma$ contains a good point. If $x\in \Gamma$ is a good point, then  either $x$ is a vertex of $P$, $x\in \Gamma_1$, or $\pi(x)$ 
is a vertex of $P$. In all cases, $T$ contains a good point. Thus, to prove that $T$ contains a good point it suffices to show that $T'$ contains a good point. We will show this next.

We note that the set of good points is non-empty, since the vertices of  any $(n-1)$-gon inscribed in $P$ and containing $P'$ must be good points by Lemma \ref{m+1verticeslemma} (ii). Let $x_1\in \partial P$ be a good point. Let  $Q=(x_k)_{k=1}^{n-1}$ be the $(n-1)$-gon obtained from the BLC with starting point $x_1$. 
Since no vertex of $Q$ is a vertex of $P$, and no two vertices of $Q$ belong to the same edge of $P$ (by our earlier assumption), exactly  one edge of $P$ is missing a vertex of $Q$. Moving $x_1$ forward along the BLC, let us replace it by the point right before the edge of $P$ that is skipped, so that the edge is skipped going from $x_1$ to $x_2$. This new location for $x_1$ is still that of a good point, by Lemma \ref{Fxgood}. Let us now assume, without loss of generality, that $x_1\in (p_n,p_1)$, so that the skipped edge is $[p_1,p_2]$. Then we have 
\[
x_k\in (p_k,p_{k+1})
\] 
for  all $k=2, \ldots,n-1$.

Let $c$ and $c'$ be consecutive points of $T'$ in the counterclockwise order such that $x_1\in [c,c']$. Notice that $c$ and $c'$ belong to the edge $[p_n,p_1]$, since $T'$ contains the vertices of $P$. We will be done once we have shown that either $c$ or $c'$ is a good point.

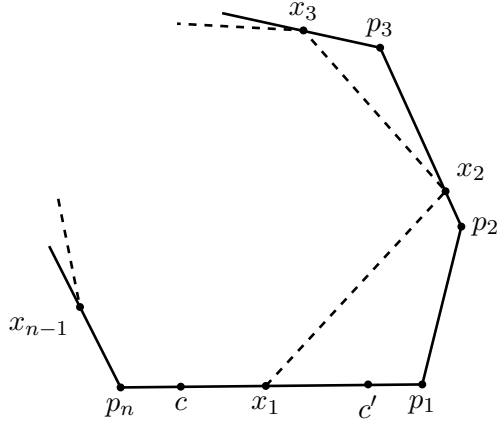
\begin{figure}[h]
\centering

\begin{tikzpicture}[line width = 1pt, x=1.0cm, y=1.0cm]

\coordinate (pn) at (0.64,-1.08);
\coordinate (p1) at (4.65,-1.04);
\coordinate (p2) at (5.17,1.06);
\coordinate (p3) at (4.09,3.44);
\coordinate (p4) at (1.99,3.9);
\coordinate (p5) at (-0.31,0.8);

\coordinate (O) at (1.39, 3.76);
\coordinate (P) at (-0.19, 1.43);

\coordinate (x1) at (2.57,-1.06);
\coordinate (x2) at (4.96,1.53);
\coordinate (x3) at (3.07,3.67);
\coordinate (x4) at (0.1, -0.01); 

\coordinate (c)  at (1.44, -1.07);
\coordinate (c') at (3.93,-1.04);

\draw (p5) -- (pn) -- (p1) -- (p2) -- (p3) -- (p4);

\draw[dashed] (x1) -- (x2) -- (x3) -- (O);
\draw[dashed] (P) -- (x4);

\foreach \point in {pn, x1, c, c', p1, p2, x2, p3, x3, x4}
{
\draw[fill=black] (\point) circle (1pt);
}

\foreach \point/\label in {pn/$p_n$, c/$c$, x1/$x_1$, c'/$c'$, p1/$p_1$}{
\node[below] at (\point) {\label};
}

\node[right] at (p2) {$p_2$};
\node[above right] at (x2) {$x_2$};
\node[above] at (p3) {$p_3$};
\node[above] at (x3) {$x_3$};
\node[below left] at (x4) {$x_{n-1}$};

\end{tikzpicture}

\caption{The BLC starting at the good point $x_1\in [c,c']$ skips $[p_1,p_2]$ and stops at $x_{n-1}\in [p_{n-1},p_n]$.}
\end{figure}

By the stopping condition of the BLC, $\pi^{(n-1)}(x_1)=\pi(x_{n-1})$ belongs to $\arc[x_1,x_2]$ (Lemma \ref{pockets}). Since $x_1$ lies on the edge $[p_n,p_1]$, and $x_2$ on the edge $[p_2,p_3]$, this leaves three possibilities for the edge of $P$ containing $\pi^{(n-1)}(x_1)$: $[p_n,p_1]$, $[p_1,p_2]$, or $[p_2,p_3]$. 
Suppose that $\pi^{(n-1)}(x_1)\in [p_1,p_2]$. Since $c$ and $c'$ are consecutive juncture points of the (piecewise projective) function $\pi^{(n-1)}$, 
the point $\pi^{(n-1)}(x)$ does not change edges while  $x$ ranges in $[c,c']$. It follows that  $\pi^{(n-1)}(c)$ belongs to $[p_1,p_2]$.
This, in turn, implies that $c$ is a good point, as the BLC starting at $c$ would stop before reaching $\pi^{(n-1)}(c)$, 
  completing the proof in this case. The same analysis shows that if   $\pi^{(n-1)}(x_1)\in [p_2,p_3]$, then
$c$ is a good point, and we are again done. Let us thus assume, for the remainder of the proof, that $\pi^{(n-1)}(x_1)$ belongs to $[p_n,p_1]$. Consequently, 
$\pi^{(n-1)}(x)$ belongs to $[p_n,p_1]$ for all $x\in [c,c']$.

Denote by $f\colon \ell_{p_np_1}\to \ell_{p_np_1}$ the projectivity agreeing with  $\pi^{(n-1)}$  on the interval $[c,c']$.
Note that a point $y\in [c,c']$ is good if and only if $f(y)\geq y$ in the order of $\ell_{p_np_1}$. In particular, since $x_1$ is a good point, $f(x_1)\geq x_1$. The proof will be complete once we have shown that either $f(c)\geq c$ or $f(c')\geq c'$.

Write 
\begin{equation*}
f=\beta_{n-2}\circ \cdots \circ\beta_{1}\circ\alpha,
\end{equation*} 
where $\alpha=\ell_{p_np_1}\wedge\ell_{p_2p_3}$ is the first perspectivity used in the BLC starting at $x_1$, and 
$\beta_k=\ell_{p_{k+1}p_{k+2}}\wedge\ell_{p_{k+2}p_{k+3}}$ for $k=1,\ldots,n-2$ are the subsequent perspectivities.  Each of these perspectivities is orientation preserving, since in each case the center of perspectivity lies to the left of the domain and codomain lines (Lemma \ref{aboutperspectivities}). Thus, $f$ is an orientation preserving projectivity.

\emph{Claim: If $y\in  \ell_{p_np_1}$ is such that $\alpha(y)\notin (x_2,p_3)$, then $f(y)\in [p_n,f(x_1)]$.}

Proof of the claim: The center of perspectivity of  $\beta_{k}$ belongs to the segment $(x_{k+1},x_{k+2})$ for all $k$. By  Lemma \ref{pocketperspectivities} (ii), keeping in mind that $\beta_1(p_3)=p_3$, the pole of $\beta_1$ belongs to the segment $(x_2,p_3)$. It follows that $\beta_1$ maps $\ell_{p_2p_3}\backslash (x_2,p_3)$ into $[p_3,x_3]$. By a similar argument, $\beta_{k}$ maps  $[p_{k+1},x_{k+1}]$ into $[p_{k+2},x_{k+2}]$ for $k=2,\ldots, n-2$. So, if  $\alpha(y)$ belongs to $\ell_{p_2p_3}\backslash (x_2,p_3)$,  then continuing to evaluate  on each $\beta_k$ we obtain that
\[
f(y)=(\beta_{n-2}\circ \cdots \circ\beta_{1}\circ\alpha)(y)
\in [p_n,\beta_{n-2}(x_{n-1})]=[p_n,f(x_1)].
\]
This completes the proof of the claim.

If $f$ has no pole, then it is an  affine function on $\ell_{p_np_1}$. In this case, $f(c)<c$ and $f(x_1)\geq x_1$ together  imply that $f(c')\geq c'$. Thus, either $c$ or $c'$ is a good point, and we are done.  

Suppose, then, that $f$ has a pole, and let us denote it by $y_0$. Since $y_0\notin [c,c']$, we have either that $y_0<c$ or $y_0>c'$ in the order of $\ell_{p_np_1}$. Suppose that $y_0> c'$. See Figure \ref{polerightofcprime}. Since  $f$ is a convex function on $(-\infty,y_0)$ (see Section \ref{geoprelims}), we cannot have that $f(c)<c$, $f(x_1)\geq x_1$, and $f(c')<c'$. Thus, either $f(c)\geq c$ or $f(c')\geq c'$, i.e., either $c$ or $c'$ is good, and we are  done.

\begin{figure}[h]

\centering

\begin{tikzpicture}[line width = 1pt,x=1.0cm,y=1.0cm]

\coordinate (s) at (-1, -2);
\coordinate (pn) at (0.5,-2);
\coordinate (p1) at (6.5,-2);
\coordinate (t) at (8, -2);
\coordinate (c) at (1.89,-2);
\coordinate (x1) at (3.33,-2);
\coordinate (c') at (4.56,-2);
\coordinate (y0) at (5.77,-2);

\draw  (pn)-- (p1);
\draw [dashed] (s) -- (pn);
\draw [dashed] (p1) -- (t);
\draw  (pn)-- (-0.5,-0.5);
\draw  (p1)-- (6.9,-0.5);

\foreach \point in {pn, p1, c, x1, c', y0} {
    \draw [fill=black] (\point) circle (1.5pt);
}

\foreach \point/\label in {pn/$p_n$, p1/$p_1$, c/$c$, x1/$x_1$, c'/$c'$, y0/$y_0$}{
\node[below] at (\point) {\label};}

\node[above] at (y0) {pole};
\end{tikzpicture}

\caption{Case when the pole of $f$ is to the right of $c'$.}
\label{polerightofcprime}
\end{figure}
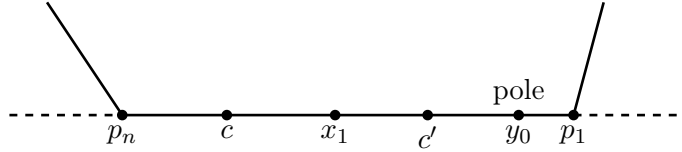

Let us henceforth assume that $f$ has a  pole $y_0$, and that  $y_0<c$ in the order of $\ell_{p_np_1}$.
We split the remainder of the proof into  three cases, according to the relative position of the lines $\ell_{p_np_1}$ and $\ell_{p_2p_3}$.

First, suppose that $\ell_{p_2p_3}$ is parallel to $\ell_{p_np_1}$. In this case $\alpha$ is affine and increasing,  so we can find 
a point $w\in \ell_{p_np_1}$, with $w>x_1$, such that $\alpha(w)=p_3$. By the claim proved above, $f(w)\in[p_n,f(x_1)]$.
 The inequalities $f(w)\leq f(x_1)$ and $x_1< w$ contradict that $f$ is increasing on any segment not containing its pole.  Thus, 
 $\ell_{p_2p_3}$ is not parallel to $\ell_{p_np_1}$ under our present assumptions.

Let $z$ be the intersection point of  $\ell_{p_np_1}$ and $\ell_{p_2p_3}$.  Since $n\geq 4$, $z$ is outside of $P$. Thus, $\alpha(z)=z\notin (x_2,p_3)$. By the claim proved above,  $f(z)\in [p_n,f(x_1)]$.  

Since $z$ is outside of $P$, we either have that $z>p_1$
or $z<p_n$ in the order of $\ell_{p_np_1}$.  If $z>p_1$, then the inequalities $z>x_1$ and $f(z)\leq f(x_1)$ contradict that $f$ is increasing on the segment $[c,z]$, which does not contain its pole. We thus have that  $z<p_n$. See Figure \ref{zlessthanpn}.
It follows that $z$ is a good point, since $f(z)\in[p_n,f(x_1)]$ implies that $z<p_n\leq f(z)$. Observe that  $y_0\notin [z,c]$,  since any  value of $f$ before its pole is larger than any of its values after the pole, but $f(z)\leq f(x_1)$. Thus,  $y_0<z$.

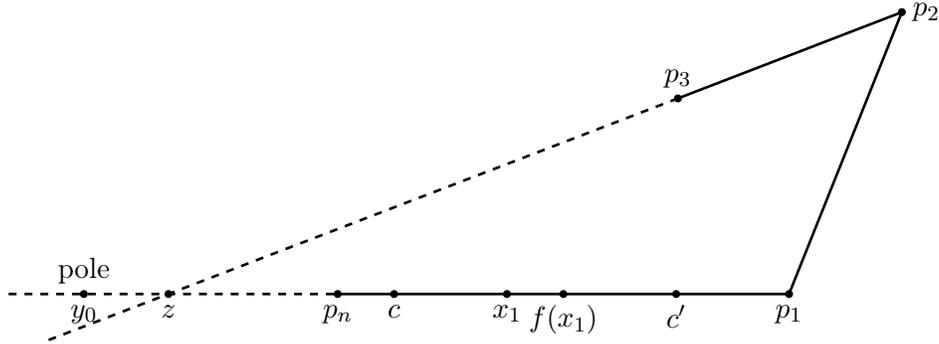
\begin{figure}[h]
\centering

\begin{tikzpicture}[line width = 1pt,x=0.75cm,y=0.75cm]

\coordinate (y0) at (-7.5, -2);
\coordinate (z) at (-6, -2);
\coordinate (pn) at (-3, -2);
\coordinate (c) at (-2, -2);
\coordinate (x1) at (-0, -2);
\coordinate (fx1) at (1, -2);
\coordinate (c') at (3, -2);
\coordinate (p1) at (5, -2);
\coordinate (p3) at (3.03, 1.47);
\coordinate (p2) at (7, 3);

\draw (pn) -- (p1) -- (p2) -- (p3);
\draw[dashed, shorten <= -1cm] (y0) -- (pn);
\draw[dashed, shorten >= -1.8cm] (p3) -- (z);

\foreach \point in {z, p1, p2, p3, pn, c, x1, fx1, c', y0}{
\draw[fill=black] (\point) circle (1pt);}

\foreach \point/\label in {z/$z$, p1/$p_1$, pn/$p_n$, c/$c$, x1/$x_1$, fx1/$f(x_1)$, c'/$c'$, y0/$y_0$}{
\node[below] at (\point) {\label};}

\node[right] at (p2) {$p_2$};
\node[above] at (p3) {$p_3$};
\node[above] at (y0) {pole};
\end{tikzpicture}

\caption{Case $y_0<z<p_n$.}
\label{zlessthanpn}
\end{figure}

 Since $f$ is concave down on $(y_0,+\infty)$, and   $c$ is between good points $z$ and $x_1$, it follows that $c$ is a good point.

We have thus shown that either $c$ or $c'$ is a good point, as desired.
\end{proof}

Let us outline how our results give an effective way of testing for degenerate containment given a pair of $n$-gons $P$ and $P'$, with $P'\precsim P$:
\begin{enumerate}
\item Check whether $P$ is set-convex. If it is not, then $P'$ is degenerately contained in $P$.

\item Check whether the vertices of $P'$ are all collinear. If they are, then $P'$ is degenerately contained in $P$.

\item If $P$ is set-convex, re-index its vertices so that it is convex and oriented counterclockwise. Note that the relation of degenerate containment is
unaffected by this, so we can test whether $P'$ is degenerately contained in the new re-indexed $P$.

\item If the vertices of $P'$ are non-collinear, apply the Broken Line Construction (BLC) starting from points on the set $T$. If, for some such point, the result is an $l$-gon with $l<n$, then $P'$ is degenerately contained in $P$. Otherwise, $P'$ is non-degenerate in $P$.
\end{enumerate}

\section{The Threshold}
In this section we begin the investigation of the Attainability and Bang-Bang problems for $n$-gons  $P'$ non-degenerately contained in a convex $n$-gon $P$.

\begin{definition}
The threshold set associated to a convex polygon $P$, denoted by $\TT_P$, is defined as the set of all polygons  $P'$ attainable from $P$ and such that at least one vertex of $P'$ belongs to $\partial P$. 
\end{definition}

Throughout this section we let $P$ be a convex $n$-gon oriented counterclockwise. In Theorem \ref{thresholdthm} below we characterize the 
polygons in  $\TT_P\backslash \D_P$ (in the threshold and non-degenerate)  by an  explicit criterion involving the BLC.

\begin{lemma}\label{orientation}
Let $\{P(t):c_0\leq t\leq c_1\}$ be a continuous path of set-convex polygons. If $P(c_0)$ is convex oriented counterclockwise, then so is $P(c_1)$.
\end{lemma}

\begin{proof}
 Suppose, for a contradiction, that $P(c_0)$ is convex and oriented counterclockwise but $P(c_1)$ is not. Then, for some pairwise distinct indices 
 $i,i+1,j$,  we have that $p_j(c_0)$ lies strictly to the left of the directed line through $p_i(c_0)$ and $p_{i+1}(c_0)$ while $p_j(c_1)$ does not lie strictly to the left of the directed line through $p_{i}(c_1)$ and $p_{i+1}(c_1)$. By continuity, there exists $t\in[c_0,c_1]$  such that  $p_j(t)$, $p_{i}(t)$, and $p_{i+1}(t)$ are collinear. This violates that $P(t)$ is set-convex.
\end{proof}

\begin{lemma}\label{pionedge}
Let $\{P(t):0\leq t\leq c\}$ be a decreasing path of polygons such that $P(0)=P$ and $P(t)$ is non-degenerate in $P$ for all $t$. Suppose that $p_i(c)\in \partial P$ for some $i$. The following statements hold:
\begin{enumerate}[(i)]
\item
Either $p_i(t)\in (p_{i-1},p_i]$ for all $t\leq c$, or $p_i(t)\in [p_i,p_{i+1})$ for all $t\leq c$.
\item
Let $t_0\in [0,c]$ be the least $t$ such that $p_i(t)$  remains fixed thereafter, i.e.,
\[
t_0=\min\{t\in [0,c]:p_i(t')=p_i(c)\hbox{ for all }t'\in [t,c]\}.
\]
If $p_i(c)\in (p_{i-1},p_i]$, then $p_{i-1}(t_0)\in [p_{i-1},p_i(t_0))$, and if $p_i(c)\in [p_i,p_{i+1})$, then $p_{i+1}(t_0)\in (p_i(t_0),p_{i+1}]$.
\end{enumerate}
\end{lemma}

\begin{proof}
Assume without loss of generality that $i=1$. 

(i) Suppose $p_1(c)$ belongs to the edge $[p_k,p_{k+1}]$ of $P$. Let us show that $p_1(t)\in [p_k,p_{k+1}]$ for all $t\leq c$. Suppose for contradiction that there exists $0\leq c'<c$ such that $p_1(c')\notin [p_k,p_{k+1}]$. Let $c_0$ be the first time after $c'$ that $p_1(t)$ reaches the interval $[p_k,p_{k+1}]$, i.e.,
\[
c_0=\min \{t\in [c',c]: p_1(t)\in [p_k,p_{k+1}]\}.
\]

Consider the sets $I(t) = \co (P(t)) \cap [p_k,p_{k+1}]$ for $0\leq t\leq c_0$. The mapping $t\mapsto I(t)$ is obtained by intersecting a decreasing family of convex sets with a segment. It is thus a decreasing family of segments. Since $p_1(c_0)\in I(c_0)$, these segments are non-empty. Moreover, since $[p_k,p_{k+1}]$ is a face of $\co(P)$, we have 
\[
I(t) = \co (P(t)) \cap [p_k,p_{k+1}] = \co(P(t)\cap [p_k,p_{k+1}]), 
\]
and thus the endpoints of $I(t)$ are vertices of $P(t)$. Since the vertices of $P(t)$ are finitely indexed, there exists an infinite sequence of times $t_1<t_2<\cdots$ with $t_N\to c_0$ and two fixed indices $i,j$ such that $I(t_N)=[p_i(t_N),p_j(t_N)]$ for all $N$. Note that $i,j\neq 1$, since $p_1(t)\notin [p_k,p_{k+1}]$ for $c'\leq t<c_0$. It follows by continuity that $p_i(c_0),p_j(c_0)\in I(c_0)$. If $i=j$, i.e., $I(t_N)$ is a singleton for all $N$, then  $p_i(c_0)=p_1(c_0)$, which contradicts the non-degeneracy of $P(c_0)$ in $P$. If $i\neq j$, then $p_i(c_0),p_j(c_0),p_1(c_0)$ are collinear vertices of $P(c_0)$ with pairwise distinct indices, which again contradicts the non-degeneracy of $P(c_0)$ in $P$. Hence, we conclude that $p_1(t)\in [p_k,p_{k+1}]$ for all $t\leq c$.

Since $p_1=p_1(0)\in [p_{k},p_{k+1}]$, we have either $k=1$ or $k=n$. It is not possible to have $p_1(t)=p_2$ for some $t>0$, as this would entail that $p_1\in [p_2,p_3]$. Similarly, we have $p_1(t)\neq p_n$ for all $t$. This completes the proof of (i).

(ii) Suppose $p_1(c)\in (p_n,p_1]$. Our goal is to show that $p_n(t_0)\in [p_n,p_1(t_0))$.
The case $p_1(c)\in [p_1,p_2)$ can be dealt with similarly. 

If $t_0=0$ the result holds trivially, so assume that $t_0>0$. Note that, by part (i), $p_1(t)\in (p_n,p_1]$ for all $t$. 
Consider the decreasing intervals $I(t)=\co(P(t))\cap [p_n,p_1]$, for $0\leq t\leq c$. As argued in part (i), the endpoints of $I(t)$ are vertices of $P(t)$. Moreover, by (i),  these vertices can only be $p_n(t)$ and $p_1(t)$. Since $p_1(t)\in I(t)$ for all $t$, we either have that $I(t)=[p_n(t),p_1(t)]$ or $I(t)=\{p_1(t)\}$, for each $t$.

By the definition of $t_0$, and since $t_0>0$, there exist $t_1<t_2<\cdots$ such that $t_N\to t_0$ and $p_1(t_N)\neq p_1(t_0)$($=p_1(c)$) for all $N$. Observe  that $I(t_N)$ is not a singleton, as it contains both $p_1(t_N)$ and $p_1(t_0)$. It follows that $I(t_N) = [p_n(t_N), p_1(t_N)]$ for all $N$, and by continuity, that  $I(t_0)=[p_n(t_0), p_1(t_0)]$. This shows that  $p_n(t_0)\in [p_n,p_1]$. Since $p_n(t)\in [p_n,p_1]$ for all $t\leq t_0$, and $p_n(t)\neq p_1(t)$ for all $t$,
we further deduce that $p_n(t_0)\in [p_n,p_1(t_0))$.
\end{proof}

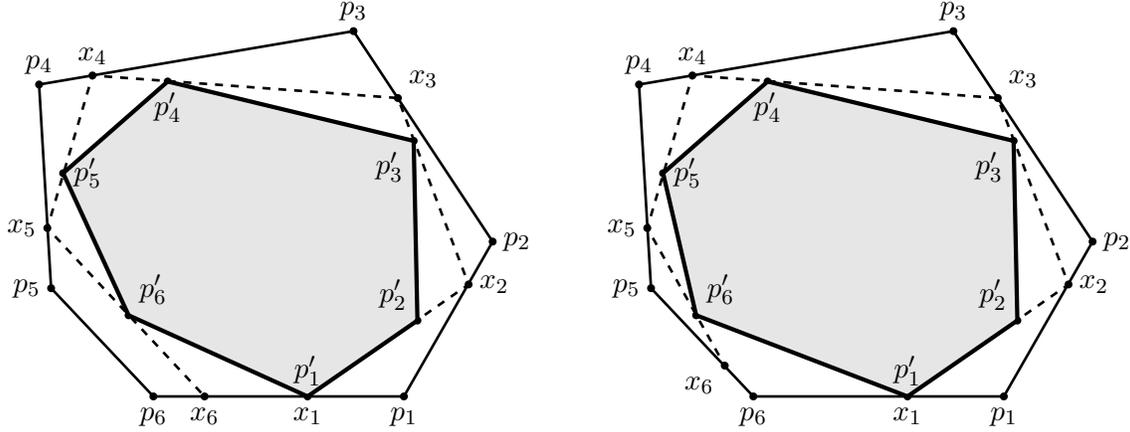
\begin{figure}
\centering

\begin{tikzpicture}[line width = 1pt, x=1.0cm, y=1.0cm]

\coordinate (p6) at (0.66,-1);
\coordinate (p1) at (3.99,-1);
\coordinate (p2) at (5.17,1.06);
\coordinate (p3) at (3.32,3.86);
\coordinate (p4) at (-0.86,3.15);
\coordinate (p5) at (-0.70,0.44);
\coordinate (p1') at (2.71,-1);
\coordinate (p2') at (4.17,0.01);
\coordinate (p3') at (4.12,2.40);
\coordinate (p4') at (0.85,3.19);
\coordinate (p5') at (-0.54,1.97);
\coordinate (p6') at (0.33,0.08);
\coordinate (x2) at (4.85,0.49);
\coordinate (x3) at (3.91,2.97);
\coordinate (x4) at (-0.15,3.27);
\coordinate (x5) at (-0.75,1.24);
\coordinate (x6) at (1.34,-1);

\foreach \point in {p6, p1, p2, p3, p4, p5, p1', p2', p3', p4', p5', p6', x2, x3, x4, x5, x6} {
    \draw[fill=black] (\point) circle (1pt);
}

\draw (p6) -- (p1) -- (p2) -- (p3) -- (p4) -- (p5) -- cycle;
\draw[line width=1.6pt,fill=black!10!white] (p1') -- (p2') -- (p3') -- (p4') -- (p5') -- (p6') -- cycle;
\draw [dashed] (p1') -- (x2) -- (x3) -- (x4) -- (x5) -- (x6);

\foreach \point/\label in {p6/$p_6$, x6/$x_6$, p1'/$x_1$, p1/$p_1$}{
\node[below] at (\point) {\label};
}

\node[right] at (p2) {$p_2$};
\node[above] at (p3) {$p_3$};
\node[above] at (p4) {$p_4$};
\node[left] at (p5) {$p_5$};
\node[above] at (p1') {$p_1'$};
\node[above left] at (p2') {$p_2'$};
\node[below left] at (p3') {$p_3'$};
\node[below] at (p4') {$p_4'$};
\node[right] at (p5') {$p_5'$};
\node[above right] at (p6') {$p_6'$};
\node[right] at (x2) {$x_2$};
\node[above right] at (x3) {$x_3$};
\node[above] at (x4) {$x_4$};
\node[left] at (x5) {$x_5$};
\end{tikzpicture}\qquad
\begin{tikzpicture}[line width = 1pt, x=1.0cm, y=1.0cm]

\coordinate (p6) at (0.66,-1);
\coordinate (p1) at (3.99,-1);
\coordinate (p2) at (5.17,1.06);
\coordinate (p3) at (3.32,3.86);
\coordinate (p4) at (-0.86,3.15);
\coordinate (p5) at (-0.70,0.44);
\coordinate (p1') at (2.71,-1);
\coordinate (p2') at (4.17,0.01);
\coordinate (p3') at (4.12,2.40);
\coordinate (p4') at (0.85,3.19);
\coordinate (p5') at (-0.54,1.97);
\coordinate (p6') at (-0.1,0.08);
\coordinate (x2) at (4.85,0.49);
\coordinate (x3) at (3.91,2.97);
\coordinate (x4) at (-0.15,3.27);
\coordinate (x5) at (-0.75,1.24);
\coordinate (x6) at (0.28,-0.59);

\foreach \point in {p6, p1, p2, p3, p4, p5, p1', p2', p3', p4', p5', p6', x2, x3, x4, x5, x6} {
    \draw[fill=black] (\point) circle (1pt);
}

\draw (p6) -- (p1) -- (p2) -- (p3) -- (p4) -- (p5) -- cycle;
\draw[line width=1.6pt,fill=black!10!white] (p1') -- (p2') -- (p3') -- (p4') -- (p5') -- (p6') -- cycle;
\draw [dashed] (p1') -- (x2) -- (x3) -- (x4) -- (x5) -- (x6);

\foreach \point/\label in {p6/$p_6$, p1'/$x_1$, p1/$p_1$}{
\node[below] at (\point) {\label};
}

\node[below left] at (x6) {$x_6$};
\node[right] at (p2) {$p_2$};
\node[above] at (p3) {$p_3$};
\node[above] at (p4) {$p_4$};
\node[left] at (p5) {$p_5$};
\node[above] at (p1') {$p_1'$};
\node[above left] at (p2') {$p_2'$};
\node[below left] at (p3') {$p_3'$};
\node[below] at (p4') {$p_4'$};
\node[right] at (p5') {$p_5'$};
\node[above right] at (p6') {$p_6'$};
\node[right] at (x2) {$x_2$};
\node[above right] at (x3) {$x_3$};
\node[above] at (x4) {$x_4$};
\node[left] at (x5) {$x_5$};
\end{tikzpicture}

\caption{Left: The BLC starting at $p_1'$ shows that $P'$ is attainable. Right: the BLC shows that $P'$ is unattainable.}

\label{figthreshold}
\end{figure}

\begin{theorem}\label{thresholdthm}
Let $P$ be a convex $n$-gon oriented counterclockwise.  Let $P'$ be an $n$-gon non-degenerately contained in $P$ and such that $p_i'\in \partial P$ for some $1\leq i\leq n$. Then $P'$ is attainable from $P$ if and only if it is convex, oriented counterclockwise, and one of the two following conditions is satisfied:
\begin{enumerate}[(i)]
\item
$p_i'$ belongs to $(p_{i-1},p_i]$ and the BLC starting at $x_1=p_i'$ stops at $x_{n} \in [p_{i-1},p'_i)$.
\item
$p_i'$ belongs to $[p_i,p_{i+1})$ and the clockwise BLC starting at $x_1=p_i'$ stops at $x_{n} \in (p_i',p_{i+1}]$.
\end{enumerate}
Moreover, in either case $P'$ is attainable in $2n-1$ moves. (See Figure \ref{figthreshold}.)
\end{theorem}

\begin{proof}
We assume without loss of generality that $i=1$.

Suppose that $P'$ is convex, oriented counterclockwise, and one of the two conditions of the theorem is satisfied. Let's assume that we are in the first case, as the situation is entirely symmetric in the second case.  Let $(x_k)_{k=1}^n$ be the polygon resulting from the BLC starting at 
$x_1=p_1'$. We explain how to go from $P'$ to $P$ via push-outs: First, push $p_2'$ out onto $x_2$ with $p_1'$. Then push $p_3'$ out onto $x_3$ with $x_2$. Continue in this way until $p_n'$ has been pushed out onto $x_n\in [p_n,p_1')$. Next, push $p_1'$ out with $x_n$ onto $p_1$.
This takes $n$ push-out moves. At this point, $P'$ is inscribed in $P$ and  $p_1'=p_1$. By Lemma \ref{whenmeqn}, we can push out the remaining vertices of $P'$ to occupy all the vertices of $P$ using $n-1$ push-out moves. Thus, after a total of $2n-1$ push-out moves, the vertices of  $P'$ occupy all the vertices of $P$ and $p_1'$ occupies $p_1$. By Lemma \ref{orientation}, $P'$ remains non-degenerate and oriented counterclockwise throughout this process. It follows that $p_i'$ has been  pushed out onto $p_i$ for all $i$.

Let us now prove the converse. Assume that $P'$ is attainable from $P$, non-degenerately contained in $P$, and that $p_1\in \partial P$. Let $\{P(t):0\leq t\leq c\}$ be a path decreasing from $P$ to $P'$. Since $P(t)$ remains set-convex for all $t$, $P'$ is convex and oriented counterclockwise, by Lemma \ref{orientation}.  By Lemma \ref{pionedge} (i), $p_1'$ belongs to either  $(p_n,p_1]$ or $[p_1,p_2)$. Let's assume that we are in the first case, as the situation is entirely symmetric in the second case. 

Let $t_0$ be the smallest number such that $p_1(t) = p_1'$ for $t\in [t_0,c]$. By Lemma \ref{pionedge} (ii), we have $p_n(t_0)\in [p_n,p_1)$. Observe that $P(t_0)$ interpolates between $P'$ and $P$, and its vertices $p_n(t_0)$ and $p_1(t_0)$ belong to the edge $[p_n,p_1]$. We can perform push-outs on $P(t_0)$ to move its remaining vertices onto $\partial P$. The resulting polygon $Q=(q_k)_{k=1}^n$ is inscribed in $P$,
oriented counterclockwise, and interpolates between $P'$ and $P$. Moreover, we have $q_1=p_1'$ and $q_n\in [p_n,p_1')$.

 The vertices of $Q$ divide $\partial P$ into $n$ partition arcs   $\arc[q_{i-1},q_i)$, with $i=1,\ldots,n$.
The BLC (relative to $P'$) starting at $x_1=q_1$  results in an interpolating polygon $(x_k)_{k=1}^l$ with at most $n$ vertices, by Lemma \ref{m+1verticeslemma} (ii), and also with at least $n$ vertices,  since $P'$ is non-degenerate in $P$. Therefore, by Lemma \ref{pockets}, the last point in the construction  $x_n$ belongs to the last partition arc created by $Q$: $\arc[q_n,q_1)\subseteq [p_n, p'_1)$. This shows that condition (i) holds.
\end{proof}

\begin{corollary}\label{compactT}
The threshold set $\TT_P$ is a compact subset of $\PP_n$.
\end{corollary}
\begin{proof}
A polygon in $\TT_P$ is either degenerate in $P$, in which case it is attainable in fewer than  $5n$ pull-in moves, or non-degenerate in $P$, in which case it is attainable in at most $2n-1$ pull-in moves, by the previous theorem. Thus $\TT_P$ is the intersection of the set of $n$-gons with at least one vertex in $\partial P$, a closed set, and the set of $n$-gons attainable from $P$ in at most $5n$ pull-in moves, a compact set by Proposition \ref{compactNmoves}. Thus, $\TT_P$ is compact.
\end{proof}

\begin{corollary}\label{p1'corollary} 
Let  $P'$ be non-degenerate in $P$, attainable from $P$, and such that $p_i'\in (p_{i-1},p_i]$. If we change $p_i'$ anywhere along $[p_i',p_i]$, the resulting $n$-gon is still attainable. 
\end{corollary}

\begin{proof}
Assume without loss of generality that $i=1$. Let $p_1''\in [p_1',p_1]$ and set  
$P''=(p_1'',p_2',\ldots,p_n')$.   If $P''$ is degenerate in $P$, then it is attainable by Theorem \ref{in5nmoves}.  Let us henceforth assume that $P''$ is non-degenerate in $P$. This implies that the BLC  relative to $P''$ starting at $y_1=p_1''\in [p'_1,p_1]$ results in an $n$-gon $(y_k)_{k=1}^n$, and that the pivots used in this construction are $p_2',p_3',\ldots,p_n'$ (the vertices of $P''$ other than $p_1''$). We wish to show that $y_n\in [p_n,y_1)$, as by Theorem  \ref{thresholdthm} this means that $P''$ is attainable.

Observe that, since the pivots used in the the construction of $(y_k)_{k=1}^n$ are  $p_2',p_3',\ldots,p_n'$, the BLC  with starting point $y_1=p_1''$ and relative to the $(n-1)$-gon $\tilde P:=(p_2',p_3',\ldots,p_n')$ results also in the polygon  $(y_k)_{k=1}^n$. 

By Theorem \ref{thresholdthm}, the BLC relative to $P'$ starting at $x=p'_1$ results in an $n$-gon $(x_k)_{k=1}^n$ with $x_n\in [p_n,p'_1)$. These points create $n$ partition arcs $\arc[x_k,x_{k+1})$, $k=1,\ldots,n$, along $\partial P$. By construction, $\tilde P$ lies to the left of   $\ell_{x_kx_{k+1}}$ for all $k$. Thus, we can apply Lemma \ref{pockets} to the BLC relative $\tilde P$ starting at  $y_1=p_1''$
and the aforementioned partition arcs. This
leaves only two possibilities for the location of $y_n$: either $y_n\in [x_1,y_1)$ (the last point in the construction  belongs to the same partition arc as the first), or each partition arc contains exactly one point, in which case $y_n\in [x_n,x_1)$ (the last point belongs to the last partition arc).  In either case, $y_n\in [p_n,y_1)$ as desired.
\end{proof}

\section{The Vestibule}
Throughout this section we fix a convex $n$-gon $P$ oriented counterclockwise.  Recall that we denote by $\D_P$ the set of $n$-gons that are degenerately contained in $P$, and by $\TT_P$ the set of  $n$-gons attainable from $P$ and with at least one vertex in $\partial P$, i.e., the threshold set associated to $P$.

We call \emph{neighbor pull-in move} a pull-in  of a vertex $p_i$ toward one of its two neighboring vertices  $p_{i-1}$ or $p_{i+1}$.

\begin{definition}
For a fixed convex $n$-gon $P$,  we define the vestibule set $\V_P$ as the set of $n$-gons $P'$ obtained 
by applying one single neighbor pull-in move to an $n$-gon in $\TT_P$. (\emph{Note}: 
We allow the parameter of the pull-in move to be 0. So we have $\TT_P\subseteq \V_P$.) 
\end{definition}

In Theorem \ref{attainableVuD} below we show that, for $n\geq 4$, $\V_P\cup \D_P$ is precisely the set of polygons attainable from $P$.
The main step in the proof consists in showing that $\V_P\cup\D_P$ is closed under pull-in moves (Theorem \ref{closedundermoves}). 
Before proving this, we need several preparatory  lemmas.

\begin{lemma}\label{compactV}
The set $\V_P$ is compact. 
\end{lemma}
\begin{proof}
Let us denote by $\C_P$ the set of $n$-gons contained in $P$.
For each $i=1,\ldots,n$, let $f_{i,i+1}\colon \C_P\times [0,1]\to \C_P$ be the function implementing a pull-in of $p_i'$ toward $p_{i+1}'$:
\[
f_{i,i+1}(P',t)=(p_1',\ldots, (1-t)p_i'+tp_{i+1}',\ldots,p_n').
\]
Define similarly $f_{i+1,i}\colon \C_{P}\times [0,1]\to \C_P$ for a pull-in of $p_{i+1}'$ toward $p_i'$. Then 
\[
\V_P=\bigcup_{i=1}^n f_{i,i+1}(\TT_P \times [0,1])\cup f_{i+1,i}(\TT_P \times [0,1]). 
\] 
It follows that $\V_P$ is compact, as $f_{i,i+1}$ and $f_{i+1,i}$ are continuous for all $i$ and $\TT_P$ is compact (Corollary \ref{compactT}).
\end{proof}

The following two lemmas specify which pivots (Lemma \ref{sameBLC}) and edges of $P$ (Lemma \ref{lastwallpockets}) the BLC uses under certain restrictive conditions. 

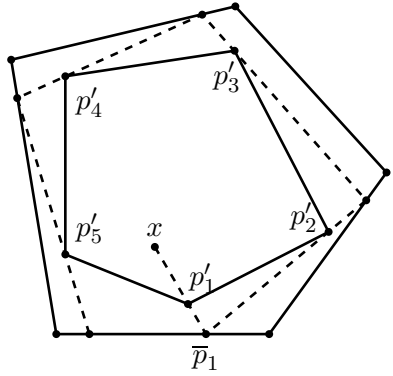
\begin{figure}[h]
\centering
\begin{tikzpicture}[line width = 1pt, x=1.0cm, y=1.0cm]

\coordinate (p1) at (2.62, -0.36);
\coordinate (p2) at (5.45, -0.36);
\coordinate (p3) at (7.01, 1.79);
\coordinate (p4) at (5, 4);
\coordinate (p5) at (2.02, 3.29);
\coordinate (p1p) at (4.37, 0.04);
\coordinate (p2p) at (6.24, 1);
\coordinate (p3p) at (4.99, 3.41);
\coordinate (p4p) at (2.74, 3.07);
\coordinate (p5p) at (2.74, 0.70);
\coordinate (x) at (3.93, 0.80);
\coordinate (p1bar) at (4.61, -0.36);
\coordinate (q2) at (6.74, 1.42);
\coordinate (q3) at (4.56, 3.89);
\coordinate (q4) at (2.1, 2.78);
\coordinate (q5) at (3.06, -0.36);

\draw (p1) -- (p2) -- (p3) -- (p4) -- (p5) -- cycle;
\draw (p1p) -- (p2p) -- (p3p) -- (p4p) -- (p5p) -- cycle;
\draw[dashed] (x) -- (p1bar) -- (q2) -- (q3) -- (q4) -- (q5);

\foreach \point in {p1, p2, p3, p4, p5, x, p1p, p2p, p3p, p4p, p5p, p1bar, q2, q3, q4, q5}
{
\draw[fill = black] (\point) circle (1pt);
}

\node[above] at (x) {$x$};
\node[above, xshift = 2mm] at (p1p) {$p_1'$};
\node[left, yshift = 2mm] at (p2p) {$p_2'$};
\node[below, xshift = -1mm] at (p3p) {$p_3'$};
\node[below right] at (p4p) {$p_4'$};
\node[above right] at (p5p) {$p_5'$};
\node[below] at (p1bar) {$\overline{p}_1$};

\end{tikzpicture}

\caption{$P'$ in $P$ as in Lemma \ref{sameBLC}, with push-out to $\overline p_1$ and BLC starting at $\overline p_1$ denoted by a dashed line.}
\label{figLemma103}
\end{figure}

\begin{lemma}\label{sameBLC}
Let $P'$ be an $n$-gon  whose vertices all lie in the interior of $P$. Suppose that $P'$ is non-degenerately contained in $P$ and oriented counterclockwise.  For a fixed $i$, let $x\in \co(P')\backslash\{p_i'\}$ and let $\overline p_i\in \partial P$ denote the push-out of $p_i'$  with $x$ onto $\partial P$.  Let $\overline P$ denote the polygon resulting from this push-out move, i.e., $\overline P=(p_1',\ldots,p_{i-1}',\overline p_i,p_{i+1}',\ldots,p_n')$.
The following are true:

\begin{enumerate}[(i)]
\item
The polygon resulting from the BLC relative to $\overline P$ with initial point $\overline p_i$ is an $n$-gon, and  the construction uses the pivots  $p_{i+1}',p_{i+2}',\ldots,p_{i-1}'$, in this order.

\item 
The BLCs relative to  $P'$ and $\overline P$, with initial point $\overline p_i$,   result in the same $n$-gon.
\end{enumerate}
(See Figure \ref{figLemma103}.)
\end{lemma}

\begin{proof}
We assume without loss of generality that $i=1$.

(i) The path connecting $P'$ to $\overline P$, obtained from the push-out of $p_1'$ onto $\partial P$ with $x$, consists of set-convex polygons, as $P'$ is non-degenerate in $P$. Thus,  $\overline P$ is convex and  oriented counterclockwise by Lemma \ref{orientation}. 

By the convexity of $\overline P$, the ray $r_{\overline p_1 p_2'}$ is right tangent  to $\overline P$. Thus, the evaluation of $\pi_{\overline P}$ (the Poncelet map relative to $\overline P$) at $\overline p_1$ uses $p_2'$ as a pivot. This is the first pivot used in the BLC relative to $\overline P$ starting at $\overline p_1$. On the other hand, the stopping condition of the BLC prevents $\overline p_1$ from being used as pivot after the first step.  Thus, since $\overline P$ is non-degenerate in $P$, each of its vertices other than $\overline p_1$ must be used as a pivot in the BLC.  Further, by the counterclockwise orientation of $\overline P$,  the pivots  used are precisely $p_2',p_3',\ldots,p_n'$, in this order.

(ii) Given $y\in \partial P$, the evaluations of the Poncelet maps $\pi_{P'}$ and $\pi_{\overline P}$ at $y$ agree as long as a vertex of $P'$ is used as pivot. Thus, the BLC relative to $P'$ starting at $\overline p_1$ also uses the pivots $p_2',p_3',\ldots,p_n'$ in this order. It remains to  show that the BLC relative to $P'$ also stops  after using $p_n'$ as a pivot. It suffices to show  that $p_1'$ is never used as a pivot after the first step of the construction. Suppose that $p_1'\in [y,\pi_{P'}(y)]$ for some $y\in \partial P$ with $y\neq \overline p_1$. Then $x$ is to the left of $\ell_{y\pi_{P'}(y)}$, as it is in $\co(P')$ by assumption. It follows that $\overline p_1$ is to the right of $\ell_{y\pi_{P'}(y)}$. Thus, by the convexity of $P$ and the fact that the chord $[y,\pi_{P'}(y)]$ intersects the interior of $P$,  $\pi_{P'}(y)\notin \arc(y,\overline p_1)$ and the stopping condition of the BLC is triggered. This shows that $p_1'$ is not used as a pivot in the BLC relative to $P'$ after the first step.
\end{proof}

\begin{lemma}\label{lastwallpockets}
Let $P'$ be an $n$-gon  whose vertices all lie in the interior of $P$ and suppose that the BLC starting at $x_1\in \partial P$ results in an $n$-gon.

If $x_1\in  (p_{i-1},p_{i}]$ , but $x_n\notin  (p_{i-1},x_1)$, then
 $x_k\in (p_{k+i-2},p_{k+i-1})$ for $k=1, \ldots, n-1$ and
$x_n\in (p_{i-2},p_{i-1}]$. In particular, we must have that $x_1\neq p_i$.

If $x_1=p_i$, then $x_n\in (p_{i-1},x_1)$.
\end{lemma}

\begin{proof}
We assume without loss of generality that $i=1$.
The assumptions  then take the form $x_1\in (p_n,p_1]$ and $x_n\notin (p_n,p_1]$. 

Since $P'$ is in the interior of $P$, $\pi_{\cw}=\pi^{-1}$. Thus, the clockwise BLC starting at $x_n$ results in the points 
$x_n,x_{n-1},\dots,x_1$. By Theorem \ref{allaboutQ} (the clockwise version), each of these points belongs to a different half-open edge $(p_{k-1},p_{k}]$, except possibly $p_n$ and $p_1$. But, by assumption, we have ruled out the latter possibility. Thus, there is exactly one point $x_k$ in each half-open edge $(p_{k-1},p_{k}]$. Further, by the counterclockwise order of the BLC points,  $x_k\in (p_{k-1},p_{k}]$ for all $k$. Notice now that for $k=1,\ldots,n-1$ we cannot have $x_k=p_k$, since this would mean that the edge  $[p_k,p_{k+1}]$ contains a vertex of $P'$,  namely the pivot $p(x_{k})$ in the evaluation $\pi(x_{k})=x_{k+1}$.
Thus,  $x_k\in (p_{k-1},p_{k})$ for $k=1,\ldots,n-1$ and $x_n\in (p_{n-1},p_n]$, as desired.
\end{proof}

Consider the set $\V_P\cup\D_P$, which is compact in the space of $n$-gons $\PP_n$, as $\V_P$ and $\D_P$ are both compact sets (Lemmas \ref{compactD} and \ref{compactV}).  Let $\partial (\V_P\cup \D_P)$ denote the boundary of $\V_P\cup\D_P$. Theorem \ref{closedundermoves} below states that for $n\geq 4$ the set $\V_P \cup \D_P$ is closed under neighbor pull-in moves. To prove it, it will suffice to focus on polygons in $\partial (\V_P\cup\D_P)$, due to the following simple observation, which we state as a lemma:

\begin{lemma}\label{reducetoboundary}
Suppose that for every ``boundary polygon" $P'\in \partial (\V_P\cup\D_P)$ and path $t\mapsto P'(t)$ generated by a neighbor pull-in move applied to $P'$, there exists $\epsilon>0$ such that $P'(t)\in \V_P\cup \D_P$ for all $0\leq t \leq \epsilon$. Then $\V_P\cup\D_P$ is closed under neighbor pull-in moves.
\end{lemma}

\begin{proof}
Suppose for the sake of contradiction that there exists $P'\in \V_P\cup\D_P$ and a neighbor pull-in move $[0,c]\ni t\to P'(t)$ such that $P'(0)=P'$ and $P'(c)\not\in \V_P\cup \D_P$. 
Let 
\[
t_0=\max \{t\in [0,c]:P'(t)\in \V_P\cup\D_P\}.
\]
Then $P'(t_0)$ is in $\partial(\V_P\cup\D_P)$ and $P(t_0+\epsilon)$ is not in $\V_P\cup\D_P$ for all $\epsilon>0$ where $P(t_0+\epsilon)$ is defined,
contradicting the hypotheses of the lemma. 
\end{proof}

In the next lemma we explore some of the constraints placed on the positioning  of the vertices of $P'$, when 
$P'$  is non-degenerately contained in $P$ and a boundary polygon. These constraints will simplify the number
of cases that we need to consider in  the proof of Theorem \ref{closedundermoves} below.

\begin{lemma}\label{boundaryP'}
Let $P'\in \partial(\V_P\cup\D_P)$ be a boundary polygon non-degenerately contained in $P$ and  whose vertices all lie in the interior of $P$. For a fixed $i$, let $\overline p_i$ be the push-out of $p_i'$ with $p_{i-1}'$ onto $\partial P$, and let $\overline P$ denote the polygon resulting from this move on $P'$. 
Suppose that $\overline P$ belongs to $\TT_P$. The following are true:
\begin{enumerate}[(i)]
\item
 $\overline p_i\in (p_{i-1},p_i)$.
 
 \item
If $Q=(q_k)_{k=1}^n$ denotes the polygon resulting from the  BLC (relative to $P'$) starting at $\overline p_i$, then
$q_{k}\in (p_{k+i-2},p_{k+i-1})$ for $k=1,2,\ldots,n-1$ and $q_n=p_{i-1}$.
\end{enumerate}
\end{lemma}

\begin{proof}
We assume without loss of generality that $i=1$. 

(i) Since $\overline P$ is non-degenerate in $P$ and in $\TT_P$, we  have that $\overline p_1\in \arc(p_{n},p_2)$ by the Threshold Theorem (Theorem \ref{thresholdthm}). 
Thus, either $\overline p_1\in (p_n,p_1)$ or $\overline p_1\in [p_1,p_2)$. Using that $P'$ non-degenerately contained in $P$ \emph{and a boundary polygon}, 
we will rule out that $\overline p_1\in [p_1,p_2)$. 

Suppose, for  contradiction, that $\overline p_1 \in (p_1,p_2)$. Let $R=(r_k)_{k=1}^n$ be the polygon obtained from the clockwise BLC relative to $P'$ starting at $r_1=\overline p_1$.  By the clockwise analog of Lemma \ref{sameBLC} (ii), $R$ is also the  polygon resulting from the clockwise BLC relative to $\overline P$ starting at $\overline p_1$. Since $\overline P\in \TT_P$, we have by  the Threshold Theorem that $r_n\in (\overline p_1,p_2]$. 
If $r_n\neq p_2$, then by the continuous dependence of the vertices of $R$ with respect to $P'$, any small perturbation of the vertices of $P'$ still yields a  polygon $\overline P\in \TT_P$, violating that $P'$ is in $\partial (\V_P\cup\D_P)$. Therefore, we must have that $r_n=p_2$. 

The first pivot used in the construction of $R$ is $p_n'$ (Lemma \ref{sameBLC} applied to the clockwise BLC). Thus, $r_2$ is the push-out of $p_n'$ with $p_1'$ onto $\partial P$. By Lemma \ref{lastwallpockets} (applied to the clockwise BLC), $r_2\in (p_n,p_1)$.  

Consider the polygon
\[
\widetilde P= (p_1',p_2',\ldots,p_{n-1}',r_2),
\] 
obtained from $P'$ by the push-out of $p_n'$ with $p_1'$ onto $\partial P$.
Let us prove that $\widetilde P$ belongs to $\TT_P$. (See Figure \ref{figboundarypoly}.)
Since $r_2=\pi_{\cw}(r_1)$, the clockwise BLC (relative to $P'$) starting at $r_2$ results in an $n$-gon $R'=(r_2,\ldots,r_{n-1},p_2,s)$, where the last point $s=\pi_{\cw}(p_2)$ belongs to  $\arc(r_2,p_2)$. But $s\notin [p_1,p_2]$, since this would imply the existence of a vertex of $P'$ on $[p_1,p_2]$ (pivot in the evaluation $\pi_{\cw}(p_2)=s$). Thus, $s\in (r_2,p_1)$. By the clockwise analogue of Lemma \ref{sameBLC} (ii), $R'$ is also the result of the clockwise BLC relative to $\widetilde P$ starting at $r_2$. Hence $\widetilde P$ is attainable by the Threshold Theorem (Theorem \ref{thresholdthm}). 

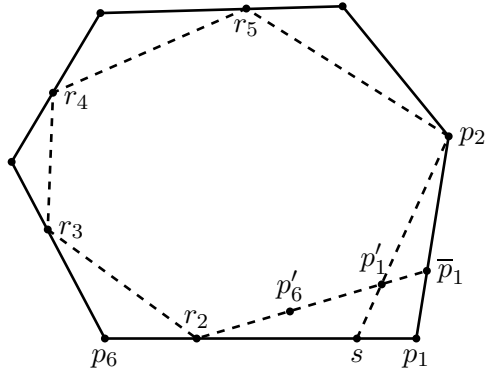
\begin{figure}
\centering

\begin{tikzpicture}[line width = 1pt, x=1.0cm, y=1.0cm]

\coordinate (p6) at (0.61,-1.65);
\coordinate (p1) at (4.75,-1.65);
\coordinate (p2) at (5.18,1.04);
\coordinate (p3) at (3.77,2.77);
\coordinate (p4) at (0.55,2.68);
\coordinate (p5) at (-0.63,0.70);
\coordinate (p1p) at (4.29,-0.93);
\coordinate (p1bar) at (4.89,-0.75);
\coordinate (r2) at (1.83,-1.65);
\coordinate (p6p) at (3.07,-1.29);
\coordinate (r3) at (-0.15,-0.20);
\coordinate (r4) at (-0.08,1.62);
\coordinate (r5) at (2.49,2.74);
\coordinate (s) at (3.96,-1.65);

\foreach \point in {p6, p1, p2, p3, p4, p5, p1p, p1bar, r2, p6p, r3, r4, r5, s} {
    \draw[fill=black] (\point) circle (1pt);
}

\draw  (p6) -- (p1) -- (p2) -- (p3) -- (p4) -- (p5) -- cycle;

\draw[dashed] (s) -- (p2) -- (r5) -- (r4) -- (r3) -- (r2) -- (p1bar); 

\node[below] at (p6) {$p_6$};
\node[below] at (p1) {$p_1$};
\node[below] at (s) {$s$};
\node[right] at (p1bar) {$\overline{p}_1$};
\node[right] at (p2) {$p_2$};

\node[above] at (r2) {$r_2$};
\node[above] at (p6p) {$p_6'$};
\node[above, xshift = -1mm] at (p1p) {$p_1'$};
\node[right] at (r3) {$r_3$};
\node[right, yshift = -1mm] at (r4) {$r_4$};
\node[below] at (r5) {$r_5$};

\end{tikzpicture}

\caption{The clockwise BLC starting at $r_2$ ends at $s\in (r_2,p_1)$, confirming that $\widetilde P$ is attainable.}

\label{figboundarypoly}
\end{figure}

Since $s\in (r_2,p_1)$ and the vertices of $R'$ depend continuously on $P'$, any sufficiently small perturbation of the vertices of $P'$ results in a $\widetilde P\in \TT_P$. This contradicts that $P'$ is in the boundary of $\V_P\cup \D_P$. We conclude that $\overline p_1 \notin (p_1,p_2)$.

We rule out that $\overline p_1=p_1$ by a similar argument:
Suppose that $\overline{p}_1=p_1$. Let $Q$ and $R$ be the polygons resulting from the BLC, and clockwise BLC respectively, relative to $\overline P$ and starting at
$\overline p_1$. These constructions can also be taken relative to $P'$, by Lemma \ref{sameBLC} (ii). It follows by Lemma \ref{lastwallpockets}  (applied to both constructions) that  $q_n\in(p_n,p_1)$ and $r_n\in(p_1,p_2)$. Thus, by the continuous dependence of $Q$ and
 $R$ on $P'$, for any small enough perturbation of the vertices of $P'$ the polygon $\overline P$ is in $\TT_P$ (as confirmed by either the BLC or the clockwise BLC starting at $\overline p_1$). This contradicts that $P'$ is in the boundary of $\V_P\cup \D_P$.

(ii)  By Lemma \ref{sameBLC} (ii), $Q$ is also the polygon resulting from the BLC starting at $q_1=\overline p_1$ relative to $\overline P$. Since $\overline P$ is non-degenerate and in $\TT_P$, we have $q_n\in\ [p_n,\overline p_1)$ by the Threshold Theorem. If $q_n\neq p_{n}$, then by the continuous dependence of the vertices of $Q$ with respect to $P'$, any small perturbation of the vertices of $P'$ results in an attainable polygon $\overline P$, violating that $P'$ is in $\partial (\V_P\cup\D_P)$. Therefore, we must have that $q_n=p_n$. That $q_k\in (p_{k-1},p_k)$ for $k=1,\ldots,n-1$ now follows from Lemma \ref{lastwallpockets}.
\end{proof}

\begin{theorem}\label{closedundermoves}
For $n\geq 4$ the set $\V_P\cup \D_P$ is closed under neighbor pull-in moves.
\end{theorem}

\begin{proof}
By Lemma \ref{reducetoboundary}, it suffices to show that if $P'\in \partial(\V_P\cup \D_P)$, and a neighbor pull-in move is applied to $P'$, then for sufficiently
small values of the parameter of the move the polygons on the path belong to $\V_P\cup\D_P$. If $P'\in \D_P$, then this is indeed the case, as $\mathcal D_P$ is obviously closed under pull-in moves. If $P'\in \mathcal T_P$, then again by the very definition of $\mathcal V_P$, a neighbor pull-in move on $P'$ results in a polygon in $\mathcal V_P$. Let us thus assume that $P'\in \V_P \setminus (\mathcal T_P \cup \mathcal D_P)$. That is, $P'$ is non-degenerate in $P$, all its vertices lie in the interior of $P$, and there exists a neighbor push-out move on $P'$ which produces a polygon in $\mathcal T_P$. 

Without loss of generality, we can assume that the push-out of vertex $p_1'$ with $p_n'$ onto $\overline p_1\in\partial P$ results in an attainable polygon $\overline P$. By Lemma \ref{boundaryP'}, we have that $\overline p_1\in (p_n,p_1)$ and  that the BLC relative to $P'$ starting at $\overline p_1$ ends at $p_n$.

Let us apply a neighbor pull-in move on $P'$. Let $P'(t)=(p_1'(t),\ldots,p_n'(t))$,  for $0\leq t\leq c$, 
denote the path generated by this move, where $P'(0)=P'$. Let $\overline P(t)$ denote the polygon resulting from applying on $P'(t)$ the push-out of $p_1'(t)$ with $p_n'(t)$ onto $\partial P$. Since $\overline p_1\in (p_n,p_1)$, for small enough $t$ we have $\overline p_1(t)\in (p_n,p_1)$. Since our goal is to show that $P'(t)\in \V_P\cup \D_P$ for small enough $t$, we may assume without loss of generality that $\overline p_1(t)\in (p_n,p_1)$ for all $0\leq t\leq c$. Further, since $\D_P$ is a closed set and
$P'\notin \D_P$, we may also assume that $P'(t)\notin  \D_P$ for all $0\leq t\leq c$, i.e., $P'(t)$ remains non-degenerate along the path, whence also convex and oriented counterclockwise, by Lemma \ref{orientation}.

We consider several cases.

\textbf{Case 0: pull-in moves not changing the line $\ell_{p_n'p_1'}$}. More concretely, we mean by this any of the following pull-in moves:
\begin{enumerate}
\item a pull-in  of $p_i'$ toward any vertex of $P'$, for $i\neq 1$ and $i\neq n$,
\item a pull-in  of $p_1'$ toward $p_n'$, 
\item a pull-in  of $p_n'$ toward $p_1'$. 
\end{enumerate}
In all these cases the directed line $\ell_{p_n'p_1'}$ does not change. Thus, the push-out $\overline p_1$ remains constant, i.e., $\overline p_1(t)=\overline p_1$ for all $t$.  In these cases the polygon $\overline P(t)$ is in $\TT_P$ for all $t\geq 0$. Indeed, $\overline P(t)$ is attained from $\overline P$ by the same pull-in move being applied on $P'$ (see Lemma \ref{commutingmoves}.). It follows that $P'(t)\in \V_P$ for all $0\leq t\leq c$, as desired.

The only moves not covered by the above case are the pull-ins of $p_n'$ toward $p_{n-1}'$ and of $p_1'$ toward $p_2'$.

\textbf{Case 1: pull-in of $p_1'$ toward $p_2'$}.  Observe that $p_n'$ lies to the left of $\ell_{p_1'p_2'}$, by the counterclockwise orientation of $P$, and to the left of $\ell_{p_np_1}$, since it is an interior point of $P$. Thus, the perspectivity with center $p_n'$ from $\ell_{p_1'p_2'}$ to $\ell_{p_np_1}$ is locally increasing by Lemma \ref{aboutperspectivities}. Since $p_1'(t)\geq p_1'$, in the order of the line $\ell_{p_1'p_2'}$, we have   $\overline p_1(t)\geq \overline p_1$ in   $\ell_{p_np_1}$ for small enough $t>0$. Thus, $\overline p_1(t)\in [\overline p_1,p_1]$ for all $0\leq t < c'$ and some $c'>0$. It now follows from Corollary \ref{p1'corollary} that 
\[
\overline P(t)=(\overline p_1(t),p_2',\ldots,p_n')
\] 
is in $\TT_P$ for all such $t$. Thus,  $P'(t)=(p_1'(t),p_2',\ldots,p_n')$ belongs to $\V_P$ for all small enough  $t>0$, as desired. See Figure \ref{figCase1}.

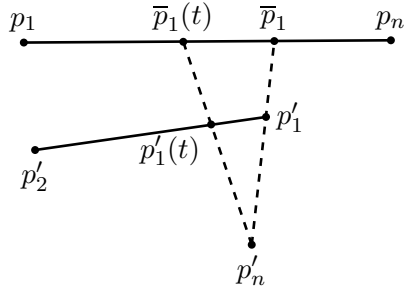
\begin{figure}[h]
\centering

\begin{tikzpicture}[line width = 1pt,x=1.0cm,y=1.0cm]

\coordinate (p_n) at (7.05,3.52);
\draw [fill=black] (p_n) circle (1pt);
\node[above] at (p_n) {$p_n$};

\coordinate (p_1) at (2.17,3.49);
\draw [fill=black] (p_1) circle (1pt);
\node[above]  at  (p_1) {$p_1$};

\coordinate (p1bar) at (5.5, 3.51);
\draw[fill=black] (p1bar) circle (1pt);
\node[above] at (p1bar) {$\overline{p}_1$};

\coordinate (p_n') at (5.2,0.8);
\draw [fill=black] (p_n') circle (1pt);
\node[below] at (p_n') {$p_n'$};

\coordinate (p_1') at (5.39,2.5);
\draw [fill=black] (p_1') circle (1pt);
\node[right] at (p_1') {$p_1'$};

\coordinate (p_2') at (2.32,2.06);
\draw [fill=black] (p_2') circle (1pt);
\node[below] at (p_2') {$p_2'$};

\coordinate (p1pt) at (4.66,2.4);
\draw [fill=black] (p1pt) circle (1pt);
\node[below left] at (p1pt) {$p_1'(t)$};

\coordinate (p1bart) at (4.29, 3.5);
\draw[fill=black] (p1bart) circle (1pt);
\node[above] at (p1bart) {$\overline{p}_1(t)$};

\draw (p_n) -- (p_1);
\draw[dashed] (p_n')-- (p1bar);
\draw (p_1')-- (p_2');
\draw[dashed] (p_n')-- (p1bart);

\end{tikzpicture}

\caption{Case 1: pull-in of $p_1'$ toward $p_2'$.}
\label{figCase1}
\end{figure}

\textbf{Case 2: pull-in of $p_n'$ toward $p_{n-1}'$}. The analysis of this case takes up the main bulk of the proof. Let us introduce some notation.

We denote by $\ell$ the line $\ell_{p_{n-1}'p_n'}$. We denote by $\overline p_n$ the push-out of $p_n'$ with $p_{n-1}'$ onto $\partial P$, and by $\overline p_{n-1}$ the  push-out of $p_{n-1}'$ with $p_n'$ onto $\partial P$. Note that $\overline p_n$ and $\overline p_{n-1}$ are the two points of intersections of $\ell$ with $\partial P$.

Let $Q=(q_k)_{k=1}^n$ denote the polygon resulting from the BLC (relative to $P'$) with starting point $q_1=\overline p_1$. We have by Lemma \ref{sameBLC}  that the pivots used in this construction are $p_2',p_3',\ldots,p_n'$, in this order. Since $P'$ is in $\partial(\mathcal V_P\cup\D_P)$, we also have by Lemma \ref{boundaryP'} that $q_k\in (p_{k-1},p_k)$ for all $k=1,\ldots,n-1$ and $q_n=p_n$.

\emph{Claim 1}: $\overline p_n\in (p_{n-2} ,p_{n-1}]\cup [p_{n-1},p_n)$ and $\overline p_{n-1}\in (p_{n-3},p_{n-2}]\cup [p_{n-2},p_{n-1})$. 

This can be thought of as a limit on the amount of ``twist" of $P'$ relative to $P$. See Figure \ref{figCase2Claim1}.

\begin{figure}[h]
\centering

\begin{tikzpicture}[line width = 1pt, ,x=1.0cm,y=1.0cm]

\coordinate (pn1) at (7.84, -0.01);
\coordinate (qn) at (8.69, 4.46);
\coordinate (pn2) at (2.63, 0.01);
\coordinate (pn3) at (0.91, 3.59);
\coordinate (qn1) at (5.81, 0);
\coordinate (pnp) at (6.51, 1.09);
\coordinate (pnpn1) at (3.41, 1.28);
\coordinate (qn2) at (1.55, 2.27);
\coordinate (pbar) at (8.03, 1);
\coordinate (pnpbar) at (1.98, 1.36);

\draw   (pn1)-- (qn);
\draw   (pn1)-- (pn2);
\draw   (pn2)-- (pn3);
\draw [domain=-0.04283695856100411:10.39517152478831] plot(\x,{(--4.605722461160488-0.18560853808428734*\x)/3.106470664337245});
\draw [dashed] (qn1)-- (qn);
\draw [dashed] (qn1)-- (qn2);
\draw [dashed] (qn2)-- (1.96, 4.19);

\foreach \point/\label in {pn1/p_{n-1}, qn/p_n, qn1/q_{n-1}, pnp/p_n', pbar/\overline{p}_n} {
    \draw [fill=black] (\point) circle (1pt);
    \node [anchor=north west, inner sep=1pt] at (\point) {$\label$};
}

\draw [fill=black] (pnpn1) circle (1pt);
\node [anchor=north, inner sep=1pt] at (pnpn1) {$p_{n-1}'$};

\foreach \point/\label in {pn3/$p_{n-3}$, qn2/$q_{n-2}$, pn2/$p_{n-2}$, pnpbar/$\overline{p}_{n-1}$}{
\node [anchor=north east, inner sep=1pt] at (\point) {\label};
 \draw [fill=black] (\point) circle (1pt);
}

\node at (10, 1.2) {$\ell$}; 
\end{tikzpicture}

\caption{Case 2, Claim 1.}
\label{figCase2Claim1}
\end{figure}
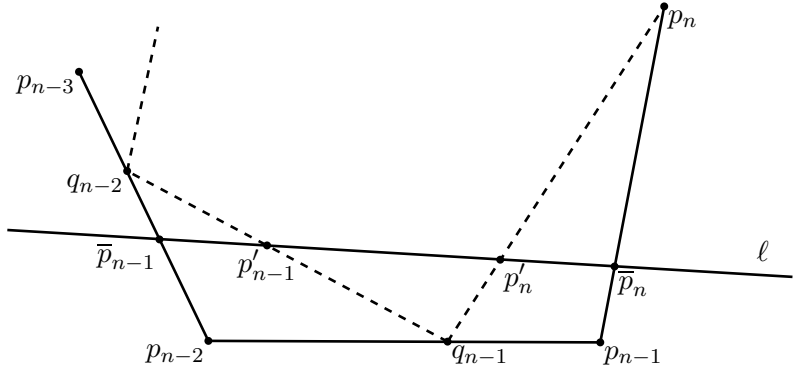

{\it Proof of Claim 1.} Since $p_{n-1}'$ lies strictly to the left of $\ell_{q_{n-1}q_n}$, and $p_n'$ belongs to $\ell_{q_{n-1}q_n}$, the push-out $\overline p_n$ of $p_n'$ with  $p_{n-1}'$ lies strictly to the right of  $\ell_{q_{n-1}q_n}$. It thus belongs to the counterclockwise arc $\arc(q_{n-1},q_n)$ of $\partial P$. Since $q_{n-1}\in (p_{n-2},p_{n-1})$ and $q_n=p_n$, we obtain that  $\overline p_n\in (p_{n-2} ,p_{n-1}]\cup [p_{n-1},p_n)$. By a similar argument,  $\overline p_{n-1}$ must belong to $\arc(q_{n-2},q_{n-1})$, and so it belongs to $(p_{n-3},p_{n-2}]\cup\ [p_{n-2},p_{n-1})$. This proves the claim.

We note that, since $n\geq 4$, the claim  just proved implies that the points $p_n,p_1,\ldots,p_{n-3}$ all lie strictly  to the left of the line $\ell$.

Let us define perspectivities between directed lines as follows: 
\begin{align*}
\alpha &=\ell\stackrel{p_1'}{\wedge}\ell_{p_np_1},\\
\beta_k &=\ell_{p_{k-1}p_k}\stackrel{p_{k+1}'}{\wedge}\ell_{p_kp_{k+1}},\hbox{ for }k=1,\ldots,n-2,\\
\gamma &=\ell_{p_{n-2}p_{n-1}}\stackrel{p_n}{\wedge}\ell.
\end{align*}
(Recall that $\ell_1\stackrel{o}{\wedge} \ell_2$ denotes the perspectivity from $\ell_1$ to $\ell_2$ with center $o$.) The center of  each of these perspectivities lies to the left of both the domain and codomain line. Indeed, in the case of $\beta_1,\ldots,\beta_{n-2}$,   the lines are created by oriented edges of $P$, and the center of perspectivity lies in the interior of $P$. In the case of $\alpha$, the center of perspectivity, $p_1'$, lies to the left of $\ell$ by the convexity and counterclockwise orientation of $P'$, and to the left of $\ell_{p_np_1}$, since it is in the interior of $P$. As for $\gamma$, we have already argued above that $p_n$ lies to the left of $\ell$, and it lies to the left of $\ell_{p_{n-2}p_{n-1}}$ by the convexity and counterclockwise orientation of $P$. It follows from  Lemma \ref{aboutperspectivities} that all these perspectivities are orientation preserving. 

Let $g\colon \ell\to \ell$ denote the projectivity resulting from the  composition of the perspectivities introduced above:
\[
g=\gamma \circ \beta_{n-2}\cdots \beta_1\circ\alpha.
\]
Notice that $g$ is orientation preserving. We have that $\alpha(p_n') =\overline p_1=q_1$. Further, by the pivots used in the construction of $Q$, and the 
edges where its vertices lie, we have that 
\[
\beta_k(q_k) =q_{k+1}\hbox{ for }k=1,\ldots,n-2,
\]
and $\gamma(q_{n-1}) =p_n'$.  Hence,  $g(p_n')=p_n'$.

Let us call a point $x\in \ell$  \emph{good} if $g(x)\geq x$. Observe that $p_n'$ is a good point.

\emph{Claim 2}: If $p_n'(t)$ is good for some $0\leq  t\leq c$, then $\overline P(t)$ belongs to $\TT_P$.

{\it Proof of Claim 2}: We will prove the contrapositive: if $\overline P(t)$ is not in $\TT_P$,
then $p_n'(t)$ is not a good point. Suppose that $\overline P(t)\notin \TT_P$. Then, by the Threshold Theorem,  the BLC relative to $\overline P(t)$ with starting point 
\[
q_1(t):=\overline p_1(t)\in (p_n,p_1)\]
results in an $n$-gon $Q(t)=(q_k(t))_{k=1}^n$ such that $q_n(t)\notin [p_n,p_1)$. By Lemma \ref{sameBLC}, the pivots in the construction of $Q(t)$ are $p_2',p_3',\ldots,p_{n-1}'$ and $p_n'(t)$. Moreover, by Lemma \ref{lastwallpockets}, $q_k(t)$ belongs to $(p_{k-1},p_k)$ for $k=1,\ldots,n$. This means that the perspectivities used in the construction of $q_2(t),\ldots,q_{n-1}(t)$ are $\beta_1,\ldots,\beta_{n-2}$, respectively. Hence, 
\begin{align*}
\beta_k(q_{k}(t)) &=q_{k+1}(t)\hbox{ for $k=1,\ldots,n-2$,}\\
g(p_n'(t)) &=\gamma(q_{n-1}(t)).
\end{align*}

The configuration of the points $q_{n-1}(t)$, $q_n(t)$, and $p_n$, and the lines $\ell$ and $\ell_{p_{n-1}p_n}$ is shown in Figure \ref{VestCase2Claim2}.  

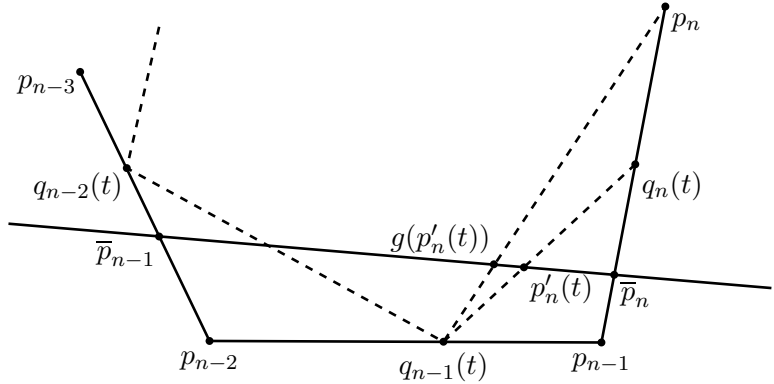
\begin{figure}[h]
\centering

\begin{tikzpicture}[line width = 1pt, ,x=1.0cm,y=1.0cm]

\coordinate (pn1) at (7.84, -0.01);
\coordinate (qn) at (8.69, 4.46);
\coordinate (pn2) at (2.63, 0.01);
\coordinate (pn3) at (0.91, 3.59);
\coordinate (qn1) at (5.74, 0);
\coordinate (pnpn1) at (3.41, 1.28);
\coordinate (qn2) at (1.53, 2.31);
\coordinate (pbar) at (8.01, 0.89);
\coordinate (pnpbar) at (1.96, 1.4);
\coordinate (qnt) at (8.29, 2.36);
\coordinate (gpnt) at (6.41, 1.03);
\coordinate (pnt) at (6.81, 0.99);

\draw (pn3) -- (pn2) -- (pn1) -- (qn) ;

\draw [domain=-0.04283695856100411:10.0999867022187] plot(\x,{(--4.713809166311734-0.25406203614764844*\x)/3.0086586755968594});

\draw [dashed] (1.96, 4.19) -- (qn2) -- (qn1) -- (qn);
\draw [dashed] (qn1)-- (qnt);

\foreach \point/\label in {qn/p_n, pbar/\overline{p}_n, qnt/q_n(t), pnt/p_n'(t)} {
    \draw [fill=black] (\point) circle (1pt);
    \node [anchor=north west, inner sep=2pt] at (\point) {$\label$};
}

\draw [fill=black] (gpnt) circle (1pt);
\node[anchor = south east, xshift = 1mm] at (gpnt) {$g(p_n'(t))$};

\foreach \point/\label in {pn2/p_{n-2}, qn1/q_{n-1}(t), pn1/p_{n-1}} {
    \draw [fill=black] (\point) circle (1pt);
    \node [anchor=north] at (\point) {$\label$};
}

\foreach \point/\label in {pn3/p_{n-3}, qn2/q_{n-2}(t), pnpbar/\overline{p}_{n-1}} {
    \draw [fill=black] (\point) circle (1pt);
    \node [anchor=north east, inner sep=1pt] at (\point) {$\label$};
}

\end{tikzpicture}

\caption{Case 2, Claim 2. If $\overline P(t)\notin \TT_P$, then $g(p_n'(t))<p_n'(t)$ in the order of  $\ell$.}
\label{VestCase2Claim2}
\end{figure}

We remark on the following properties of this figure:  $q_{n-1}(t)$ lies to the left of $\ell_{p_{n-1}p_n}$, as it is a point in $\co(P)$, and to the right of  $\ell$, since, by the BLC, $p_{n}'(t)$ lies to the left of $\ell_{p_{n-1}'q_{n-1}(t)}$.  On the other hand,   $q_n(t)$ lies to the left of $\ell$, since we cross $\ell$ when going from $q_{n-1}(t)$ to $q_n(t)$. As previously observed, $p_n$ also lies to the left of $\ell$. It follows that the perspectivity with center $q_{n-1}(t)$ from $\ell_{p_{n-1}p_n}$ to $\ell$ is decreasing on the segment 
$[q_n(t),p_n]$ and this segment does not contain the pole of the perspectivity (by Lemmas \ref{aboutperspectivities} and \ref{pocketperspectivities} (i)). Since $q_{n}(t)< p_n$ in the order of  $\ell_{p_{n-1}p_n}$,  
\[
p_n'(t)>\gamma(q_{n-1}(t))=g(p_n'(t))
\] 
in the order of  $\ell$. This contradicts that $p_n'(t)$ is a good point. The claim is thus proved.

Our goal next is to show that $p_n'(t)$ remains a good point for all $t\in [0,c]$. By the claim just proven, this will show that $\overline P(t)\in \TT_P$ for all $t\in [0,c]$, which in turn will show  that the polygons $P'(t)$ in the path  coming from the pull-in of $p'_n$ toward $p'_{n-1}$ remain in $\V_P$ for all $t\in [0,c]$, thus completing the proof of the theorem.

Recall that we have argued in Claim 1 above that either $\overline p_n\in (p_{n-2},p_{n-1}]$ or $\overline p_n\in (p_{n-1},p_n)$. We split the analysis into these two subcases, starting with the latter one.

\textbf{Subcase $\overline p_n\in (p_{n-1},p_n)$}: If the polygon $P^{\sharp}=(p_1',p_2',\ldots,p_{n-1}',\overline p_n)$ is attainable, 
then this means that the push-out of $p_n'$ with $p_{n-1}'$ onto $\partial P$ results in a polygon in $\TT_P$. As 
$p_n'$ moves towards $p_{n-1}'$, the push-out of $p_n'(t)$ with $p_{n-1}'$ onto $\partial P$ remains equal to 
$\overline p_n$. (This is analogous to Case 0 (2).) Thus, $P'(t)$ remains in $\V_P$ for all $t$, and we are done. Let us thus assume that $P^\sharp$ is not attainable.

Let $R=(r_k)_{k=1}^n$ denote the polygon resulting from the BLC relative to $P^\sharp$ with starting point $r_1=\overline p_n$. By Lemma \ref{sameBLC}, $R$ also results from the BLC relative to $P'$ starting at $\overline p_n$,  and the construction has pivots $p_1',p_2',\ldots,p_{n-1}'$. Since we have assumed that $P^\sharp$ is not attainable, $r_n\notin [p_{n-1},p_n]$, by the Threshold Theorem. It follows by Lemma \ref{lastwallpockets} that 
$r_k\in (p_{k-2},p_{k-1})$ for $k=1,\ldots,n$. Thus, the perspectivities   used in the construction of $R$ are precisely $\alpha,\beta_1,\ldots,\beta_{n-2}$, in this order. Hence,
\begin{align*}
r_2 &=\alpha(\overline p_n)\\
r_{k+1} &=\beta_{k-1}(r_k)\hbox{ for $k=2,\ldots,n-1$,}\\
g(\overline p_n) &=\gamma(r_n)
\end{align*}  

Since, by the BLC, $p_n'$ lies to strictly to the left of $\ell_{p_{n-1}'r_n}$,   $r_n$ lies strictly to the right of $\ell$. On the other hand, $p_n$ lies strictly to the left of $\ell$.  It follows that the segment $[r_n,p_n]$ intersects $\ell$, and that the intersection point, $\gamma(r_n)$, belongs to the interior of $P$. This shows that 
\[
g(\overline p_n)=\gamma(r_n)\in (\overline p_{n-1},\overline p_n).
\]
In particular,  $\overline p_n$ is not a good point.

\emph{Claim 3}: If $y\in \ell$ is such that $\alpha(y)\notin (r_2, p_1)$, then $g(y)\in [\overline p_{n-1},g(\overline p_n)]$.

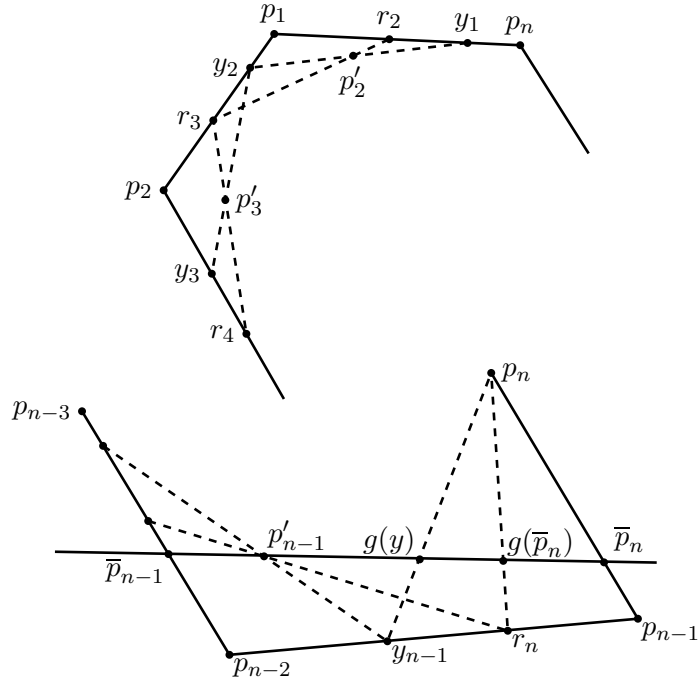
\begin{figure}[h]
\centering

\begin{tikzpicture}[line width = 1pt, x=1.0cm,y=1.0cm]
\useasboundingbox (0,0) rectangle (6, 6);

\coordinate (p2) at (1.06, 3.13);
\coordinate (p1) at (2.53, 5.21);
\coordinate (pn) at (5.80, 5.06);
\coordinate (r2) at (4.06, 5.14);
\coordinate (r3) at (1.72, 4.06);
\coordinate (y1) at (5.10, 5.09);
\coordinate (y2) at (2.21, 4.76);
\coordinate (r4) at (2.16, 1.22);
\coordinate (y3) at (1.70, 2.02);
\coordinate (p2p) at (3.58, 4.92);
\coordinate (p3p) at (1.88, 3.00);

\draw  (p2)-- (p1) -- (pn);
\draw [shorten >= -1cm] (p2)-- (r4);
\draw  (pn)-- (6.71, 3.62);
\draw [dashed] (r2)-- (r3) -- (r4);
\draw [dashed] (y1)-- (y2) -- (y3);

\foreach \point/\label in {p2/p_2, r4/r_4, y3/y_3} {
    \draw [fill=black] (\point) circle (1pt);
    \node [left] at (\point) {$\label$};
}

\foreach \point/\label in {y2/y_2, r3/r_3} {
    \draw [fill=black] (\point) circle (1pt);
    \node [left] at (\point) {$\label$};
}

\foreach \point/\label in {p1/p_1, pn/p_n, r2/r_2, y1/y_1} {
    \draw [fill=black] (\point) circle (1pt);
    \node [above] at (\point) {$\label$};
}

\draw [fill=black] (p2p) circle (1pt);
\draw [fill=black] (p3p) circle (1pt);
\node[below] at (p2p) {$p_2'$};
\node[right] at (p3p) {$p_3'$};
\end{tikzpicture}

\caption{Claim 3: Once $y_1\in [p_n,r_2]$, we deduce that $y_k\in [p_{k-1},r_{k+1}]$ for  $k=1,\ldots,n-1$.}
\label{Case2Claim3Afig}
\end{figure}

{\it Proof of Claim 3}: Set $y_1=\alpha(y)$ and $y_{k+1}=\beta_{k}(y_{k})$ for $k=1,\ldots,n-2$. Note that $g(y)=\gamma(y_{n-1})$. By Lemma \ref{pocketperspectivities} (ii) applied to the perspectivity $\beta_1$, the pole of $\beta_1$ belongs to the segment $(r_2,p_1)$. We deduce
that $\beta_1$ maps $\ell_{p_np_1}\backslash (r_2,p_1)$ into $[p_1,r_3]$. Since $y_1\notin (r_2, p_1)$,   we get that $y_2=\beta_1(y_1)\in [p_1,r_3]$ (see Figure \ref{Case2Claim3Afig}). A similar argument allows us to deduce from $y_2\notin (r_3,p_2)$ that $y_3\in [p_2,r_4]$, and inductively, that $y_k\in [p_{k-1},r_{k+1}]$ for  $k=1,\ldots,n-1$. In particular, $y_{n-1}\in [p_{n-2},r_n]$. 

Let us argue that $y_{n-1}$ lies to the right of $\ell$. If 
$p_{n-2}$ lies to the right of $\ell$, this follows from the fact that $y_{n-1}\in [p_{n-2},r_n]$, and that $p_{n-2}$ and $r_n$
both lie to the right of $\ell$. If, on the other hand, $p_{n-2}$ lies to the left of $\ell$, then again  $y_{n-1}$ lies to the right of $\ell$,
 since we cross $\ell$ when going from $y_{n-2}$ to $y_{n-1}$, and $y_{n-2}\in [p_{n-3},p_{n-2}]$ lies to the left of $\ell$.   It follows that the segment $[y_{n-1},p_n]$ intersects $\ell$. This intersection, which agrees with   $\gamma(y_{n-1})=g(y)$,  thus belongs to  $\co(P)$. Hence, $g(y)\in [\overline p_{n-1},\overline p_n]$. The perspectivity $\gamma$ is increasing on the interval $[y_{n-1},r_n]$, as this interval does not contain its pole (Lemma \ref{pocketperspectivities} (i)).
We thus get that  
\[
g(y)=\gamma(y_{n-1})\leq \gamma(r_n)=g(\overline p_n)
\] 
in the order of $\ell$. Thus,  $g(y)\in [\overline p_{n-1},g(\overline p_n)]$, completing the proof of the claim. 
See Figure \ref{Case2Claim3Bfig}.

\begin{figure}[h]
\centering

\begin{tikzpicture}[line width = 1pt, x=1.0cm,y=1.0cm]

\coordinate (pn) at (5.80, 5.06);
\coordinate (pn1) at (7.75, 1.79);
\coordinate (pn2) at (2.32, 1.31);
\coordinate (pn3) at (0.36, 4.55);

\coordinate (pnbar) at (7.30, 2.54);
\coordinate (pn1p) at (2.78, 2.62);
\coordinate (pn1bar) at (1.51, 2.65);
\coordinate (rn1) at (1.24, 3.09);
\coordinate (yn2) at (0.64, 4.09);
\coordinate (yn1) at (4.42, 1.49);
\coordinate (rn) at (6.02, 1.63);

\coordinate (gy) at (4.85, 2.58);
\coordinate (gpbar) at (5.96, 2.56);

\draw  (pn)-- (pn1) -- (pn2) -- (pn3);
\draw [domain=0:8] plot(\x,{(10.43-0.07*\x)/3.89});
\draw [dashed] (rn1)-- (rn) -- (pn);
\draw [dashed] (yn2)-- (yn1) -- (pn);

\foreach \point/\label in {pn1/p_{n-1}, pn2/p_{n-2},  
  yn1/y_{n-1}, rn/r_n} {
    \draw [fill=black] (\point) circle (1pt);
    \node [anchor=north west, inner sep=1pt] at (\point) {$\label$};
}

\draw[fill = black] (rn1) circle (1pt); 
\draw[fill = black] (yn2) circle (1pt); 

\draw[fill = black] (gy) circle (1pt);
\draw[fill = black] (gpbar) circle (1pt);

\node[anchor=south east, inner sep=1pt] at (gy) {$g(y)$};
\node[anchor=south west, inner sep=1pt] at (gpbar) {$g(\overline p_n)$};

\draw[fill = black] (pnbar) circle (1pt);
\node[above right] at (pnbar) {$\overline p_n$};

\draw[fill = black] (pn) circle (1pt);
\node[right] at (pn) {$p_n$};

\draw[fill = black] (pn1p) circle (1pt);
\node[anchor = south west, inner sep = 1pt] at (pn1p) {$p_{n-1}'$};

\draw[fill = black] (pn1bar) circle (1pt);
\node[anchor = north east, inner sep = 1pt] at (pn1bar) {$\overline p_{n-1}$};

\draw[fill = black] (pn3) circle (1pt);
\node[left] at (pn3) {$p_{n-3}$};
\end{tikzpicture}

\caption{}
\label{Case2Claim3Bfig}
\end{figure}

We split the remainder of this subcases' proof into three cases according to the relative position of the lines
$\ell$ and $\ell_{p_{n}p_1}$.

Suppose that $\ell$ and $\ell_{p_{n}p_1}$ intersect at a point $z$ in the ray $(\overline p_n,+\infty)\subset \ell$.  (See Figure \ref{Zcasefig}.)

\begin{figure}[h]
	\centering
\begin{tikzpicture}[line width = 1pt, x=1.0cm,y=1.0cm]

\coordinate (pn) at (5.80, 5.06);
\coordinate (pn1) at (7.75, 1.79);
\coordinate (pn2) at (2.32, 1.31);
\coordinate (pn3) at (0.36, 4.55);

\coordinate (pnbar) at (7.30, 2.54);
\coordinate (pn1p) at (2.78, 2.62);
\coordinate (pn1bar) at (1.51, 2.65);
\coordinate (rn1) at (1.24, 3.09);
\coordinate (yn2) at (0.64, 4.09);
\coordinate (yn1) at (4.42, 1.49);
\coordinate (rn) at (6.02, 1.63);

\coordinate (gy) at (4.85, 2.58);
\coordinate (gpbar) at (5.96, 2.56);

\coordinate (p1) at (2.96, 6.24);

\coordinate (z) at (12.08, 2.44);

\draw  (pn)-- (pn1) -- (pn2) -- (pn3);
\draw [domain=0:12.5] plot(\x,{(10.43-0.07*\x)/3.89});
\draw [dashed] (rn1)-- (rn) -- (pn);

\draw [domain=2.5:12.5] plot(\x,{(21.18-1.18*\x)/2.83});
\foreach \point/\label in {pn1/p_{n-1}, pn2/p_{n-2}}  {
    \draw [fill=black] (\point) circle (1pt);
    \node [anchor=north west, inner sep=1pt] at (\point) {$\label$};
}

\foreach \point/\label in {pn3/p_{n-3},  pn1bar/\overline{p}_{n-1}} {
    \draw [fill=black] (\point) circle (1pt);
    \node [anchor=north east, inner sep=1pt] at (\point) {$\label$};
}

\draw[fill = black] (rn1) circle (1pt); 

\draw[fill = black] (rn) circle (1pt); 
\node[anchor=north west, inner sep = 1pt] at (rn) {$r_n$};

\draw[fill = black] (gpbar) circle (1pt);

\draw[fill = black] (pn1p) circle (1pt); 
\node[anchor = south] at (pn1p) {$p_{n-1}'$};

\draw[fill = black] (p1) circle (1pt); 
\node[above] at (p1) {$p_1$};

\draw[fill = black] (z) circle (1pt); 
\node[below] at (z) {$z$};

\draw[fill=black] (pnbar) circle (1pt);
\node[anchor = south west, inner sep = 1pt] at (pnbar) {$\overline p_n$}; 

\draw[fill=black] (pn) circle (1pt);
\node[anchor = south west, inner sep = 1pt] at (pn) {$p_n$};

\node[anchor=south west, inner sep=1pt] at (gpbar) {$g(\overline p_n)$};

\end{tikzpicture}

	\caption{}
	\label{Zcasefig}
\end{figure}

We note that in this case the pole of $\alpha$ must belong to the segment $(\overline p_n,z)$. Indeed, as established above, 
$\alpha(\overline p_n)=r_2\in (p_n,p_1)$. This in turn implies that $p_1'\in (\overline p_n,\alpha(\overline p_n))$.
By Lemma \ref{pocketperspectivities} (ii) applied to the perspectivity $\alpha$ (whose center is $p_1'$) with $a=\overline p_n$ and $b=z$, the pole of $\alpha$ lies on the segment $(\overline p_n,z)$. Moreover, by the properties of perspectivities reviewed in Section \ref{geoprelims}, $\alpha$
maps the segment $[\overline p_n,z]$ onto $\ell_{p_np_1}\backslash (z,r_2)$. In particular, we find   
$w\in (\overline p_n,z]$ such that $\alpha(w)=p_1$. (We can alternatively justify the location of the pole of $\alpha$ and the existence
of $w$ first deducing from $\alpha(\overline p_n)\in (p_n,p_1)$ that $p_1'$ is inside the triangle $(\overline p_n,z,p_1)$,
and then arguing from this that the line $\ell_{p_1p_1'}$ and the line through $p_1'$ parallel to $\ell_{p_np_1}$ both intersect the segment $[\overline p_n,z]$.)

It follows by Claim 3 that $g(w)\leq g(\overline p_n)$ in the order of $\ell$. This inequality together with $\overline p_n <w$  forces $g$ to have a pole on the interval $(\overline p_n,w)$, since $g$ is increasing before and after its pole. This in turn implies, by the convexity of $g$ on $(-\infty,\overline p_n]$, that the points in $(-\infty,p_n']$ must  all be good, since $p_n'$ is good but $\overline p_n$ is not. We are thus done in this case.

Suppose now that $\ell$ and $\ell_{p_{n}p_1}$ are parallel lines. In this case the perspectivity $\alpha$ is affine and increasing. Thus,  there exists  $w\in (\overline p_n,\infty)$ such that $\alpha(w)=p_1$. As argued in the previous paragraph, this implies that $g(w)\leq g(\overline p_n)$, even though $w<\overline p_n$, which entails that $g$ has a pole on $(\overline p_n,w)$. As before, we deduce from the convexity of $g$ on $(-\infty,\overline p_n]$ that the points in $(-\infty,p_n']$ are all good, and we are done.

Suppose finally that $\ell$ and $\ell_{p_{n}p_1}$ intersect at a point $z$ in the ray $(-\infty,\overline p_n]\subset \ell$. By the convexity of $P$, and the fact that $p_n,p_1$ both lie strictly to the left of $\ell$, the point $z$ lies outside of $P$. Hence, $z\in (-\infty,\overline p_{n-1})$. Observe that $\alpha(z)=z\notin (r_2,p_1)$. By Claim 3, $g(z)\in [\overline p_{n-1},g(\overline p_n)]$. Thus, $z$ is a good point. Moreover, since  $z<\overline p_n$ and $g(z)\leq g(\overline p_n)$, the pole of $g$ is not contained in $[z,\overline p_n]$, i.e., the pole of $g$ is either in $(\overline p_n,\infty)$ or $(-\infty, z)$. In the first case, we have argued 
above that all the points in $(-\infty,p_n']$ are good. If the pole of $g$ is in $(-\infty,z)$, then $g$ is concave down on $[z,p_n']$, and since both endpoints of this segment are good points, all the points in the segment $[z,p_n']$ are good. In particular, the points in $[p_{n-1'},p_n']$ are good.

\textbf{Subcase $\overline p_n\in (p_{n-2},p_{n-1}]$}: We can deal with this subcase similarly to how we   dealt with the previous one. 
From Claim 1, we deduce that  $\overline p_{n-1}\in\ (p_{n-3},p_{n-2})$, for if $\overline p_{n-1}\in [p_{n-2},p_{n-1}]$, then
 $p_n',p_{n-1}'$ would both lie on the edge $[p_{n-2},p_{n-1}]$ of $P$; however,  they belong to the interior of $P$.

Let us show that, as in the previous subcase,  $g(\overline p_n)\in (\overline p_{n-1},\overline p_n)$. Let $S=(s_k)_{k=1}^n$ be the polygon resulting from the BLC relative to $P'$ starting at $s_1=\overline p_n$. By Lemma \ref{sameBLC}, this construction has pivots $p_1',p_2',\ldots,p_{n-1}'$. 
Since $P'$ lies to the left of $\ell_{q_{n-1}p_n}$,  and $\overline p_n\in \arc(q_{n-1},p_n)$, $P'$ lies also to the left of $\ell_{\overline p_np_n}$.
By Lemma \ref{pockets} applied to the partition arcs created by the points $\overline p_n,p_n,p_1,\ldots,p_{n-2}$, we conclude that  $s_{k}\in\ [p_{k-2},p_{k-1})$ for $k=2,\ldots,n-1$ and $s_n\in [p_{n-2},\overline p_n)$.  Together with the sequence of pivots, this implies that the perspectivities used  in the construction of $S$ are  $\alpha,\beta_1,\ldots,\beta_{n-2}$. Thus, $g(\overline p_n)=\gamma(s_n)$. Since $s_n\in [p_{n-2},\overline p_n)$ lies strictly to the right of $\ell$, the segment $[s_n,p_n]$ intersects $\ell$ in an interior point of $P$.
Thus,  
\[
g(\overline p_n)=\gamma(s_n)\in (\overline p_{n-1},\overline p_n).
\]
In particular, $\overline p_n$ is  a bad point

A variation on Claim 3 can be proved in this case as well, with minor modifications:

\emph{Claim 4}: If $y\in \ell$ is such that $\alpha(y)\notin (s_2, p_1)$, then $g(y)\in [\overline p_{n-1},g(\overline p_n)]$.

{\it Proof of Claim 4}: Set $y_1=\alpha(y)$ and $y_{k+1}=\beta_{k}(y_{k})$ for $k=1,\ldots,n-2$. Note that $g(y)=\gamma(y_{n-1})$.
From $y_1\notin (s_2, p_1)$ and Lemma \ref{pocketperspectivities} applied to $\beta_1$, we get that $y_2=\beta_1(y_1)\in [p_1,s_3]$. 
 Continuing this argument inductively, we obtain that $y_k\in [p_{k-1},s_{k+1}]$ for  $k=1,\ldots,n-1$, and in particular, $y_{n-1}\in [p_{n-2},s_n]$. Since both $p_{n-2}$ and $s_n$ lie to the right of $\ell$,  $y_{n-1}$ lies to the right of $\ell$.  
It follows that $\gamma(y_{n-1})$ agrees with the intersection of  $[y_{n-1},p_n]$ with $\ell$, and it is thus  in  $\co(P)$.
Hence, $g(y)\in [\overline p_{n-1},\overline p_n]$. Since   $\gamma$ is increasing on the interval $[y_{n-1},s_n]$, we also get that  
\[
g(y)=\gamma(y_{n-1})\leq \gamma(s_n)=g(\overline p_n)
\] 
in the order of $\ell$. Thus,  $g(y)\in [\overline p_{n-1},g(\overline p_n)]$ as desired.

We now repeat the analysis of the cases where the lines $\ell$ and $\ell_{p_np_1}$ are either parallel or intersect at a point $z$ in either
$(\overline p_n,+\infty)\subset \ell$ or $(-\infty,\overline p_{n-1})\subset \ell$. The arguments used before apply in this case   using Claim 4 in place of Claim 3. 
\end{proof}

\begin{lemma}
Let $P'$ be attainable from $P$ by a single pull-in move. Then either $P'$ is degenerate in $P$ or it is attainable from
$P$ by the application of at most two neighbor pull-in moves.
\end{lemma}

\begin{proof}
Suppose that $P'$ is obtained from $P$ by the pull-in of $p_i$ toward $p_j$, for $j\neq i$, and assume that $j\notin \{i-1,i+1\}$, as otherwise the result is obvious. If $p_i'$ belongs to $\co(p_{i-1},p_i,p_{i+1})$, then we can reach 
the final location of $p_i'$ with a pull-in of $p_i$ toward $p_{i-1}$, followed by a pull in of $p_i$ toward $p_{i+1}$. If on the other hand $p_i'\notin \co(p_{i-1},p_i,p_{i+1})$, then  $P'$ is degenerate, as in this case the $(n-1)$-gon $Q=(p_1,\ldots,p_{i-1},p_{i+1},\ldots,p_n)$ interpolates between $P'$ and $P$.
\end{proof}

\begin{theorem}\label{attainableVuD}
Let $P$ be a convex $n$-gon oriented counterclockwise, with $n\geq 4$. Then an $n$-gon $P'$ contained in $P$ is attainable if and only if
it belongs to the set $\V_P\cup \D_P$. If $P'\in \D_P$, then it is attainable in fewer than $5n$ pull-in moves, while if $P'\in \V_P\backslash \D_P$, then it is attainable in at most $2n$ pull-in moves. 
\end{theorem}

\begin{proof}
One direction of the theorem has already been proved: the polygons in $\V_P\cup\D_P$ are attainable from $P$.
Let us prove the other direction. By Theorem \ref{closedundermoves}, the set $\V_P\cup \D_P$ is closed under neighbor pull-in moves. Combined with the previous lemma, this yields that $\V_P\cup \D_P$ is in fact closed under arbitrary pull-in moves. Since $P\in \TT_P\subseteq \V_P$, any polygon attainable from $P$ in finitely many pull-in moves belongs to $\V_P\cup \D_P$. Now let $P'$ be an $n$-gon attainable from $P$. If $P'$ is degenerate in $P$, then $P'\in \V_P\cup \D_P$ and we are done. Suppose that $P'$ is non-degenerate in $P$. By Corollary \ref{Pchattering}, $P'=\lim_k P_k'$, where  each $P_k'$ is attainable from $P$ in finitely many pull-in moves. As established before, $P_k'\in \V_P\cup \D_P$ for all $k$. Since  $\V_P\cup\D_P$ is a closed set (Lemmas \ref{compactD} and \ref{compactV}), we obtain that $P'\in \V_P\cup\D_P$, as desired. 

Finally, let us discuss the estimates on the number of pull-in moves given in the theorem. 
By Theorem \ref{in5nmoves}, the polygons in $\D_P$ are attainable from $P$ in fewer than $5n$ moves. Suppose on the other hand that a polygon $P'$ belongs to $\V_P\backslash \D_P$. By the definition of $\V_P$, there exists $P''\in \TT_P$ (the threshold region)
such that $P'$ is attained from $P''$ in one pull-in move. The polygon $P''$ is necessarily non-degenerate in $P$, so it is attainable from $P$ in at most $2n-1$ pull-in moves, by Theorem \ref{thresholdthm}. Thus, $P'$ is attainable from $P$ in at most $2n$ pull-in moves.
\end{proof}

Let us outline how to effectively check whether an $n$-gon $P'$ contained in  a given $n$-gon $P$ is attainable by a decreasing by path.
\begin{enumerate}
\item
If $P$ is not set-convex, then $P'$ is attainable from $P$ (by Theorem \ref{in5nmoves}).

\item
If $P$ is set-convex, then there is a permutation of the indices $\{1,\ldots,n\}$, such that $P_{\sigma}=(p_{\sigma(i)})_{i=1}^n$
is convex oriented counterclockwise. Then $P'$ is attainable from $P$ if and only if $P'_\sigma=(p'_{\sigma(i)})_{i=1}^n$ is attainable from $P_{\sigma}$. Assume then that $P$ is convex oriented counterclockwise.

\item
Test whether  $P'\in \D_P$ (see remarks after Theorem \ref{degeneracytestthm}). If so, then $P'$ is attainable. 

\item
If $P'\notin \D_P$, test whether $P'\in \V_P\backslash \D_P$ by considering every 
neighbor push-out of a vertex of $P'$ onto $\partial P$ and testing whether the resulting polygons belong to $\TT_P$ using Theorem \ref{thresholdthm}
(the Threshold Theorem). If so, then $P'$ is attainable from $P$, and otherwise, it is not.    
\end{enumerate}


\begin{bibdiv}
\begin{biblist}

\bib{frydman80}{article}{
   author={Frydman, Halina},
   title={The embedding problem for Markov chains with three states},
   journal={Math. Proc. Cambridge Philos. Soc.},
   volume={87},
   date={1980},
   number={2},
   pages={285--294},
}


\bib{frydman2}{article}{
   author={Frydman, Halina},
   title={A structure of the bang-bang representation for $3\times 3$
   embeddable matrices},
   journal={Z. Wahrsch. Verw. Gebiete},
   volume={53},
   date={1980},
   number={3},
   pages={305--316},
}


\bib{frydman83}{article}{
   author={Frydman, Halina},
   title={On a number of Poisson matrices in bang-bang representations for
   $3\times 3$ embeddable matrices},
   journal={J. Multivariate Anal.},
   volume={13},
   date={1983},
   number={3},
   pages={464--472},
}
	
\bib{goodman}{article}{
   author={Goodman, Gerald S.},
   title={An intrinsic time for non-stationary finite Markov chains},
   journal={Z. Wahrscheinlichkeitstheorie und Verw. Gebiete},
   volume={16},
   date={1970},
   pages={165--180},
}

\bib{goodmanCIME}{article}{
   author={Goodman, G. S.},
   title={The embedding problem for stochastic matrices},
   conference={
      title={Stochastic differential equations},
   },
   book={
      series={C.I.M.E. Summer Sch.},
      volume={77},
      publisher={Springer, Heidelberg},
   },
   date={2010},
   pages={231--249},
}



\bib{johansen73}{article}{
   author={Johansen, S.},
   title={A central limit theorem for finite semigroups and its application
   to the imbedding problem for finite state Markov chains},
   journal={Z. Wahrscheinlichkeitstheorie und Verw. Gebiete},
   volume={26},
   date={1973},
   pages={171--190},
}
	
\bib{johansen-ramsey}{article}{
   author={Johansen, S.},
   author={Ramsey, F. L.},
   title={A bang-bang representation for $3\times 3$ embeddable stochastic
   matrices},
   journal={Z. Wahrsch. Verw. Gebiete},
   volume={47},
   date={1979},
   number={1},
   pages={107--118},
}

\bib{kingman-william}{article}{
   author={Kingman, J. F. C.},
   author={Williams, David},
   title={The combinatorial structure of non-homogeneous Markov chains},
   journal={Z. Wahrscheinlichkeitstheorie und Verw. Gebiete},
   volume={26},
   date={1973},
   pages={77--86},
}

\bib{lee}{book}{
   author={Lee, J. M.},
   title={Axiomatic geometry},
   series={Pure and Applied Undergraduate Texts},
   volume={21},
   publisher={American Mathematical Society, Providence, RI},
   date={2013},
   pages={xviii+469},
}
	
\end{biblist}	
\end{bibdiv}
\end{document}